\documentclass[12pt,a4paper,leqno]{article}
\usepackage{latexsym}
\usepackage{amsmath}
\usepackage{amssymb,amscd,amsthm,epic,url}
\usepackage[all]{xy}
\usepackage[dvips]{graphicx}
\usepackage[english,francais]{babel}

\DeclareFontEncoding{OT2}{}{} 
  \newcommand{\textcyr}[1]{%
    {\fontencoding{OT2}\fontfamily{wncyr}\fontseries{m}\fontshape{n}%
     \selectfont #1}}
\newcommand{\sha}{{\mbox{\textcyr{Sh}}}}

\newcommand{\ra}{\rightarrow}

\newcommand{\ol}{\overline}
\newcommand{\ul}{\underline}

\newcommand{\rS}{\mathrm{S}}

\newcommand{\rH}{\mathrm{H}}
\newcommand{\rI}{\mathrm{I}}

\newcommand{\rN}{\mathrm{N}}

\newcommand{\rR}{\mathrm{R}}

\newcommand{\rU}{\mathrm{U}}

\newcommand{\rW}{\mathrm{W}}

\newcommand{\Aut}{\mathrm{Aut}}

\newcommand{\Isom}{\mathrm{Isom}}
\newcommand{\Isomext}{\mathrm{Isomext}}

\newcommand{\img}{\mathrm{Im}}

\newcommand{\sico}{\mathrm{sc}}

\begin{document}
\centerline
{\bf Embeddings of maximal tori in classical groups and }
\smallskip
\centerline
{\bf explicit Brauer--Manin obstruction}
\bigskip
\centerline
{E. Bayer-Fluckiger, T-Y. Lee, R. Parimala}
\bigskip

\date{\today}

 \centerline {\bf   Introduction}
\bigskip
Embeddings of maximal tori into classical groups over global fields of characteristic $\not = 2$ are the subject matter of  several recent papers (see for instance Prasad and
Rapinchuk  [PR 10],  [F 12], [Lee 14], [B 12], [B 13]), with special attention to the Hasse principle. 
In particular, it is shown in [Lee 14] that the Brauer--Manin obstruction is the only one. 

\medskip
The present paper gives necessary and suffi\-cient conditions for the embedding of maximal tori in classical groups.
As in  [PR 10], the embedding problem will be described in terms of embeddings of \'etale algebras with involution into central simple algebras with
involution. Let $(E,\sigma)$ be an \'etale algebra with involution defined over a global field, satisfying certain dimension conditions (cf. \S 1). 
In \S 4, we define a group $\sha(E,\sigma)$ which plays an important role in the embedding problem. 

\medskip
Let $(A,\tau)$ be a central simple algebra with involution defined over the same global field, and assume that  everywhere locally there exists an (oriented) embedding of
$(E,\sigma)$ 
in $(A,\tau)$. Then we define a map $f : \sha(E,\sigma) \to {\bf Z}/2{\bf Z}$ such that $(E,\sigma)$ can be embedded in $(A,\tau)$
globally if and only if $f = 0$ (cf. Theorem  5.5.1.).

\medskip
By [Lee 14] we know that the  Brauer--Manin obstruction is the only one, hence we obtain
an explicit description of this obstruction.

\medskip
In addition to the Hasse principle, one also needs to know when an embedding exists
over local fields. This is done in [Lee 14] in terms of Tits indices, and in \S 3 of the present paper (see also 6.3. and 8.2.) in terms of classical invariants. 
Finally, \S 9 contains some applications and examples. In particular, we recover Theorem A of Prasad and Rapinchuk (see [PR 10], Introduction, page 584).

\medskip
The paper has two appendices. The first one outlines the relationship of the point of view and results of [Lee 14] and those of the present paper, and the
second one contains a new proof of Theorem B of Prasad and Rapinchuk (see [PR 10], Introduction, page 586).

\medskip
We thank Gopal Prasad for his interest in our results, and for encouraging us to include an alternative proof of Theorem B of [PR 10] in our paper.

\bigskip
{\bf \S 1. Definitions, notation and basic facts}
\bigskip

{\bf 1.1. Embeddings of algebras with involution}
\medskip

Let $L$ be a field, and let $A$ be a central simple algebra  over $L$. Let $\tau$ be an involution of $A$, and let $K$ be the fixed field of $\tau$ in $L$. Recall that
$\tau$ is said to be of the {\it first kind} if $K = L$ and of the {\it second kind} if $K \not = L$; in this case, $L$ is a quadratic extension of $K$. Let ${\rm dim}_L(A) = n^2$.
Let $E$ be a commutative \'etale algebra of rank $n$ over $L$, and let $\sigma: E \to E$ be a $K$--linear  involution. such that $\sigma | L = \tau | L$.  Set  $F = \{e \in E | \sigma(e) = e \}$. 
Assume that the following dimension condition holds :

\begin{center}
$\dim_K F =\left\{
         \begin{array}{ll}
             \ \ \  n & \hbox{if $L \not = K$;} \\
             \ [{{n+1} \over 2}] & \hbox{if $L = K$.}
           \end{array}
         \right.$
\end{center}

\medskip

An {\it embedding} of $(E,\sigma)$ in $(A,\tau)$ is by definition an injective homomorphism
$f : E \to A$ such that $\tau (f(e)) = f(\sigma(e))$ for all $e \in E$. It is well--known that embeddings
of maximal tori into classical groups can be described in terms of embeddings of \'etale algebras with involution 
into central simple algebras with involution satisfying the above dimension hypothesis (see for
instance [PR 10], Proposition 2.3). 

\medskip
We say that a separable field extension $E'/L$ is a {\it factor} of $E$ if $E = E' \times E''$ for some \'etale $L$--algebra $E''$. It is
well--known that  $E$ can be embedded in the  algebra $A$ if and only if each of the factors of $E$
splits $A$ : 

\medskip
\noindent
{\bf Proposition 1.1.1.} {\it The \'etale algebra $E$ can be embedded in the  central simple algebra $A$ if and only if for every factor $E'$ of $E$, the
algebra $A \otimes _L E'$ is a matrix algebra over $E'$.}

\medskip 
\noindent
{\bf Proof.} See for instance [PR 10], Proposition 2.6.

\medskip
Let $\epsilon : E \to A$ be an $L$--embedding which may not respect the given involutions. The following properties are well--known :

\medskip
\noindent
{\bf Proposition 1.1.2.} {\it There exists a $\tau$--symmetric element $\alpha  \in A^{\times}$ such that for
$$\theta =
 \tau \circ {\rm Int}(\alpha)
 $$ we have
 $$\epsilon (\sigma(e)) = \theta (\epsilon (e)) \ {\rm for \ all} \   e \in E,$$
in other words $\epsilon : (E,\sigma) \to (A,\theta)$ is an $L$---embedding of algebras with involution.}

\medskip
\noindent
{\bf Proof.} See [K 69], \S2.5. or [PR 10], Proposition 3.1. 

\medskip
Note that $\theta$ and $\tau$ are of the same type (orthogonal, symplectic or unitary), since $\alpha$ is $\tau$--symmetric.

\medskip
For all $a \in F^{\times}$, let $\theta_a : A \to A$ be the involution given by $\theta_a =  \theta  \circ {\rm Int}(\epsilon (a))$. Note that
$\epsilon : (E,\sigma) \to (A,\theta_a)$ is an embedding of algebras with involution.

\medskip
\noindent
{\bf Proposition 1.1.3.} {\it The following conditions are equivalent :

\smallskip 
{\rm (a)} There exists an $L$--embedding  $\iota : (E,\sigma) \to (A,\tau)$ of algebras with involution.

\smallskip 
{\rm (b)} There exists an $a \in F^{\times}$ such that $(A,\theta_a) \simeq (A,\tau)$ as algebras with involution.}

\medskip
\noindent
{\bf Proof.} See [PR 10], Theorem 3.2.

 \medskip
 If $\iota : (E,\sigma) \to (A,\tau)$ is an embedding of algebras with involution, and if $a \in F^{\times}$, $\alpha \in A^{\times}$ are
 such that ${\rm Int}(\alpha) : (A,\theta_a) \to  (A,\tau)$ is an isomorphism of algebras with involution
 satisfying  ${\rm Int}(\alpha) \circ  \epsilon = \iota$, then $(\iota,a,\alpha)$ are
 called {\it parameters} of the embedding.

\medskip
\noindent
{\bf Lemma 1.1.4.} {\it Let $a, b \in F^{\times}$ and let $\alpha \in A^{\times}$. Then we have :

\smallskip
{\rm (i)}  ${\rm Int}(\alpha) : (A,\theta_a) \to (A,\theta_b)$ is
an isomorphism of algebras with involution if and only if there exists  $\lambda \in L^{\times}$ such that $\theta (\alpha) \epsilon (b) \alpha  = \lambda \epsilon (a)$.

\smallskip

 {\rm (ii)} Moreover,  we have  ${\rm Int}(\alpha) \circ \epsilon = \epsilon$ if and only if there exists $y \in E^{\times}$ and  
 $\lambda \in L^{\times}$ such that $\alpha = \epsilon (y)$ and ${\rm N}_{E/F}(y) = \lambda a b^{-1}$.}

 \medskip
 \noindent
 {\bf Proof.} The proof of  {\rm (i)}  follows from a direct computation. Let us prove  {\rm (ii)}. If ${\rm Int}(\alpha) \circ \epsilon = \epsilon$, then we have $\alpha \epsilon (x) \alpha^{-1} = \epsilon (x)$ for all $x \in E$. Since $E$
 is a maximal commutative subalgebra of $A$, this implies that $\alpha \in \epsilon (E)$. Let $\alpha = \epsilon (y)$ for some $y \in E$. Then we have
 $\theta (\epsilon (y)) \epsilon (b) \epsilon (y) = \lambda \epsilon (a)$. This implies that $b \sigma(y) y = \lambda a$, hence ${\rm N}_{E/F}(y) = \lambda a b^{-1}$. The converse is clear.

 \bigskip
In particular, there exists an isomorphism of algebras with involution $(A,\theta_a) \to (A,\theta_b)$ commuting with $\epsilon$ if and only if we have   $a b^{-1} \in L^{\times}  {\rm N}_{E/F}(E^{\times})$.

\medskip
\noindent
{\bf Definition 1.1.5.} We say that $(E,\sigma)$ is {\it split} if there exists an idempotent $e \in E$ such that $e + \sigma(e) = 1$.

\medskip
Equivalently, $(E,\sigma)$ is split if $E \simeq E_1 \times E_2$ with $\sigma(E_1) = E_2$. 

\medskip
\noindent
{\bf Definition 1.1.6.} We say that $(A,\tau)$ is {\it hyperbolic} if there exists an idempotent $a \in A$ such that $a + \tau(a) = 1$.

\medskip
Equivalently, $(A,\tau)$ is hyperbolic if $A \simeq M_r(D)$ for some division algebra $D$, and $\tau$ is induced by a hyperbolic hermitian form 
over $D$ (cf. [KMRT 98],  Chapter II, (6.7) and (6.8)).

\medskip
\noindent
{\bf Proposition 1.1.7.} {\it Suppose that $(E,\sigma)$ is split. Then the following are equivalent :

\smallskip
{\rm (a)} The \'etale algebra with involution $(E,\sigma)$ can be embedded in the central simple algebra with involution $(A,\tau)$.

\smallskip
{\rm (b)} All the factors of $E$ split $A$, and the involution $(A,\tau)$ is hyperbolic.}

\medskip
\noindent
{\bf Proof.} Assume that {\rm (a)} holds. Then by Proposition 1.1.1. all the factors of $E$ split $A$.  Let $\iota : (E,\sigma) \to (A,\tau)$ be an embedding, and
 let $e \in E$ be an idempotent such that $e + \sigma(e) = 1$. Set $a = \iota(e)$. Then
$a \in A$ is an idempotent, and
we have $a + \tau(a) = 1$, hence $(A,\tau)$ is hyperbolic. Conversely, assume that {\rm (b)} holds. Since 
$E$ can be embedded in $A$, by Proposition 1.1.2.  there exists an involution $\theta : A \to A$ such that $(E,\sigma)$ embeds into $(A,\theta)$. Hence
$(A,\theta)$ is hyperbolic, and therefore $(A,\tau) \simeq (A,\theta)$. By Proposition 1.1.3. this implies that $(E,\sigma)$ embeds into $(A,\tau)$, hence
{\rm (a)} holds.

\bigskip

{\bf 1.2. Scaled trace forms}

\medskip Let us keep the notation introduced in 1.1. In particular, 
$(E,\sigma)$ is an \'etale algebra with involution. Let $a \in F^{\times}$, and let us consider the form

$$T_a : E \times E \to L$$ given by $$T_a(x,y) = {\rm Tr}_{E/L}(a x \sigma (y)).$$

Then $T_a$ is a quadratic form if $L = K$, and a hermitian form if $L/K$ is a quadratic extension. For $a = 1$, we use the notation $T = T_1$.

\medskip
\noindent
{\bf Proposition 1.2.1.} {\it Let $a \in F^{\times}$. Then we have $${\rm det}(T_a) = {\rm N}_{E/L}(a) {\rm det}(T),$$
in $K^{\times}/K^{\times 2}$ if $L = K$, and in $K^{\times}/{\rm N}_{L/K}(L^{\times})$ if $L/K$ is a quadratic extension.}

\medskip
\noindent
{\bf Proof.} Let $E^{\sharp}$ be  the $L$--vector space of $\sigma$-semilinear homomorphisms $f : E \to L$ (i.e. $f(\lambda x) = \sigma(\lambda)f(x)$  for all $x \in E$ and
$\lambda \in L$). For any quadratic or hermitian form $b : E \times E \to L$, let us
denote by ${\rm ad}(b) : E \to E^{\sharp}$ the $L$--linear map defined by  ${\rm ad}(b)(x)(y) = b(x,y)$ for all $x,y \in E$. Let $(e_1,\dots,e_n)$ be an $L$--basis
of $E$, and let $(e_1^{\sharp},\dots,e_n^{\sharp})$ be the dual basis. Then ${\rm det}(b)$ is the determinant of ${\rm ad}(b)$ in the bases $(e_1,\dots,e_n)$ and
$(e_1^{\sharp},\dots,e_n^{\sharp})$.

\medskip
Let $m_a : E \to E$ be the multiplication by $a$. By definition, we have ${\rm N}_{E/L}(a) = {\rm det}(m_a)$.  Note that we have ${\rm ad}(T_a) = {\rm ad}(T) \circ m_a$. 
This implies that  ${\rm det}(T_a) = {\rm N}_{E/L}(a) {\rm det}(T)$. 

\medskip
\noindent
{\bf Corollary 1.2.2.} {\it Let $a \in F^{\times}$. Then we have

\medskip
\noindent
{\rm (a)} If $L = K$ and $n$ is even, then ${\rm det}(T_a) = {\rm det}(T)$.

\medskip
\noindent
{\rm (b)} If $L$ is a quadratic extension of $K$, then ${\rm det}(T_a) = {\rm N}_{F/K}(a) {\rm det}(T)$. }

\medskip
\noindent
{\bf Proof.} Let us assume that $L = K$.  By Proposition 1.2.1. we have ${\rm det}(T_a) = {\rm N}_{E/K}(a) {\rm det}(T) \in K^{\times}/K^{\times 2}.$ Since $a \in F^{\times}$,
we have  ${\rm N}_{E/K}(a) =  {\rm N}_{F/K}(a)^2$, which is an element of $K^{\times 2}$, hence  we have ${\rm det}(T_a) =  {\rm det}(T) \in K^{\times}/K^{\times 2}.$ 
This proves {\rm (a)}. 

\medskip
Suppose now that $L$ is a quadratic extension of $K$. Then Proposition 1.2.1. implies that ${\rm det}(T_a) = {\rm N}_{E/L}(a) {\rm det}(T) \in K^{\times}/{\rm N}_{L/K}(L^{\times})$. Since $a \in F^{\times}$, we have ${\rm N}_{E/L}(a) = {\rm N}_{F/K}(a)$, and this implies {\rm (b)}.

\bigskip

{\bf 1.3. The discriminant of an \'etale algebra with involution}

\medskip Recall that $T : E \times E \to L$ is defined by $T(x,y) = {\rm Tr}_{E/L}(x \sigma (y))$.

\medskip
\noindent
{\bf Definition 1.3.1.} Set ${\rm disc}(E,\sigma) = {\rm det}(T)$, considered as an element of $K^{\times}/K^{\times 2}$ if $L = K$, and as an
element of $K^{\times}/{\rm N}_{L/K}(L^{\times})$ if $L/K$ is a quadratic extension. This element is called the {\it discriminant} of the \'etale algebra
with involution $(E,\sigma)$.

\medskip
\noindent
{\bf Lemma 1.3.2.} {\it Suppose that $L = K$, and that $n = 2r$. Then 

\medskip
{\rm (i)} ${\rm disc}(E,\sigma) = (-1)^r {\rm disc}(E)$.

\medskip
{\rm (ii)} For all $a \in F^{\times}$ we have  ${\rm disc}(T_a) = {\rm disc}(E)$.}

\medskip
\noindent 
{\bf Proof.} Let us denote by $T_{E/K} : E \times E \to K$, given by $(x,y) \mapsto {\rm Tr}_{E/K}(xy)$, the usual trace form.
We have ${\rm rank}(E) = 2 {\rm rank}(F) = 2r$. Writing $E = F (\sqrt d)$ for some $d \in F^{\times}$, a computation shows
that ${\rm det} (T) = (-1)^r{\rm det}(T_{E/K})$. By definition, we have 
${\rm disc}(E) = {\rm det}(T_{E/K})$, hence  ${\rm disc}(E,\sigma) = (-1)^r {\rm disc}(E)$.

\smallskip
Since  ${\rm disc}(T) = (-1)^r {\rm det}(T)$, by {\rm (i)} we have
${\rm disc}(T) = {\rm disc}(E)$. By Corollary 1.2.2. we have ${\rm disc}(T_a) =  {\rm disc}(T)$, hence ${\rm disc}(T_a) = {\rm disc}(E)$.

\bigskip

{\bf 1.4. An embedding criterion}

\medskip
Assume that $A = M_n(L)$ and that $\tau$ is an orthogonal or unitary involution. Then $\tau : A \to A$ is given by an $n$--dimensional form $b : V \times V \to L$, which is
quadratic if $L = K$ and hermitian if $L \not = K$. We have an embedding criterion, in terms of the forms introduced in 1.2 :

\medskip
\noindent
{\bf Proposition 1.4.1.} {\it There exists an embedding of algebras with involution $(E,\sigma) \to (A,\tau)$ if and only if
there exists $a \in F^{\times}$ such that $b \simeq T_a$.}

\medskip
\noindent
{\bf Proof.} If $L = K$, then this is well--known (see for instance [PR], 7.1.). The proof is similar in the case when $L \not =  K$. However,
we give a proof for the convenience of the reader.

\medskip Note that $A = {\rm End}(V)$. Since $\tau$ is induced by $b : V \times V \to L$, we have $b(ex,y) = b(x,\tau(e)y)$ for all $e \in {\rm End}(V)$ and all $x,y \in V$. 

\medskip
Suppose first that there exists $a \in F^{\times}$ such that $b \simeq T_a$. Let us identify $V$ to $E$ and $b$ to $T_a$,  and note that sending $e \in E$ to the multiplication
by $e$ gives rise to an  embedding $E \to {\rm End} (E)$. Identifying $b$ to $T_a$, we have, for all $e,x,y \in E$,
$$b(ex,y) = {\rm Tr}_{E/L}(a e x \sigma(y)) =  {\rm Tr}_{E/L}(a x \sigma(\sigma(e)y))  = b(x,\sigma(e)y).$$
Since this holds for all $x \in E$, and that $b(ex,y) = b(x,\tau(e)y)$ for all $x \in E$, we have $\sigma(e) = \tau(e)$ for all $e \in E$. Hence
 the natural embedding of $E$ in $A \simeq {\rm End}(E)$ is an embedding of algebras with involution.

\medskip
Suppose now that there exists an embedding of algebras with involution $\iota : (E,\sigma) \to (A,\tau)$. Then for all $e \in E$, we
have $$b(\iota (e)x,y) = b(x, \tau (\iota (e))y) = b(x, \iota (\sigma (e))y).$$

\medskip Let us show that there exists a hermitian form $h : V \times V \to E$ such that $b(x,y) = {\rm Tr}_{E/L}(h(x,y))$ for all $x,y \in V$.
Let us fix $x,y \in V$, and let us consider the linear form $E \to L$ such that $e \mapsto f(e) = b(\iota(e)x,y)$. Since $E$ is a separable
$L$--algebra, there exists $e' \in E$ such that ${\rm Tr}_{E/L}(e e') = f(e)$ for all $e \in E$. Set $h(x,y) = e'$. Let us check that $h$
is a hermitian form. It is easy to see that $h$ is linear in the first variable, so it remains to check that $\sigma (h(x,y)) = h(y,x)$ for
all $x,y \in V$. We have $${\rm Tr}_{E/L}(e \sigma (h(x,y))) = \sigma \ [ {\rm Tr}_{E/L}(\sigma(e) h(x,y))] = \sigma [b(\iota(\sigma(e))x, y))] =$$
$$ = \sigma [b(x,\iota(e)y)] = \sigma[\sigma[b(\iota(e)y,x)]] = b(\iota(e)y,x) = {\rm Tr}_{E/L}(e h(y,x)).$$ Since this holds for all $e \in E$, we have
$h(y,x) = \sigma (h(x,y))$, as claimed. Therefore $h : V \times V \to E$ is a one dimensional hermitian form. Let us identify the 1--dimensional
$E$--vector space $V$ with $E$. Then there
exists $a \in F^{\times}$ such that $h(x,y) = a x \sigma(y)$. Hence we have $b \simeq T_a$, and this completes the proof
of the Proposition. 

\medskip Note that if $(A,\theta)$ is the involution induced by $T$ and if $a \in F^{\times}$, then $T_a$ induces the involution $(A,\theta_a)$. 

\medskip

{\bf  1.5. Invariants of central simple algebras with involution}

\medskip
If $(A,\tau)$ is of orthogonal type and $n$ is even, we denote by ${\rm disc}(A,\tau)$ its {\it discriminant} (cf. [KMRT 98], Chap II. (7.2), and by $C(A,\tau)$ its {\it Clifford algebra} (cf. [KMRT 98], Chap II.  (8.7)).  We denote by $Z(A,\tau)$ the center of the algebra $C(A,\tau)$. Then 
$Z(A,\tau)$ is a quadratic \'etale algebra over $K$. If $(A,\tau)$ is unitary, then we denote by ${D}(A,\tau)$ its {\it discriminant algebra} (cf. [KMRT 98], Chap II, (10.28)).
The {\it  signature} of $(A,\tau)$ is defined in [KMRT 98], Chap II. (11.10) and (11.25).

\medskip
If moreover $A \simeq M_n(L)$, then $\tau$ is induced by a symmetric, skew--symmetric or hermitian form, according as $\tau$ is of orthogonal,
symplectic or unitary type. In this case, we have some additional invariants, such as the Hasse invariant (in the orthogonal case), as well
as the determinant  (in the unitary case). In particular, if $\tau$ is unitary and induced by a hermitian form $h$ over $L/K$, then we set ${\rm det}(A,
\tau) =
{\rm det}(h) \in K^{\times}/{\rm N}_{L/K}(L^{\times})$. Let us write $L = K(\sqrt \delta)$, and let us denote by ${\rm Br}(K)$ the Brauer group of $K$. Then we have $D(A,\tau) = ({\rm disc}(h),\delta) \in {\rm Br}(K)$, where ${\rm disc}(h) = (-1)^{n(n-1)/2} {\rm det}(h)$
(cf. [KMRT 98], Chap II. (10.35)). If $\tau$ is orthogonal and induced by a quadratic form $q$ over $K$, then we denote by $w(q) \in {\rm Br_2}(K)$ its
Hasse invariant. 

\medskip
Let $d \in F^{\times}$ be such that
$E = F(\sqrt d)$. The following result is due to Brusamarello, Chuard--Koulmann and Morales (cf.  [BCM 03], Theorem 4.3.) :

\medskip
\noindent
{\bf Lemma 1.5.1.} {\it Let $(A,\theta)$ be an orthogonal involution. Assume that $n$ is even, and let $a \in F^{\times}$. Then we have $w(T_a) = w(T) + {\rm cor}_{F/K}(a,d).$}

\medskip
\noindent
{\bf Lemma 1.5.2.} {\it Let $(A,\theta)$ be a unitary involution, and let $a \in F^{\times}$. Then $D(A,\theta_a) = D(A,\theta) + {\rm cor}_{F/K} (a,d)$.}

\medskip
\noindent
{\bf Proof.} By  [KMRT 98], Chap. II, (10.36), we have $D(A,\theta_a) = D(A,\theta) + ({\rm N}_{F/K} (a),L/K)$. We have $ ({\rm N}_{F/K} (a),L/K) =
{\rm cor}_{F/K}(a,E/F) = {\rm cor}_{F/K}(a,d)$, hence the lemma is proved.

\medskip
{\bf 1.6. Some necessary embedding conditions}

\medskip
The existence of an embedding of algebras with involution $(E,\sigma) \to (A,\tau)$ implies the following relationship between the
discriminants of $E$ and $(A,\tau)$ :

\medskip
\noindent
{\bf Proposition 1.6.1.} {\it  Suppose that the degree of $A$ is even, and that $(A,\tau)$ is of the orthogonal type. If there exists an embedding of algebras with
involution $(E,\sigma) \to (A,\tau)$, then we have  ${\rm disc} (E) = {\rm disc} (A,\tau) \in K^{\times}/K^{\times 2}.$}

\medskip
\noindent
{\bf Proof.} Let $M$ be the function field of the Severi--Brauer variety of the algebra $A$. Then we have $A \otimes_K M \simeq M_n(M)$, and the
involution $\tau$ is induced by a quadratic form $q$ over $M$. By Proposition  1.4.1. and Lemma 1.3.2.  {\rm (ii)} (see also [B 12], Lemma 1.4.1.) we have 
${\rm disc} (E \otimes_K M) = {\rm disc} (q) \in M^{\times}/M^{\times 2}.$ Since the natural map  $K^{\times}/K^{\times 2} \to  M^{\times}/M^{\times 2}$ 
is injective, we have ${\rm disc} (E) = {\rm disc} (A,\tau) \in K^{\times}/K^{\times 2}$.

\medskip
\noindent
{\bf Proposition  1.6.2.} {\it Suppose that  $A \simeq M_n(L)$, and that $(A,\tau)$ is of the unitary type. If there exists an embedding of algebras with
involution $(E,\sigma) \to (A,\tau)$, then we have  ${\rm det}(A,\tau) {\rm disc}(E,\sigma) ^{-1} \in {\rm N}_{F/K}(F^{\times})  
 {\rm N}_{L/K}(L^{\times}).$}

\medskip
\noindent
{\bf Proof.} Since $A \simeq M_n(L)$, the involution $\tau$ is induced by a hermitian form $h$. By Proposition 1.4.1. there exists $a \in F^{\times}$
such that $h \simeq T_a$. By Corollary 1.2.2. {\rm (b)} we have  ${\rm det}(T_a) =  {\rm N}_{F/K} (a) {\rm det}(T) $.
Recall that ${\rm  disc}(E,\sigma)$ is by definition equal to  ${\rm det}(T) \in K^{\times}/{\rm N}_{L/K}(L^{\times})$. We have
${\rm det}(A,\tau) = {\rm det}(h) = {\rm det}(T_a)$. This implies that ${\rm det}(A,\tau) = {\rm N}_{F/K} (a) {\rm disc}{(E,\sigma)} \in K^{\times}/{\rm N}_{L/K}(L^{\times}) $,
hence we have
$ {\rm det}(A,\tau) {\rm disc}(E,\sigma) ^{-1} \in {\rm N}_{F/K}(F^{\times})  
 {\rm N}_{L/K}(L^{\times}).$ 
 
 \bigskip

 {\bf  \S 2. Orientation}

\medskip
In order to treat the {\it non--split orthogonal} case, we need an additional tool, namely the notion of {\it orientation}. Assume that $(A,\tau)$ is an
orthogonal involution, and that the degree of $A$ is even. Let us set ${\rm deg}(A) = 2r$. 

\medskip
We have seen that the existence of an embedding of algebras with involution $(E,\sigma) \to (A,\tau)$ implies that ${\rm disc}(E) = {\rm disc}(A,\tau) \in K^{\times}/K^{\times 2}$ (see Proposition 1.6.1.). Therefore the discriminant algebra of $E$ (see below) is isomorphic to the $K$--algebra $Z(A,\tau)$. However, such an
isomorphism is not unique. This leads to the notions of {\it orientation}, and of  {\it oriented embedding}, needed for the analysis of
the Hasse principle (see 6.1.). 

\medskip

{\bf  2.1. Discriminant algebra}

\medskip
We have $E \simeq  F[X]/(X^2 - d)$ for some $d \in F^{\times}$. Let 
us consider the $F$--linear involution  $\sigma' : F[X]/(X^2 - d) \to F[X]/(X^2 - d)$ determined by $\sigma'(X) = - X$. Then we have an isomorphism of
algebras with involution $(E,\sigma) \simeq (F[X]/(X^2 - d), \sigma')$. Let $x$ be the image of $X$ in $E$, and note that we have $\sigma (x) = -x$. Let $\Delta (E)$ be the discriminant algebra of $E$ (cf. [KMRT 98], Chapter V, \S 18,
p. 290). 

\medskip
\noindent
{\bf Lemma 2.1.1.} {\it We have an isomorphism of $K$--algebras $$\Delta (E) \simeq K[Y]/(Y^2 -  (-1)^r {\rm N}_{E/K}(x)).$$}
\noindent
{\bf Proof.}  Recall that  $T_{E/K} : E \times E \to K$, defined by $T_{E/K}(e,f) = {\rm Tr}_{E/K}(ef)$, is the trace form of $E$. Then by [KMRT 98],  Proposition (18.2) we have $$\Delta (E) \simeq K[Y]/(Y^2 - {\rm det}(T_{E/K})).$$ Note that ${\rm Tr}_{E/K} = { \rm Tr}_{F/K} \circ {\rm Tr}_{E/F}$, and that the trace form $T_{E/F} : E \times E \to F$,
defined by $T_{E/F}(e,f) = {\rm Tr}_{E/F}(ef)$, is isomorphic to $<2,2d>$. Further, we have $d = -  {\rm N}_{E/F}(x)$ and hence  ${\rm N}_{F/K}(d) = (-1)^r  {\rm N}_{E/K}(x)$. Therefore we have  ${\rm det}(T_{E/K}) =  (-1)^r {\rm N}_{E/K}(x) \in K^{\times}/K^{\times 2}$, and
this concludes the proof of the lemma.

\medskip
Let us denote by $y$ the image of $Y$ in $\Delta (E)$. The elements $x$ and $y$ will be fixed in the sequel. Let $\rho : \Delta (E) \to \Delta(E)$ be the
automorphism of $\Delta (E)$ induced by $\sigma$. Note that we have $\rho (y) = (-1)^r$, and that hence $\rho$ is the identity if $r$ is even and the
non--trivial automorphism of the quadratic algebra $\Delta(E)$ if $r$ is odd. 

\medskip
{\bf 2.2. Generalized Pfaffian}

\medskip For any central simple algebra $A$ over $K$  of degree $2r$ with an orthogonal involution $\theta$,
let us denote by ${\rm Skew}(A,\theta)$ the set $\{a \in A \ | \ \theta (a) = -a \} $ of skew elements
of $A$ with respect to the involution $\tau$.  Recall that $C(A,\theta)$ is the Clifford algebra 
of $(A,\theta)$, and that  $Z(A,\theta)$ is the center of $C(A,\theta)$. Recall that $Z(A,\theta)$
is a quadratic \'etale algebra over $K$. 
Let us denote by $\gamma$ the
non--trivial automorphism  of $Z(A,\theta)$ over $K$.

\medskip
The {\it generalized Pfaffian} (cf. [KMRT 98], Chapter II, \S 8) of $(A,\theta)$ is a homogeneous polynomial map of degree $r$,
denoted by $$\pi_{\theta}  : {\rm Skew}(A,\theta) \to Z(A,\theta)$$  such that for all $a \in {\rm Skew}(A,\theta)$, we have  
$\gamma(\pi_{\theta} (a))  = - \pi_{\theta}(a)$, and
$\pi_{\theta}(a)^2 = (-1)^r {\rm Nrd}(a)$; for all $x \in A$ and $a \in {\rm Skew}(A,\theta)$,  we have $(\pi_{\theta}(x a \theta (x)) = {\rm Nrd}_A(x) \pi_{\theta}(a)$
(cf. [KMRT 98], Proposition (8.24)).

\medskip
{\bf 2.3. Orientation}

\medskip For any orthogonal involution $(A,\tau)$, an isomorphism of $K$--algebras $$\Delta(E) \to Z(A,\tau)$$ will be called an {\it orientation}.

\medskip Let us assume that the \'etale algebra $E$ can be embedded in the central simple algebra $A$, and
let us fix an embedding $\epsilon : E \to A$. By Proposition 1.1.2.  there exists an involution $\theta : A \to A$ of orthogonal type such
that $\epsilon : (E,\sigma) \to (A,\theta)$ is an embedding of algebras with
involution. 

\medskip Let us fix such an involution $(A,\theta)$. We now define an orientation $u : \Delta(E) \to Z(A,\theta)$ that will be fixed in
the sequel. Fix a generalized Pfaffian map $\pi_{\theta}  : {\rm Skew}(A,\theta) \to Z(A,\theta)$ as above. Recall that 
$E \simeq F[X]/(X^2 - d)$, that $\Delta(E) \simeq K[Y]/(Y^2 - (-1)^r {\rm N}_{E/K}(x)),$
and that we have fixed the images $x$ of $X$ in $E$ and $y$ of $Y$ in 
$\Delta(E)$.
Let $$u : \Delta(E) \to Z(A,\theta)$$ be defined by $$y \mapsto \pi_{\theta}( \epsilon (x)).$$

\noindent
{\bf Lemma 2.3.1.} {\it The map $u$ is an isomorphism of $K$--algebras.}

\medskip
\noindent
{\bf Proof.} 
We have $\gamma (\epsilon (x) ) = - \epsilon (x)$. Further, $(\pi_{\theta} (\epsilon (x))^2 = (-1)^r {\rm Nrd}_A(\epsilon (x)) = (-1)^r {\rm N}_{E/K}(x) = y^2.$
This implies that $u$ is an isomorphism of $K$--algebras.

\medskip

{\bf 2.4. Similitudes}

\medskip
Let $\alpha \in A^{\times}$. Following [KMRT], Definition (12.14), page 158, we say that $\alpha$ is a  {\it similitude} of $(A,\tau)$ if $\alpha \tau (\alpha) \in K^{\times}$.  For a similitude $\alpha \in A^{\times}$, the scalar $\alpha \tau (\alpha)$ is called the {\it multiplier} of the similitude. We say that
$\alpha$ is a {\it proper}
similitude if ${\rm Nrd}(\alpha) =  (\alpha \tau (\alpha))^r$; otherwise, ${\alpha}$ is called an {\it improper}  similitude. Note that ${\alpha}$
is a  similitude if and only if  ${\rm Int}(\alpha) : (A,\tau) \to (A,\tau)$ is an isomorphism of algebras with involution.
If $A$ is split, then $(A,\tau)$ admits improper similitudes (indeed, any reflection is an improper similitude).

\medskip

Any isomorphism of algebras with involution ${\rm Int}(\alpha) : (A,\tau) \to (A,\tau')$ induces an isomorphism of the Clifford algebras $C(A,\tau) \to C(A,\tau')$.
Let us denote by $$c(\alpha) : Z(A,\tau) \to Z(A,\tau')$$ the restriction of this isomorphism to the centers of the Clifford algebras. The following property
will be important in the sequel.

\medskip
\noindent
{\bf Lemma 2.4.1.} {\it Let $(A,\tau)$ be an orthogonal involution, and let $\alpha \in A^{\times}$ be a similitude. Then $\alpha$ is a proper similitude
if and only if $c(\alpha)$ is the identity}.

\medskip
\noindent
{\bf Proof.} See for instance [KMRT 98], Proposition (13.2), page 173.

\medskip

{\bf 2.5. Compatible orientations}

\medskip
Recall that $\epsilon : E \to A$ is an embedding of algebras, that $\theta : A \to A$ is an orthogonal involution such that
$\epsilon : (E,\sigma) \to (A,\theta)$ is an embedding of algebras with involution, and that we are fixing an orientation
$u : \Delta (E) \to Z(A,\theta)$. We now define a notion of {\it compatibility} of orientations.

\medskip
\noindent
{\bf Lemma 2.5.1.} {\it Let $(A,\tau)$ be a central simple algebra with an orthogonal involution, and let $\iota : (E,\sigma) \to (A,\tau)$ be a embedding
of algebras with involution. Let $\alpha \in A^{\times}$ be such that ${\rm Int}(\alpha) :(A,\tau) \to (A,\tau)$ is an automorphism of algebras with involution, and
${\rm Int}(\alpha) \circ \iota = \iota$. Then

\medskip
{\rm (a) } There exists  $x \in E^{\times}$ such that $\alpha = \iota (x)$, and ${\rm N}_{E/F}(x) \in K^{\times}$. 

\smallskip
{\rm (b)} The map $c(\alpha)$ is the identity.}

\medskip
\noindent
{\bf Proof.} Since ${\rm Int}(\alpha) \circ \iota = \iota$, the restriction of ${\rm Int}(\alpha)$ to $\iota (E)$ is the identity. Note that $\iota (E)$ is a maximal
commutative subalgebra of $A$. Hence we have $\alpha = \iota (x)$ for some $x \in E^{\times}$. As ${\rm Int}(\alpha): (A,\tau) \to (A,\tau)$ is an
automorphism of algebras with involution, we have $\alpha \tau (\alpha) = \lambda$ for some $\lambda \in K^{\times}$. Hence we have  $(\iota x) \tau (\iota x) = 
\lambda$. Since $\iota : (E,\sigma) \to (A,\tau)$ is an embedding of algebras with involution, we have $\iota (x \sigma(x) ) = \lambda$. This completes the
proof of {\rm (a)}. 

\medskip
Let us prove (b). By part (a), we have $\alpha \tau \alpha = \iota (x \sigma (x))  = \iota (\lambda) = \lambda$. This implies that
$\alpha$ is a similitude. Moreover, we have  ${\rm Nrd}(\alpha) = {\rm N}_{E/K}(x) = {\rm N}_{F/K}(\lambda)
= \lambda^r$. Hence $\alpha$ is a proper similitude, and by lemma 1.7.1. this implies that $c(\alpha)$ is the identity.

\medskip

\noindent
{\bf Definition 2.5.2.} Let $\theta' : A \to A$ be an orthogonal involution such that
$\epsilon: (E,\sigma) \to (A,\theta')$ is an embedding of algebras with involution, and let $u' : \Delta (E) \to Z(A,\theta')$ be an orientation. We say
that the orientations $u$ and $u'$ are {\it compatible} if for every isomorphism of algebras with involution ${\rm Int}(\alpha) :(A,\theta) \to (A,\theta')$
such that ${\rm Int}(\alpha) \circ \epsilon = \epsilon$, we have $u ' = c(\alpha) \circ u$. 

\medskip
Recall that for all $a \in F^{\times}$, we define an involution $\theta_a : A \to A$ by $\theta_a = \theta \circ {\rm Int}(\epsilon (a))$. Note that the embedding $\epsilon : (E,\sigma) \to (A,\theta)$ induces
an embedding of algebras with involution $\epsilon : (E,\sigma) \to (A,\theta_a)$. Our next aim is to define an orientation of $(A,\theta_a)$ compatible with
the orientation $u$ of $(A,\theta)$. 
Let $K_s$ be a separable closure of $K$, and set $A_s = A \otimes _K K_s$.

\medskip
\noindent
{\bf Proposition 2.5.3.} {\it Let $a \in F^{\times}$. Then there exists a unique isomorphism $\phi_a : Z(A,\theta) \to Z(A,\theta_a)$ such that
for all $\alpha \in A_s^{\times}$ which gives an isomorphism of algebras with involution ${\rm Int}(\alpha) : (A_s,\theta) \to (A_s,\theta_a)$  with ${\rm Int}(\alpha) \circ \epsilon = \epsilon$, we have $c(\alpha)  = \phi_a$.}

\medskip
\noindent
{\bf Proof.} 
Let $d \in K^{\times}$ represent the square class of ${\rm disc}(A, \theta)$, and let us write $Z(A,\theta) = K \oplus Kz$ with $z^2 = d$. Note that
$d$ also represents the square class of ${\rm disc}(A, \theta_a)$, since $a \in F^{\times}$.
Let us write 
$Z(A,\theta_a) = K \oplus Kz_a$ with $z_a^2 = d$.

\medskip
 Let $b \in (E \otimes_K K_s)^{\times}$ be such that $b \sigma(b) = a^{-1}$. 
 Then ${\rm Int}(\epsilon(b)) : (A_s,\theta) \to (A_s, \theta_a)$ is an isomorphism of algebras with involution commuting with $\epsilon$,  and it induces an
isomorphism of the Clifford algebras  $C(A_s,\theta) \to  C(A_s, \theta_a)$.

 \medskip
 We have $A_s  = M_{2r}(K_s)$, and $\theta :  A_s \to A_s$ is induced by a
 quadratic form $q : V \times V \to K_s$.  Let $(e_1,\dots,e_{2r})$ be an orthogonal basis for $q$. Since $Z(A,\theta) = K \oplus K(e_1\dots e_{2r})$,
 we have $z = \mu (e_1\dots e_{2r})$ for some $\mu \in K_s^{\times}$. Let us replace $e_1$ by $\mu^{-1}e_1$. Then we have $z = e_1\dots e_{2r}$. 
 
 \medskip
 Set $q = \epsilon(b)^t q_a \epsilon (b)$. Since $a^{-1} = b \sigma(b)$ and $a$ is $\theta$--symmetric, the involution induced by $q_a$ is $\theta_a$. Let us consider the isometry $\epsilon(b) : (V,q) \to (V,q_a)$. Then $\epsilon(b)$ induces a map $c(\epsilon(b)) : C(V,q) \to C(V,q_a)$ which sends $e_1\dots e_{2r}$ to
 $(\epsilon(b)e_1)\dots (\epsilon(b)e_{2r})$. Therefore we have $(\epsilon(b)e_1)\dots (\epsilon(b)e_{2r})^2 = q_a(\epsilon(b)e_1)\dots q_a(\epsilon(b)e_{2r})
 = q(e_1) \dots q(e_{2r}) = (e_1 \dots e_{2r})^2 = d$. This implies that
 $\epsilon(b) (e_1) \dots \epsilon(b) (e_{2r}) = \pm z_a$ and $c(\epsilon (b))(z) = \pm z_a$. Hence the restriction of the map $c(\epsilon (b))$ to $Z(A_s,\theta)$ is defined over $K$.

 \medskip
 Set $\phi_a = c( \epsilon (b))$, and note that $\phi_a : Z(A,\theta) \to Z(A,\theta_a)$ is an isomorphism. 
 
 \medskip
 Let us show that $\phi_a$ is independent of the choice of $b$. Let $b' \in A_s$ such that $b' \sigma (b') = a$. Then we have $c({\rm Int} (\epsilon (b'))) = c({\rm Int}
 (\epsilon (b))). $ We have an isomorphism of algebras with involution ${\rm Int}( \epsilon (b^{-1}b')) : (A,\theta) \to (A,\theta)$ satisfying 
 ${\rm Int}( \epsilon (b^{-1}b')) \circ \epsilon = \epsilon$. Hence by Lemma 2.4.1.  the map $$c({\rm Int}( \epsilon (b^{-1}b'))) : Z(A,\theta) \to Z(A,\theta)$$ is
 the identity. Therefore $c( \epsilon (b))  = c( \epsilon (b'))$, hence $c( \epsilon (b)) $ is independent of the choice of $b$. 
 
 \medskip
Let $\alpha \in A_s^{\times}$ be such that ${\rm Int}(\alpha) : (A_s,\theta) \to (A_s,\theta_a)$  is an isomorphism of algebras with involution with ${\rm Int}(\alpha) \circ \epsilon = \epsilon$. Then by Lemma 2.4.1. there exists $x \in (E \otimes_KK_s)^{\times}$ such that  $\alpha = \epsilon (x)$. This implies that $c({\rm Int} (\epsilon (x))) = c(\epsilon (b)) =
\phi_a$. Hence $c(\alpha) = \phi_a$, as required. This also shows the uniqueness of $\phi_a$, and  completes the proof of the proposition.

\bigskip

Recall that we have fixed an isomorphism $u : \Delta(E) \to Z(A,\theta)$. For all $a \in F^{\times}$, let us define an orientation by 
$u_a = \phi_a \circ u : \Delta(E) \to Z(A,\theta_a)$. Then $u_a$ is compatible with $u$. Note that $\phi_1$  is the identity, hence $u_1 = u$. 

\medskip
For all $a \in F^{\times}$, let us identify $\Delta(E)$ with $Z(A,\theta_a)$ via the orientation $u_a$. This endows the Clifford algebra $C(A,\theta_a)$
with a structure of $\Delta (E)$--algebra. We have the following

\medskip
\noindent
{\bf Lemma 2.5.4.} {\it For all  $a \in F^{\times}$ we have 
$$C(A,\theta_a)  = C(A,\theta) + {\rm res}_{\Delta (E)/K}{\rm cor}_{F/K}(a,d)$$ in ${\rm Br}(\Delta (E))$.}

\medskip
\noindent
{\bf Proof.} This follows from   [BCM 03], Proposition 5.3. 

\medskip
{\bf 2.6. Oriented embeddings}

\medskip
Recall that the existence of  an embedding of algebras with involution $(E,\sigma) \to (A,\tau)$ is equivalent with the existence of  an element $a \in F^{\times}$ such that
the algebras with involution  $(A,\theta_{a})$ and $(A,\tau)$ are isomorphic. We need the stronger notion of oriented embedding, defined
as follows :

\medskip
\noindent
{\bf Definition 2.6.1.} Let $(A,\tau)$ be an orthogonal involution, and let $\nu : \Delta (E) \to Z(A,\tau)$ be an orientation. An embedding $ \iota : (E, \sigma) \to (A, \tau)$ is called an {\it oriented
embedding} with respect to $\nu$ if there exist $a \in F^{\times}$ and $\alpha \in A^{\times}$  satisfying the following conditions :

\medskip

{\rm (a)} ${\rm Int}(\alpha) : (A,\theta_{a}) \to (A,\tau)$   is an isomorphism of algebras with involution such that ${\rm Int}(\alpha)  \circ  \epsilon = \iota$.

\medskip

{\rm (b)} The induced automorphism  $c(\alpha) :  Z(A,\theta_{a})  \to Z(A,\tau)$ satisfies $$c(\alpha) \circ   u_{a} = \nu.$$

\medskip
We say that there exists an oriented embedding of algebras with involution with respect to $\nu$ if there exists $ (\iota, a, \alpha)$
as above. The elements $ (\iota, a, \alpha,\nu)$ are called {\it parameters} of the oriented embedding.

\medskip

 {\bf 2.7. Changing the orientation -- improper similitudes}
 
 \medskip
Let  $\nu :   \Delta (E) \to Z(A,\tau)$ be an orientation. We have
 
\medskip
\noindent
{\bf Proposition 2.7.1.} {\it Suppose that $(A,\tau)$ admits an improper similitude. Assume  that there exists an embedding of algebras with involution $(E,\sigma) \to (A,\tau)$. Then there
exists an oriented embedding $(E,\sigma) \to (A,\tau)$ with respect to $\nu$.  Moreover, if $(\iota,a,\alpha)$ are parameters of an embedding of $(E,\sigma)$ in  $(A,\tau)$, then 
there
exist $\iota'$ and $\beta$ such that $(\iota',a,\beta,\nu)$ are parameters of an oriented embedding.}

\medskip
\noindent
{\bf Proof.} If $c(\alpha) \circ u_a = \nu$, then $({\rm Int}(\alpha) \circ \epsilon, a, \alpha)$ are parameters of an oriented embedding
$(E,\sigma) \to (A,\tau)$. Suppose that $c(\alpha) \circ u_a \not = \nu$. Let $\gamma \in A^{\times}$ be an improper similitude. Then $c(\gamma)$ is not
the identity, and hence we have $c(\gamma \alpha) \circ u_a = \nu$. Set  $\beta = \gamma \alpha$.
Then $({\rm Int}(\beta) \circ \epsilon ,a,\beta)$ are parameters of an oriented embedding, as claimed.

\medskip
\noindent
{\bf Lemma 2.7.2.} {\it Let Suppose that $K$ is  a local field or the field of real numbers, and let $(A,\tau)$ be an orthogonal involution. Assume that  if $A$ is non--split, then $ {\rm disc}(A,\tau) \not = 1
\in K^{\times}/K^{\times 2}$. Then $(A,\tau)$ admits improper similitudes.}

\medskip
\noindent
{\bf Proof.} If $A$ is split, then any reflection is an improper similitude. Suppose now that $A$ is not split.
Then we have $A \simeq M_r(H)$, where $H$ is a quaternion division algebra. Let $Z = Z(A,\tau)$. Set
$D = {\rm disc}(A,\tau)$, and note that $Z \simeq K(\sqrt {D})$. Then $Z$ is a quadratic extension of $K$, since $D  \not  \in K^{\times 2}$. Hence
$H$ is split by $Z$. The involution $\tau$ is induced by an $r$--dimensional hermitian form $h$ over $H$.
If $r > 3$, then  the hermitian form $h$ is isotropic (see [T 61], Theorem 3,  if $K$ is a local field, and [Sch 85], Theorem 10.3.7. if $K$ is
the field of real numbers). Therefore $h \simeq h' \oplus h''$, where $h'$ and $h''$ are
hermitian forms over $H$ with ${\rm dim}(h') \le 3$ and $h''$ hyperbolic. Let $r' = {\rm dim}(h')$, and let $B = M_{r'}(H)$. Let $\tau'$ be the
involution of $B$ induced by $h'$, and note that ${\rm disc}(B,\tau') = {\rm disc}(A,\tau) = D$. Since $H$ is split by $Z$, 
we have 
$H= (\lambda,D) \in {\rm Br}(K)$ for some $\lambda \in K^{\times}$. 

\medskip
We claim that  $\lambda$ is a multiplier of a similitude of $(B,\tau')$. Indeed,
since $r' \le 3$, we may apply the criterion of [PT 04], Theorem 4. Let $\gamma (B,\tau') \in {\rm Br}(K)$ such that $\gamma_Z = C(B',\tau')$ in ${\rm Br}(Z)$
(cf. [PT 04], Theorem 2). Then by [PT 04], Theorem 4, the element $\lambda$ is the multiplier of a similitude of $(B,\tau')$ if and only if  $\lambda . \gamma = 0$
in $H^3(K)/(K^{\times}.A)$. If $K$ is a local field, then $H^3(K) = 0$, hence the condition is fulfilled. Assume that $K$ is the field of real numbers. Then
either $\gamma = 0$, or $\gamma = H$ in ${\rm Br}(K)$. Since $A$ is non-split, we have $A = H$ in ${\rm Br}(K)$. Therefore we have $\lambda . \gamma = 0$
n $H^3(K)/(K^{\times}.A)$ in both cases.

\medskip

Therefore by [PT 04], Theorem 4,  the element $\lambda$ is the multiplier of a similitude of
$(B,\tau')$, therefore also of the hermitian form $h'$. The hermitian form $h''$ is hyperbolic, therefore $h''$ has a similitude of multiplier $\lambda$. Thus
the hermitian form $h$ has a similitude of multiplier $\lambda$ as well, and hence $(A,\tau)$ has a similitude of multiplier $\lambda$. By [PT 04], Theorem 1,
using the fact that $A = H= (\lambda,D) \in {\rm Br_2}(K)$, we see that $\lambda$ is the multiplier of an improper similitude.

\medskip
\noindent
{\bf Corollary 2.7.3.} {\it Suppose that there exists an embedding of algebras with involution $(E,\sigma) \to (A,\tau)$, and that one of the
following holds :

\medskip
{\rm (i)} $A$ is split.

\smallskip
{\rm (ii)} $K$ is a local field, or the field of real numbers, and ${\rm disc}(A,\tau) \not = 1$ in $K^{\times}/K^{\times 2}$.

\medskip
Then there
exists an oriented embedding $(E,\sigma) \to (A,\tau)$ with respect to $\nu$. Moreover, if  $(\iota,a,\alpha)$ are parameters of an embedding of $(E,\sigma)$ in  $(A,\tau)$ and if  then there
exist $\iota'$ and $\beta$ such that $(\iota',a,\beta,\nu)$ are parameters of an oriented embedding.}

\medskip
\noindent
{\bf Proof.} In both cases, $(A,\tau)$ admits an improper similitude. If $A$ is split, then any reflection in ${\rm U}(A,\tau)$ is an improper similitude.
If $K$ is local or the field of real numbers, then Lemma 2.7.2. implies that $(A,\tau)$ has an improper similitude. 
Hence the Corollary follows from Proposition 2.7.1.

\medskip

{\bf 2.8. Changing the orientation -- r odd}

\medskip 
Recall that 
$E \simeq F[X]/(X^2 - d)$, that $\Delta(E) \simeq K[Y]/(Y^2 - (-1)^r {\rm N}_{E/K}(x)),$
and that we have fixed the images $x$ of $X$ in $E$ and $y$ of $Y$ in 
$\Delta(E)$. Recall that $\rho : \Delta(E) \to \Delta(E)$ is the 
automorphism of $\Delta(E)$ induced by $\sigma : E \to E$, and  that $\rho$ is the identity if $r$ is even, and the non--trivial automorphism of $\Delta(E)$ over $K$ if $r$ is odd.

\medskip
Recall also that $u : \Delta(E) \to Z(A,\theta)$ is defined by $y \mapsto \pi_{\theta}( \epsilon (x)).$

\medskip
\noindent
{\bf Lemma 2.8.1.} {\it Let ${\rm Int}(\gamma) : (A,\theta) \to (A,\theta)$ be an isomorphism of algebras with involution safisfying ${\rm Int}(\gamma) \circ \epsilon \circ \sigma = \epsilon$. Then we have $c(\gamma) \circ u \circ \rho = u$.}

\medskip
\noindent
{\bf Proof.} It suffices to prove that this is true over a separable closure. Therefore we may assume that $A = M_{2r}(K)$ and that $\theta : A \to A$ is the transposition.
We have $\gamma \theta(\gamma) = \gamma \gamma^t = \lambda$ for some $\lambda \in K^{\times}$. Recall that  ${\rm Nrd}(\gamma) = \eta \lambda^r$, where
$\eta = 1$ if $\gamma$ is a proper similitude, and $\eta = -1$ if $\gamma$ is an improper similitude.
We have $\epsilon (x) = {\rm Int}(\gamma) \circ \epsilon \circ \sigma (x) = \gamma \epsilon (\sigma (x)) \gamma^{-1} = \lambda^{-1} \gamma \epsilon (\sigma (x)) \gamma^t $.

\medskip
On the other hand, we have $\pi_t(\lambda^{-1} \gamma (\epsilon (\sigma (x)) \gamma ^t) = \lambda^{-r} {\rm Nrd}(\gamma) \pi_t(\epsilon (\sigma(x)) = \eta \pi_t  (- \epsilon (x)) = (-1)^r \eta \pi_t(\epsilon (x))$. Hence we have $(-1)^r \eta \pi_t(\epsilon (x)) = \pi_t(\epsilon (x))$, thus $\eta = (-1)^r $. This implies that $\gamma$ is
a proper similitude if $r$ is even, and an  improper similitude if $r$ is odd.  By Lemma 2.4.1.  this implies that $c(\gamma)$ is the identity if $r$ is even, and the non--trivial automorphism of $Z(A,\theta)$ if $r$ is odd. 
Therefore we have $c(\gamma) \circ u \circ \rho (y) = u(y)$, and hence $c(\gamma) \circ u  \circ \rho = u$. 

\medskip
\noindent
{\bf Proposition 2.8.2.} {\it Let $a, b \in F^{\times}$, and let ${\rm Int}(\alpha) : (A,\theta_a) \to (A,\tau)$ and ${\rm Int}(\beta) : (A,\theta_b) \to (A,\tau)$ be  isomorphisms of algebras with involution
such that ${\rm Int}(\alpha) \circ \epsilon \circ \sigma = {\rm Int}(\beta) \circ  \epsilon$. Then we have  $c(\alpha) \circ u_a  \circ \rho =
c(\beta) \circ u_b$.}

\medskip
\noindent
{\bf Proof.} Let $K_s$ be a separable closure of $K$, and let $\gamma_a, \gamma_b \in K_s^{\times}$ be such that ${\rm Int}(\gamma_a)  : (A,\theta) \to (A,\theta_a)$ 
and ${\rm Int}(\gamma_b)  : (A,\theta) \to (A,\theta_b)$ are isomorphisms of algebras with involution commuting with $\epsilon$. Then we have
$u_a = c(\gamma_a) \circ u$ and $u_b = c(\gamma_b) \circ u$. We have ${\rm Int}(\gamma_b^{-1} \beta^{-1} \alpha \gamma_a) \circ \epsilon \circ \sigma =
{\rm Int} (\gamma_b^{-1} \beta^{-1} \alpha) \circ ({\rm Int}(\gamma_a))  \circ \epsilon  \circ \sigma = {\rm Int}(\gamma_b^{-1} \beta^{-1}) \circ {\rm Int}(\alpha) \circ \epsilon \circ \sigma =
{\rm Int}(\gamma_b^{-1} \beta^{-1}) \circ {\rm Int}(\beta) \circ \epsilon = {\rm Int}(\gamma_b^{-1}) \circ \epsilon = \epsilon$.  By Lemma 2.8.1.  this implies
that $c(\gamma_b^{-1} \beta^{-1} \alpha \gamma_a) \circ u \circ \rho = u$, hence we have $c(\alpha) \circ u_a \circ  \rho = c(\beta) \circ u_b$. 

\medskip
Let $\nu :   \Delta (E) \to Z(A,\tau)$ be an orientation.

\medskip
\noindent
{\bf Corollary 2.8.3.} {\it Suppose that $r$ is odd, and that there exists an embedding of algebras with involution $(E,\sigma) \to (A,\tau)$.
Then there
exists an oriented embedding $(E,\sigma) \to (A,\tau)$ with respect to $\nu$. Moreover, if  $(\iota,a,\alpha)$ are parameters of an embedding of $(E,\sigma)$ in  $(A,\tau)$, then there
exist $\iota'$, $b$ and $\beta$ such that $(\iota',b,\beta,\nu)$ are parameters of an oriented embedding.}

\medskip
\noindent
{\bf Proof.} 
Let 
$(\iota,a,\alpha)$ be parameters of an embedding of $(E,\sigma)$ in  $(A,\tau)$. If $c(\alpha) \circ u_a = \nu$, then $(\iota,a,\alpha,\nu)$ are parameters of an 
oriented embedding with respect to $\nu$. Otherwise, we have  $c(\alpha) \circ u_a  \circ \rho = \nu$. Set $\iota' = \iota \circ \sigma$. Then there exist $b \in F^{\times}$ and  $\beta \in A^{\times}$ such that $\iota' = {\rm Int}(\beta) \circ \epsilon$. By
Proposition 2.8.2. we have $c(\beta) \circ u_b =  c(\alpha) \circ u_a  \circ \rho =   \nu$, and hence 
$(\iota',b,\beta,\nu)$ are parameters of an oriented embedding.

\bigskip
{\bf \S 3. Local conditions}

\medskip
The aim of this section is to give necessary and sufficient conditions for an embedding of $(E,\sigma)$ in $(A,\tau)$ to exist when $K$
is a local field of characteristic $\not =2$ or the field of real numbers. This is done in [Lee14] in terms of Tits indices - however, the results of [Lee 14] are not
used here. We assume that all the factors of $E$ split $A$. Hence there exists an embedding of algebras $\epsilon : E \to A$, and an involution $\theta$
of $A$ of the same type as $\tau$ such that $\epsilon : (E,\sigma) \to (A,\theta)$ is an embedding of algebras with involution (cf. Proposition 1.1.2).

\medskip
Let $E =  E_s \times E_n$ where $E_s$ and $E_n$ are \'etale $K$--algebras stable under $\sigma$ with $(E_s,\sigma)$ split and of maximal rank for
this property. Let $2 \rho$ be the rank of $E_s$.

\bigskip
{\bf 3.1. Orthogonal involutions -- the even dimensional case} 

\medskip
Assume that  $(A,\tau)$ is an orthogonal involution.

\medskip
\noindent
{\bf Proposition 3.1.1.} {\it Assume that $n$ is even, and that $K$ is a local field. Then there exists an embedding of algebras with involution of $(E,\sigma)$ into $(A,\tau)$ if
and only if  one of the following conditions holds :

\smallskip
{\rm (i)} $(E,\sigma)$ is split and $(A,\tau)$ is hyperbolic.

\smallskip
{\rm (ii)} $(E,\sigma)$ is not split, and ${\rm disc}(A,\tau) = {\rm disc}(E) \in K^{\times}/K^{\times 2}$.}

\medskip
\noindent
{\bf Proof.} 
{\rm (i)} Suppose that $(E,\sigma)$ is split. Then   $(E,\sigma)$ embeds into $(A,\tau)$
if and only if  $(A,\tau)$ is hyperbolic (cf. Proposition 1.1.7.). 

\medskip
{\rm (ii)}
Suppose that $(E,\sigma)$ is not split. By Proposition 1.6.1. if $(E,\sigma)$ can be embedded  into $(A,\tau)$, then
we have ${\rm disc}(A,\tau) = {\rm disc}(E) \in K^{\times}/K^{\times 2}$. 
Conversely, 
assume that ${\rm disc}(A,\tau) = {\rm disc}(E) \in K^{\times}/K^{\times 2}$. 

\medskip Suppose first that $A$ is split, in other words that $A \simeq M_n(K)$. Then $\tau$ is induced by an $n$--dimensional quadratic form $q$ over $K$, and we have
${\rm disc}(q) = {\rm disc}(A,\tau) = {\rm disc}(E) \in K^{\times}/K^{\times 2}$. By [B 12], Proposition 2.2.1. there
exists $a \in F^{\times}$ such that $w(T_a) = w(q) \in {\rm Br_2}(K)$.  We have ${\rm disc}(T_a) = {\rm disc}(E) \in K^{\times}/K^{\times 2}$ (cf. [B 12], Lemma 1.3.2.)
Therefore the quadratic forms $q$ and $T_a$ have the same dimension, discriminant and Hasse invariant,
hence they are isomorphic. Thus $(E,\sigma)$ embeds  into $(A,\tau)$ (cf. Proposition 1.4.1.). 

\medskip Suppose now that $A$ is not split. Since  $K$ is a local field, we have  $A = M_r(H)$
with $H$ a quaternion division algebra. By Proposition 1.6.1.  we have
${\rm disc}(A,\theta) = {\rm disc}(E) \in K^{\times}/K^{\times 2}$. Therefore ${\rm disc}(A,\tau) = {\rm disc}(A,\theta)$. 
By  [T 61], Theorem 3, this implies that $(A,\tau) \simeq (A,\theta)$. Therefore $(E,\sigma)$ embeds into $(A,\tau)$.
This completes the proof of the Proposition. 

\medskip
\noindent
{\bf Proposition 3.1.2.} {\it Suppose that $K = {\bf R}$ and that $n$ is even. Then there exists an embedding of $(E,\sigma)$ in $(A,\tau)$ if and only
one of the following conditions hold :

\smallskip  {\rm (i)}  $A \simeq M_n({\bf R})$, the involution $\tau$ is induced by the quadratic form $q$, and 
the signature of $q$ is of the shape $(2r + \rho, 2s + \rho)$ for some non--negative integers $r$ and $s$.

\smallskip  {\rm (ii)} $A \simeq M_r(H)$, where $H$ is a quaternion division algebra.}

\medskip
\noindent
{\bf Proof.} If $A \simeq M_n({\bf R})$, then the result follows from  [B 12], Proposition 2.3.2. Assume that $A \simeq M_r(H)$. Then by [Sch 85], Theorem 10.3.7.  we have $(A,\tau) \simeq (A,\theta)$. Therefore $(E,\sigma)$ can be embedded in  $(A,\tau)$.

\bigskip
{\bf 3.2. Orthogonal involutions -- the odd dimensional case} 

\medskip
Assume that  $(A,\tau)$ is an orthogonal involution, and  that $n$ is odd. 
Then we have $A \simeq M_n(K)$, and $\tau$ is induced by an $n$--dimensional quadratic form $q$. We have 
$E = E' \times K$, where $E'$ is a rank $n-1$ \'etale $K$--algebra invariant by $\sigma$. If $n = 1$, then $E = K$ and $\sigma$ is 
the identity. In this case, it is clear that there exists an embedding of $(E,\sigma)$ into $(A,\tau)$. Suppose that $n \ge 3$, and let $F' = (E')^{\sigma}$ be the subalgebra of
$E'$ composed of the elements fixed by the restriction of $\sigma$ to $E'$. We have the following :

\medskip
\noindent
{\bf Proposition 3.2.1.} {\it Assume that $n$ is odd and $n \ge 3$, and that $K$ is a local field. Then there exists an embedding of $(E,\sigma)$ in $(A,\tau)$ if
and only if one of the following holds :

\smallskip
{\rm (i)} $(E',\sigma)$ is split, and $q \simeq q' \oplus q''$ with ${\rm dim}(q') = n-1$ and $q'$ hyperbolic.

\smallskip
{\rm (ii)} $(E',\sigma)$ is not split.}

\medskip
\noindent
{\bf Proof.} Suppose that  $(E,\sigma)$ embeds into  $(A,\tau)$. Then by Proposition 1.4.1. there exists $a \in F^{\times}$ such
that $q \simeq T_a$. We have $F = F' \times K$, and $a = (a',a'')$ with $a' \in (F')^{\times}$ and $a'' \in K^{\times}$. Note that we have  $T_a \simeq T_{a'} \oplus T_{a''}$.
Let $A' = M_{n-1}(K)$, and let $\tau'$ be the involution of $A'$ induced by $T_{a'}$. Then $(E',\sigma)$ embeds into
$(A',\tau')$. If $(E',\sigma)$ is split, then by Proposition 3.1.1.  {\rm (i)}, the quadratic form $T_{a'}$ is hyperbolic. Set $q' = T_{a'}$ and $q'' = T_{a''}$; then  we have $q \simeq q' \oplus q''$ with ${\rm dim}(q') = n-1$ and $q'$ hyperbolic, as claimed.

\medskip
Conversely, suppose that if $(E',\sigma)$ is split, then $q \simeq q' \oplus q''$ with ${\rm dim}(q') = n-1$ and $q'$ hyperbolic, and let us prove that
embeds $(E,\sigma)$ into $(A,\tau)$. Let us show that we have $q \simeq T_{a'} \oplus T_{a''}$ with $a' \in (F')^{\times}$ and $a'' \in K^{\times}$. 
Suppose first that $(E',\sigma)$ is split. Then we have $q \simeq q' \oplus q''$ with ${\rm dim}(q') = n-1$ and $q'$ hyperbolic. Let
$A' = M_{n-1}(K)$ and let $\tau'$ be the involution induced by $q'$. Then  Proposition 3.1.1.  {\rm (i)} implies that  $(E',\sigma)$ embeds into
$(A',\tau')$. Therefore by Proposition 1.4.1. there exists $a' \in (F')^{\times}$ such that $q' \simeq T_{a'}$. Note that as ${\rm dim}(q'') = 1$, there exists
$a'' \in K^{\times}$ such that $q'' \simeq T_{a''}$, hence the statement is proved in this case. Suppose now that $(E',\sigma)$ is not split, and set
$a'' = (-1)^{n-1 \over 2} {\rm det}(q) {\rm disc}(E')$ $
 \in K^{\times}/K^{\times 2}$. Since $K$ is a local field, there exist quadratic forms $q'$ and $q''$ with $q \simeq q' \oplus q''$,  ${\rm dim}(q') = n-1$, ${\rm dim}(q'') = 1$, and
 ${\rm det}(q'') = a''$. This is clear if $n \ge 5$, since a non--degenerate  quadratic form of dimension $\ge 5$ over a local field represents all non--zero elements.
 Assume that $n = 3$. Then $a'' = - {\rm det}(q) {\rm disc}(E')$.  Since $(E',\sigma)$ is not split, we have ${\rm disc}(E') \not \in K^{\times 2}$. The quadratic
 form $q \oplus < {\rm det}(q) {\rm disc}(E')>$ has dimension 4 and non--square discriminant, hence it is isotropic (see for instance [Sch 85], Theorem 6.4.2. page 217). Hence $q$
 represents $a'' = - {\rm det}(q) {\rm disc}(E')$, as claimed. 
 This implies that ${\rm disc}(q') = {\rm disc}(E')$. 
 Set $A' = M_{n-1}(K)$, and let $\tau'$ be  the involution of $A'$ induced by $q'$. Then Proposition 3.1.1.  {\rm (ii)} implies that 
 $(E',\sigma)$ embeds into $(A',\tau')$. By Proposition 1.4.1. we have $q' \simeq T_{a'}$ for some $a' \in (F')^{\times}$. 
Note that as ${\rm dim}(q'') = 1$, we have $q'' \simeq T_{a''}$. 

\medskip
Therefore in both cases we have $q \simeq T_{a'} \oplus T_{a''}$ with $a' \in (F')^{\times}$ and $a'' \in K^{\times}$. Set $a = (a',a'')$. Then $q \simeq q_a$,
and by Proposition 1.4.1. there exists an embedding of $(E,\sigma)$ in $(A,\tau)$. This completes the proof of the Proposition.

\medskip
\noindent
{\bf Proposition 3.2.2.} {\it Suppose that  $K = {\bf R}$ and that $n$ is odd. Then there exists an embedding of $(E,\sigma)$ in $(A,\tau)$ if and only
if the signature of $q$ is of the shape $(r + \rho,  s + \rho)$ for some non--negative integers $r$ and $s$.}

\medskip
\noindent
{\bf Proof.} We have $E =  E' \times {\bf R}$, where $E'$ is a rank $n-1$ \'etale $\bf R$--algebra invariant by $\sigma$. Let $F'$ be the subalgebra
of $E'$ of the elements fixed by $\sigma$. Assume that 
there exists an embedding of $(E,\sigma)$ in $(A,\tau)$. Then by Proposition 1.4.1.  there exists $a \in F^{\times}$ such
that $q \simeq T_a$. We have $F = F' \times {\bf R}$, and $a = (a',a'')$ with $a' \in (F')^{\times}$ and $a'' \in {\bf R}^{\times}$. 
Let $A' = M_{n-1}(K)$, and let $\tau'$ be the involution of $A'$ induced by $T_{a'}$. Then by Proposition 1.4.1. there exists an embedding of $(E',\sigma)$ in
$(A',\tau')$. By Proposition 3.1.2. this implies that the signature of $T_{a'}$ is of the shape $(2r' + \rho, 2s'+ \rho)$ for some non--negative integers $r'$ and $s'$. 
Therefore the signature of $q$ is  $(2r'  + 1 + \rho, 2s'+ \rho)$ or $(2r' + \rho, 2s'+ 1 +  \rho)$.

\medskip
Conversely, assume that the signature of $q$ is of the shape $(r + \rho,  s + \rho)$ for some non--negative integers $r$ and $s$. Then $r + s + 2 \rho = n$,
hence one of  $r$ or $s$ is odd and the other is even. Let $q'$ be a quadratic form of signature $(r - 1 + \rho,  s + \rho)$ if $r$ is odd, and $(r + \rho,  s  - 1 + \rho)$ if
$s$ is odd. Then the dimension of $q'$ is even, and hence by Proposition 2.1.3. there exists an embedding of $(E',\sigma)$ in $(A',\tau')$, where
$A' = M_{n-1}(K)$, and where $\tau'$ is the involution of $A'$ induced by $q'$. Therefore by Proposition 1.4.1. there exists $a' \in (F')^{\times}$ such that
$q' \simeq T_{a'}$. Set $a'' = 1$ if $r$ is odd and $a'' = -1$ if $s$ is odd, and let $a = (a',a'')$. Then $q$ and $T_a$ have the same signature, hence
they are isomorphic. By Proposition 3.1.2.. this implies that there exists an embedding of $(E,\sigma)$ in $(A,\tau)$.

\bigskip
{\bf 3.3. The symplectic case}

\medskip Assume that $(A,\tau)$ is a symplectic involution. If $A \simeq M_r(D)$ for some quaternion division algebra $D$, let $h$ be a
hermitian form with respect to the canonical involution of $D$ which induces $\tau$. The signature of $h$ is defined as in [Sch 85], 10.1.8. {\rm (i)}. 
Let $E = E_s \times E_n$, where $E_s$ and $E_n$ are stable under $\sigma$ such that $(E_s,\sigma)$ is split and $(E_n,\sigma)$ is non-split.
Let $4 \rho$ be the rank of $E_s$. 

\medskip
\noindent
{\bf Theorem 3.3.1.} {\it Suppose that $K$ is a local field, or $K = {\bf R}$.  Then  $(E,\sigma)$ can be embedded in  $(A,\tau)$ if and only if
one of the following holds :

\medskip
{\rm (i)} $K$ is a local field, or $A$ is split.

\medskip
{\rm (ii)} $K = {\bf R}$, $A$ is non--split, and  
${\rm sign}(h)$ 
is of the shape $(s+ \rho, s' + \rho)$, where $s$ and $s'$ are non--negative integers.}

\medskip
\noindent
{\bf Proof.}{\rm (i)}  If $A$ is split, then the involutions are  given by  skew--symmetric matrices with
coefficients in $K$. All non--degenerate skew--symmetric matrices of the same dimension are isomorphic. Hence the algebras with involution $(A,\tau)$ and $(A,\theta)$ are isomorphic,  therefore  $(E,\sigma)$ can be embedded in  $(A,\tau)$.
Suppose  that $A = M_r(D)$, where $D$ is the unique quaternion division  algebra over $K$, and that $K$ is a local field. 
By [Sch 85], 10.1.7. the algebras with involution $(A,\tau)$ and $(A,\theta)$ are isomorphic. Therefore $(E,\sigma)$ can be embedded in  $(A,\tau)$. 

\medskip {\rm (ii)} 
Suppose that $K$ is the field of real numbers and that $A = M_r(D)$, where $D$ is the unique quaternion division  algebra over $K$.
Since all the factors of $E$ split $A$,  the \'etale algebra
$E$ is isomorphic to the direct product of $r$ copies of ${\bf C}$. Hence we have $E_s = ({\bf C} \times {\bf C})^{\rho}$, and $\sigma$ acts on each copy of
${\bf C} \times {\bf C}$ by exchanging the two factors, and we have $E_n = {\bf C}^{r - 2 \rho}$, and $\sigma$ acts on each copy of $\bf C$ by complex
conjugation. Set $F_s = E_s^{\sigma}$ and $F_n = E_n^{\sigma}$. Then we have $F_s \simeq {\bf C}^{\rho}$ and $F_n \simeq {\bf R}^{r - 2 \rho}$. 

\medskip 

Let us consider the following hermitian forms with respect to the canonical involution of $D$ : let $h_1$ be the $2 \rho$--dimensional hyperbolic form,
and let $h_2$ be the $r - 2 \rho$--dimensional unit form. Let $h_0$ be the orthogonal sum of $h_1$ and $h_2$, and let $\theta : A \to A$ be the involution
induced by $h_0$. 

\medskip
Let us denote by $x \mapsto \overline x$ the canonical involution of $D$, and let $\epsilon : E_s \times E_n \to A$ be the map defined by 
$$\alpha (x_1,y_1, \dots, x_{\rho},y_{\rho},z_1,\dots,x_{r-2 \rho}) = {\rm diag}(x_1,\overline y_1, \dots, x_{\rho}, \overline y_{\rho},z_1,\dots,z_{r - 2 \rho}).$$
It is easy to check that $\epsilon$ is an embedding of algebras with involution $(E,\sigma) \to (A,\theta)$.

\medskip If $c_1,\dots,c_{r - 2 \rho} \in  {\bf R}$, let us denote by $h_c$ the $r - 2 \rho$--dimensional diagonal hermitian form  $<c_1,\dots,c_{r - 2 \rho}>$. Let
$H_c$ be the orthogonal sum of the $2 \rho$--dimensional hyperbolic hermitian form with $h_c$. 
For  $a = (b_1,\dots,b_{\rho}) \times (c_1,\dots,c_{r - 2 \rho}) \in {F_s}^{\times} \times {F_n}^{\times}$, we see that the involution $\theta_a : A \to A$ is
induced by the hermitian form $H_c$. Note that the signature of $H_c$ is equal to $(\rho + s, \rho + s')$, where $s$ is the number of positive $c_i$'s, and
$s'$ the number of negative $c_i$'s. By Proposition 1.1.3. there exists an embedding of algebras with involution $(E,\sigma) \to (A,\tau)$ if and only
if there exists $a \in F^{\times}$ such that $(A,\theta_a) \simeq (A,\tau)$. By [Sch 85], 10.1.7. and 10.1.8. {\rm (i)}, this is equivalent with the existence
of $a \in F^{\times}$ such that ${\rm sign}(h_a) = {\rm sign}(h)$. Therefore there exists an embedding of algebras with involution $(E,\sigma) \to (A,\tau)$
if and only if ${\rm sign}(h)$ is of the shape $(s+ \rho, s' + \rho)$, where $s$ and $s'$ are non--negative integers.

\bigskip
{\bf 3.4. The unitary case}

\medskip 
Assume that $L$ is a quadratic extension of $K$, and suppose that $(A,\tau)$ is an $L/K$ unitary involution. 
Assume that $A = M_n(L)$, and that $\tau$ is induced by an $n$--dimensional hermitian form $h$ over $L/K$ (note that
when $K$ is a local field or $K = {\bf R}$, then this hypothesis is always fulfilled).

\medskip
\noindent
{\bf Proposition 3.4.1.} { {\it Suppose that $K$ is a local field. Then there exists an embedding of algebras with involution of $(E,\sigma)$ into $(A,\tau)$ if
and only if  one of the following conditions holds :

\smallskip
{\rm (i)} $(E,\sigma)$ is split and $(A,\tau)$ is hyperbolic.

\smallskip
{\rm (ii)} $(E,\sigma)$ is not split, and ${\rm det}(A,\tau) {\rm disc}(E,\sigma)^{-1} \in {\rm N}_{F/K}(F^{\times}) 
 {\rm N}_{L/K}(L^{\times}).$}
 
 \medskip
 \noindent
 {\bf Proof.} {\rm (i)} follows from Proposition 1.1.7. Suppose that $(E,\sigma)$ is not split, and that $(E,\sigma)$ embeds into $(A,\tau)$. Then by Proposition 1.6.2. we have 
  $${\rm det}(A,\tau) {\rm disc}(E,\sigma)^{-1} \in {\rm N}_{F/K}(F^{\times}) 
 {\rm N}_{L/K}(L^{\times}).$$  Conversely, assume that  ${\rm det}(A,\tau) {\rm disc}(E,\sigma)^{-1}  \in {\rm N}_{F/K}(F^{\times}) 
 {\rm N}_{L/K}(L^{\times}).$ Then there exists $a \in F^{\times}$ such that ${\rm det}(A,\tau)= {\rm N}_{F/K}(a)  {\rm disc} (E,\sigma) \in 
 K^{\times}/{\rm N}_{L/K} (L^{\times}).$ We have ${\rm det}(T_a) = {\rm det}(A,\tau)$, hence  $h$ and $T_a$ have
 equal dimension and determinant. Since $K$ is a local field, this implies that  $T_a \simeq h$. By Proposition 
 1.4.1.  there exists an embedding of $(E,\sigma)$ into $(A,\tau)$. This completes the proof of the proposition.
 
 \medskip
 Recall that  $E =  E_s \times E_n$ where $E_s$ and $E_n$ are \'etale $K$--algebras stable under $\sigma$ with $(E_s,\sigma)$ split and of maximal rank for
this property, and that we denote by $2 \rho$ the rank of $E_s$.

 \medskip
 \noindent
 {\bf Proposition 3.4.2.} {\it Suppose that $K = {\bf R}$. Then there exists an embedding of algebras with involution of $(E,\sigma)$ into $(A,\tau)$ if and
 only if the signature of $h$ is of the shape $(r + \rho,s+\rho)$ for some non--negative integers $r$ and $s$.}
 
 \medskip
 \noindent
 {\bf Proof.} Indeed, since $L/K$ is a quadratic field extension, we have $L = {\bf C}$.  Therefore $E$ is isomorphic to the direct product of $n$ copies of  ${\bf C}$. 
 Let us denote by $\sigma_0 : {\bf C} \to {\bf C}$ the complex conjugation, and by $\sigma_1 : {\bf C} \times {\bf C} \to {\bf C} \times {\bf C}$ the map defined by
 $\sigma_1(a,b) = (\sigma_0(b), \sigma_0(a))$. Then we have $E_s = ({\bf C} \times {\bf C})^{\rho}$, and the restriction of $\sigma : E_s \to E_s$ 
 to each copy of  ${\bf C} \times {\bf C}$  is equal to $\sigma_1$; we have $E_n = {\bf C}^{n - 2\rho}$, and the restriction of $\sigma : E_n \to E_n$ to each
 copy of ${\bf C}$ is equal to $\sigma_0$. 
  
 \medskip

 Set $F_s = E_s^{\sigma}$ and $F_n = E_n^{\sigma}$. Note that $F = F_s \times F_n$, and that $F_s = {\bf C}^{\rho}$ and $F_n = {\bf R}^{n-\rho}$. 
 Let $a = (a_s,a_n) \in F_s^{\times} \times F_n^{\times}$.  Then the restriction of $T_a : E \times E \to K$ to $E_s$ is hyperbolic with signature $(\rho,\rho)$, and its restriction to $E_n$ has
 signature $(r_a,s_a)$, where $r_a$ (respectively $s_a$) is the number of positive (respectively negative) coefficients of $a_n \in {\bf R}^{n-2\rho}$. Hence
 the signature of $T_a$ is $(\rho + r_a,\rho + s_a)$. By Proposition 1.4.1. there exists an embedding of $(E,\sigma)$ into $(A,\tau)$ if and only if
 $h \simeq T_a$ for some $a \in F^{\times}$. Hence $(E,\sigma)$ can be embedded into $(A,\tau)$ if and only if the signature of $h$ is of the shape
 $(\rho +r,\rho + s)$ for some non-negative integers $r$ and $s$.

\bigskip
{\bf  \S 4. The Tate--Shafarevich group}

\medskip
We keep the notation of the previous sections, and suppose that $K$ is a global field. Recall that either $L = K$, or $L$ is a quadratic
extension of $K$. The aim of this section is to define a group that measures the failure of the Hasse principle. 

\medskip

Let us denote by
$\Omega_K$ the set of places of $K$. For
all $v \in \Omega_K$, we denote by $K_v$ the completion of $K$ at $v$.  For all $K$--algebras 
$B$, set $B^v = B \otimes_K K_v$. 

\medskip
The commutative \'etale algebra $E$ is by definition a product of separable field extensions of $L$. Let us write
$E = E_1 \times \dots \times E_m$, with $\sigma(E_i) = E_i$ for all $i = 1,\dots, m$, and such that $E_i$ is either
a  field stable by $\sigma$  or a product of  two fields exchanged by $\sigma$.  Recall that $F = E^{\sigma}$. 

\medskip Set $I = \{1,\dots ,m \}$. 
We have $F = F_1 \times \dots \times F_m$, where $F_i$ is the fixed field of $\sigma$ in $E_i$ for all $i \in I$. Note that
either $E_i = F_i = K$,  $E_i =   F_i \times F_i$ or $E_i$ is a quadratic field extension of $F_i$. For all
$i \in I$, let $d_i \in F_i^{\times}$ such that $E_i = F_i (\sqrt {d_i})$ if $E_i/F_i$ is a quadratic extension, and $d_i = 1$ otherwise. Set
$d = (d_1,\dots,d_m)$. 

\medskip
If $i \in I$ is such that $E_i$ is a quadratic extension of $F_i$, let $\Sigma_i$ be the set of places $v \in \Omega_K$ such that all the places of $F_i$ over $v$ split in $E_i$.
If $E_i =  F_i \times F_i$ or if $E_i = K$,  set $\Sigma_i = \Omega_K$.  

\medskip
If $L \not = K$, let $\Sigma(L/K)$ be the set of places of $K$ that split in $L$. If $L = K$, then we set $\Sigma(L/K) = \emptyset$.

\bigskip
Given an m-tuple $x = (x_1,..., x_m)\in({\bf Z}  /2 {\bf Z})^{m}$, set $$I_0 = I_0(x) = \{ i \ | \ x_i = 0 \},$$ $$I_1 = I_1(x) = \{ i \ | \ x_i = 1 \}.$$ Note that
$(I_0,I_1)$ is a partition of $I$.
Let $S'$ be the set
$$ S' = \{(x_1,...,\ x_m)\in ({\bf Z}  /2 {\bf Z})^{m} \ | \Sigma(L/K) \cup 
(\underset{i\in I_0}{\cap}\Sigma_i)\cup(\underset{j\in I_1}{\cap}\Sigma_j)=\Omega_K \} ,$$ and set $$S = S' \cup(0,...,0)\cup(1,...,1).$$

We define an equivalence relation on $S$ by
\begin{center}
$(x_1,..., x_m)\sim(x'_1,..., x'_m)$ if $(x_1,..., x_m)+(x'_1,...,
x'_m)=(1,..., 1)$. \end{center} 

Let us denote by $\sha = \sha(E,\sigma)$ the set of equivalence classes of $S$ under the above equivalence relation. 
\medskip

For all $x \in S$, we denote by $x$ its class in $\sha$, and by $(I_0(x), I_1(x))$ the corresponding partition of $I$. 
Let us denote by $P'$ the set of non--trivial partitions $(I_0,I_1)$ of $I$ such that $\Sigma (L/K) \cup (\underset{i\in I_0}{\cap}\Sigma_i)\cup(\underset{j\in I_1}{\cap}\Sigma_j)=\Omega_K$,
and set $P = P' \cup \{  (I, \emptyset) \} \cup \{ \emptyset , I)  \}$. 
Let us define an equivalence relation on $P$ by $(I_0,I_1) \sim (I_1,I_0)$. Sending $x$ to $(I_0(x),I_1(x))$ induces a bijection between $\sha$ and
the set of equivalence classes of $P$ under this equivalence relation.

\medskip

Componentwise addition gives a group structure on the set of equivalence classes of $({\bf Z}  /2 {\bf Z})^{m}$.  Let us denote this group by $(C_m,+)$. We have

\medskip
\noindent
{\bf Lemma 4.1.1.} {\it  The set $\sha$ is a subgroup of $C_m$.}

\medskip
\noindent
{\bf Proof.} It is clear that the class of $(0,\dots,0)$ is the neutral element, and that every element is its own opposite, so we only need to check that the sum of two elements
of $\sha$ is again in $\sha$. If $J$ is a subset of $I$, set $\Omega(J) =  \underset{i\in J} \cap \Sigma_i$. 
As we have seen above, the set $\sha$ is in bijection with the set of equivalence classes of partitions $P/ \sim$. Moreover, the transport of structure induces
$$(I_0,I_1) + (I_0',I_1') = ((I_0 \cap I_0')  \cup (I_1 \cap I_1') , (I_0 \cap I_1') \cup( I_1  \cap I_0')).$$ Let us show that this is an element of $P/\sim$. 
This is equivalent with proving that  $\Omega_K$ is equal to $$\Sigma(L/K) \cup [(\Omega(I_0 \cap I_0') ) \cap (\Omega(I_1 \cap I_1') )] \cup[(\Omega(I_0 \cap I_1') ) \cap (\Omega(I_0' \cap I_1) )] ,$$ and this follows from the equalities 
$$\Sigma(L/K) \cup \Omega(I_0) \cup \Omega(I_1) = \Omega_K,$$ and 
$$\Sigma(L/K) \cup  \Omega (I'_0) \cup \Omega(I'_1) = \Omega_K,$$  which hold as $(I_0,I_1)$ and $(I_0',I_1')$ are in $P/\sim$. 

\medskip

The following propositions will be used in Sections 6 and 8 in order to give necessary and sufficient conditions for the Hasse principle to hold. Let us
start with introducing some notation.

\bigskip
Set ${\cal C}_I = \{ (i,j) \in I \times I \ |  \ i \not = j \ {\rm and} \ \  \Sigma (L/K) \cup \Sigma_i \cup \Sigma_j  \not =\Omega_K \}$. For any subset $J$ of $I$, we say that $i,j \in J$ are
{\it connected in $J$} if there exist  $j_1,\dots,j_k \in J$ with $j_1 = i$, $j_k = j$ and $(j_r,j_{r+1}) \in {\cal C}_I$ for all $r = 1,\dots,k-1$. 

\bigskip
\noindent
{\bf Lemma 4.1.2.} {\it  
Let $(i,j) \in {\cal C}_I$, and let $v \in \Omega_K$ such that $v \not \in  \Sigma (L/K) \cup \Sigma_i \cup \Sigma_j$. 
Let $a^u_r \in  (F^u_r){^\times}$ for all $r \in I$ and $u \in \Omega_K$.  Then there exist $b^u_r \in  (F^u_r){^\times}$  such that

\medskip
\noindent
$\bullet$ $b^u_r = a^u_r$  whenever $u \not =v$ or $r \not = i,j$, and

\medskip
\noindent
$\bullet$ 
$ {\rm cor}_{F^v_i/K_v}  (b^v_i,d_i) \not =  {\rm cor}_{F^v_i/K_v}  (a^v_i,d_i)$,   $ {\rm cor}_{F^v_i/K_v}  (b^v_j,d_j) \not =  {\rm cor}_{F^v_i/K_v}  (a^v_j,d_j)$.

\medskip
\noindent
In particular, we have $$\Sigma_{v \in \Omega_K}  \ {\rm cor}_{F^v_i/K_v} (a^v_i,d_i)  \not = \Sigma_{v \in \Omega_K}  \ {\rm cor}_{F^v_i/K_v} (b^v_i,d_i),$$ 
$$\Sigma_{v \in \Omega_K}  \ {\rm cor}_{F^v_j/K_v} (a^v_j,d_j)  \not = \Sigma_{v \in \Omega_K}  \ {\rm cor}_{F^v_j/K_v} (b^v_j,d_j).$$}

\noindent{\bf Proof.} Since $(i,j) \in {\cal C}_I$, we have $\Sigma (L/K) \cup \Sigma_i \cup \Sigma_j  \not =\Omega_K$. Hence by Chebotarev's density theorem, the complement of  the set $\Sigma (L/K) \cup \Sigma_i \cup \Sigma_j $ contains
finite places.  Let us choose a finite place $v$ of $K$ such that $v \not \in \Sigma (L/K) \cup \Sigma_i \cup \Sigma_j $. As $v \not \in \Sigma_i$, we
have $E_i^v = E'_i \times M$, where $M$ is a field stable by $\sigma$, and $M^{\sigma} \not = M$. Set $M_0 = M^{\sigma}$. Similarly, we have $E_j^v = E'_j \times N$, where $N$ is
a field stable by $\sigma$, and $N^{\sigma} \not = N$. Set $N_0 = N^{\sigma}$. Then $M/M_0$ and $N/N_0$ are quadratic extensions of local fields. Let $\gamma \in M_0$
such that $\gamma \not  \in {\rm N}_{M/M_0}(M)$, and let $\delta \in N_0$ such that $\delta \not \in {\rm N}_{N/N_0}  (N)$. 
Let us write $a_i^v = (\alpha_1,\alpha_2)$
with $\alpha_1 \in (E'_i)^{\sigma}$, $\alpha_2 \in M_0$, and  $a_j^v = (\beta_1,\beta_2)$ with $\beta_1 \in (E'_j)^{\sigma}$, $\beta_2 \in N_0$. 

\medskip
Set $b_i^v = (\alpha_1, \alpha_2 \gamma)$ and $b_j^v = (\beta_1, \beta_2 \delta)$. If  $r \in I$ is such that $r \not = i,j$, then set $b_r^v = a_r^v$. For all $u \not = v$, set
$b_r^u = a_r^u$ for all $r \in I$. Then  $b^u_r \in  (F^u_r){^\times}$ have the required properties for all $u \in \Omega_K$ and $r \in I$. This completes the
proof of the Lemma.

\medskip
\noindent
{\bf Proposition 4.1.3.} 
{Let $i,j \in I$ be connected, and let  $a^u_r \in  (F^u_r){^\times}$ for all $r \in I$ and $u \in \Omega_K$.  Then there exist $b^u_r \in  (F^u_r){^\times}$  satisfying
the following conditions

\smallskip
\noindent
{\rm (i)}
\ \ \ $\Sigma_{v \in \Omega_K}  \ {\rm cor}_{F^v_i/K_v} (a^v_i,d_i)  \not = \Sigma_{v \in \Omega_K}  \ {\rm cor}_{F^v_i/K_v} (b^v_i,d_i).$

\smallskip
\noindent
{\rm (ii)} \ \ $\Sigma_{v \in \Omega_K}  \ {\rm cor}_{F^v_j/K_v} (a^v_j,d_j)  \not = \Sigma_{v \in \Omega_K}  \ {\rm cor}_{F^v_j/K_v} (b^v_j,d_j).$

\smallskip
\noindent
{\rm (iii)}  \  \ If  $r \not = i, j$, then we have 
$\Sigma_{v \in \Omega_K}  \ {\rm cor}_{F^v_r/K_v} (a^v_r,d_r)  = \Sigma_{v \in \Omega_K}  \ {\rm cor}_{F^v_r/K_v} (b^v_r,d_r).$

\smallskip
\noindent
{\rm (iv)} \  For all $v \in \Omega_K$, we have 
$\Sigma_{i  \in I}  \ {\rm cor}_{F^v_i/K_v} (b^v_i,d_i) = \Sigma_{i  \in I}  \  {\rm cor}_{F^v_i/K_v} (a^v_i,d_i).$

\smallskip
\noindent
{\rm (v)}  \ \ If $v$ is an infinite place of $K$, then $b^v_r = a^v_r$ for all $r \in I$. }

\medskip
\noindent
{\bf Proof.}
Let $j_1,\dots,j_k \in J$ with $j_1 = i$, $j_k = j$ and $(j_s,j_{s+1}) \in {\cal C}_I$ for all $s = 1,\dots,k-1$. Starting with $a^u_r \in  (F^u_r){^\times}$,
let us apply Lemma 4.1.2. successively to 
each of the pairs $(j_s,j_{s+1})$, and let $b^u_r \in (F^u_r){^\times}$ be the elements obtained at the end of the process. 

\medskip
Note that if $s \not = 1, k$, then we applied  Lemma 4.1.2.  twice. Hence we have
$${\rm cor}_{F^v_{j_s}/ K_v} (b^v_{j_s},d_{j_s}) = 
{\rm cor}_{F^v_{j_s}/ K_v} (a^v_{j_s},d_{j_s})$$ 
for $s \not  = 1,k$ and for all $v \in \Omega_K$.

\medskip
On the other hand,  if $s = 1$ or $s = k$, then we applied  Lemma 4.1.2. only once. Note also that $j_1 = i$ and $j_k = j$. Therefore we have  $${\rm cor}_{F^v_i/K_v} (b^v_i,d_i) \not = {\rm cor}_{F^v_i/K_v} (a^v_i,d_i)$$ for a certain $v \in \Omega_K,$ and $${\rm cor}_{F^v_i/K_v} (b^u_i,d_i) = {\rm cor}_{F^v_i/K_v} (a^u_i,d_i)$$ for all  $u \in \Omega_K$ with $u \not = v$. Similarly, we have 
$${\rm cor}_{F^w_j/K_w} (b^w_j,d_j) \not = {\rm cor}_{F^w_j/K_w} (a^w_j,d_j)$$ for a certain $w \in \Omega_K,$ and 
$${\rm cor}_{F^u_j/K_u} (b^u_j,d_j)  = {\rm cor}_{F^u_j/K_u} (a^u_j,d_j)$$ 
for all  $u \in \Omega_K$ with $u \not = w$. Therefore we have 
$$\Sigma_{v \in \Omega_K}  \ {\rm cor}_{F^v_i/K_v} (a^v_i,d_i)  \not = \Sigma_{v \in \Omega_K}  \ {\rm cor}_{F^v_i/K_v} (b^v_i,d_i),$$ 
$$\Sigma_{v \in \Omega_K}  \ {\rm cor}_{F^v_j/K_v} (a^v_j,d_j)  \not = \Sigma_{v \in \Omega_K}  \ {\rm cor}_{F^v_j/K_v} (b^v_j,d_j).$$ 

\medskip

Note that $b^v_r = a^v_r$ for all $v \in \Omega_K$ if $r \not = i, j$, hence we have 
$$\Sigma_{v \in \Omega_K}  \ {\rm cor}_{F^v_r/K_v} (a^v_r,d_r)  = \Sigma_{v \in \Omega_K}  \ {\rm cor}_{F^v_r/K_v} (b^v_r,d_r).$$ 

Moreover, all  the applications of  Lemma 4.1.2.  concern a place $v \in \Omega_K$ and two distinct indices $(j_s,j_{s+1}) \in {\cal C}_I$. This implies
that for all $v \in \Omega_K$, we have 
$$\Sigma_{i  \in I}  \ {\rm cor}_{F^v_i/K_v} (b^v_i,d_i) = \Sigma_{i  \in I}  \  {\rm cor}_{F^v_i/K_v} (a^v_i,d_i)$$

All the changes were made at finite places, hence we have $b^v_r = a^v_r$ for all $r \in I$ if $v$ is an infinite place. This completes the
proof of the Proposition.

\medskip
\noindent
{\bf Proposition 4.1.4.} {\it Let $a^v_i \in (F_i^v)^{\times}$ for all $v \in \Omega_K$, $i \in I$, such that :

\smallskip
\noindent
{\rm (i)} We have  $$\Sigma_{v \in \Omega_K} \Sigma_{i \in I} \ {\rm cor}_{F^v_i/K_v}  (a^v_i,d_i) = 0.$$
\noindent
{\rm (ii)} For all $x \in \sha$,  we have
$$\Sigma_{v \in \Omega_K} \Sigma_{i \in I_0(x)}  \ {\rm cor}_{F^v_i/K_v} (a^v_i,d_i) = 0.$$

Then there exist $b^v_i \in F_i^v$ for all $v \in \Omega_K$, $i \in I$ such that :

\smallskip
\noindent
{\rm (iii)} For all $i \in I$, we have

$$\Sigma_{v \in \Omega_K}  \ {\rm cor}_{F^v_i/K_v} (b^v_i,d_i) = 0.$$

\smallskip
\noindent
{\rm (iv){ For all $v \in \Omega_K$, we have

$$\Sigma_{i  \in I}  \ {\rm cor}_{F^v_i/K_v} (b^v_i,d_i) = \Sigma_{i  \in I}  \  {\rm cor}_{F^v_i/K_v} (a^v_i,d_i).$$

\smallskip
\noindent
{\rm (v)} if $v$ is an infinite place of $K$, then $b^v_i = a^v_i$.}

\medskip

\noindent
{\bf Proof.} For all $i \in I$, set $C_i = C_i(a) = \Sigma_{v \in \Omega_K}  \ {\rm cor}_{F^v_i/K_v} 
(a^v_i,d_i).$
If $C_i = 0$ for all $i \in I$, we set $b_i^v = a_i^v$ for all $i \in I$ and $v \in \Omega_K$. If not, then we construct a connected graph with vertex set $\cal V$ and edge set $\cal E$ in order to make successive modifications.

\medskip
Our aim is to construct a  graph containing two elements $i_0, i_k  \in I$  such that $C_{i_0} = C_{i_k} = 1$ and that $i_0$ and $i_k$ are connected within the
graph. 

\medskip
Let us now construct the desired graph with vertex set $\cal V$ and edge set $\cal E$. We start with the empty graph, and add edges and vertices as follows. Let us choose $i_0 \in I$ such that $C_{i_0} = 1$, and add $ \{ i_0 \}$ to $\cal V$. Set  $I_0 = \{ i_0 \}$ and $I_1 = I - I_0$. Note that $(I_0,I_1) \not \in \sha$.
Indeed, if $(I_0,I_1) \in \sha$, then by {\rm (ii)} we have
$\Sigma_{v \in \Omega_K} \Sigma_{i \in I_0}  \ {\rm cor}_{F^v_i/K_v} 
(a^v_i,d_i) = 0.$ But $\Sigma_{v \in \Omega_K} \Sigma_{i \in I_0}  \ {\rm cor}_{F^v_i/K_v} 
 (a^v_i,d_i)  = C_{i_0}$, and
$C_{i_0} = 1$, so this leads to a contradiction. Therefore by definition of $\sha$, we have
$$\Sigma (L/K) \cup (\underset{i\in I_0}{\cap}\Sigma_i)\cup(\underset{j\in I_1}{\cap}\Sigma_j) \not =\Omega_K.$$
Hence there exist $i_1 \in I_1$ and $v \in \Omega_K$ such that $v \not \in  \ \Sigma (L/K) \cup  \Sigma_{i_0} \cup \Sigma_{i_1}$. In other words, we have
$(i_0,i_1) \in {\cal C}_I$, hence $i_0$ and $i_1$ are connected. Add $\{  i_1 \}$ to  ${\cal V}$, and
add the edge connecting $i_0$ to $i_1$ to $\cal E$.  If $C_{i_1} = 1$, we stop. If not, set $I_0 = \{ i_0, i_1 \}$ and $I_1 = I - I_0$.  We again have 
$(I_0,I_1) \not \in \sha$.
Indeed, if $(I_0,I_1) \in \sha$, then by {\rm (ii)} we have
$\Sigma_{v \in \Omega_K} \Sigma_{i \in I_0}  \ {\rm cor}_{F^v_i/K_v} 
 (a^v_i,d_i) = 0.$ But $\Sigma_{v \in \Omega_K} \Sigma_{i \in I_0} \ {\rm cor}_{F^v_i/K_v} 
(a^v_i,d_i)  = C_{i_0} + C_{i_1}$, and
$C_{i_0} = 1$, $C_{i_1} = 0$, so this is again a contradiction. 
Therefore by definition of $\sha$, we have
$$\Sigma (L/K) \cup (\underset{i\in I_0}{\cap}\Sigma_i)\cup(\underset{j\in I_1}{\cap}\Sigma_j) \not =\Omega_K.$$ Hence there exists $i_2 \in I_1$ and $v \in \Omega_K$ such that $v \not \in 
\Sigma (L/K) \cup (\underset{i\in I_0}{\cap}\Sigma_i) \cup \Sigma_{i_2}$. This implies that at least one of $(i_0,i_2), (i_1,i_2)$ belong to ${\cal C}_I$.
We now add $i_2$ to $\cal V$, and add to $\cal E$
all the edges connecting $j$ to $i_2$ with $j \in {\cal V}$ such that $(j,i_2) \in {\cal C}_I$. Note that $i_0$ and $i_2$ are connected within the graph. We
continue this way, adding vertices to $\cal V$ and edges to $\cal E$. 
Since $I$ is finite, and since by {\rm (i)} there exists $j \in I$ with $j \not = i_0$ and $C_j = 1$, the process will stop after a finite number of steps. 

\medskip
In other words,  after a finite number of steps we find $i_k \in I$ such that $C_{i_k} = 1$, and such that the resulting graph with vertices $\cal V$ and edges $\cal E$
has the following property :  there exists a loop--free path in $\cal E$ connecting $i_0$ to $i_k$ such that  for any two adjacent vertices $i, j \in {\cal V}$ we
have $(i,j) \in {\cal C}_I$.  In other words, $i_0$ and $i_k$ are connected in $\cal V$.
By Proposition 4.1.3. this implies that there exist
 $c^v_i \in F_i^v$ for all $v \in \Omega_K$, $i \in I$ such that for $(c) = (c^v_i)$ we have $C_{i_0}(c) = C_{i_k}(c) = 0$ and $C_i(c) = C_i(a)$ for all $i \not = i_0, i_k$.
Therefore the number of $i \in I$ with $C_i(c) = 1$ is less than the number of $i \in I$ with $C_i(a) = 1$.  Moreover, for all $v \in \Omega_K$, we have
$\Sigma_{i  \in I}  \ {\rm cor}_{F^v_i/K_v} (c^v_i,d_i) = \Sigma_{i  \in I}  \  {\rm cor}_{F^v_i/K_v} (a^v_i,d_i),$  and that if $v$ is an infinite place, then 
$c^v_i = a^v_i$ for all $i \in I$. 
Continuing this way leads to the desired conclusion : we obtain $b^v_i \in F_i^v$ for all $v \in \Omega_K$, $i \in I$ such that for $(b) = (b^v_i)$ 
we have
$C_i (b) = \Sigma_{v \in \Omega_K}  \ {\rm cor}_{F^v_i/K_v} (b^v_i,d_i) = 0$,  for all $i \in I$, and this implies {\rm (iii)}.
Note that 
for all $v \in \Omega_K$, we have
$\Sigma_{i  \in I}  \ {\rm cor}_{F^v_i/K_v} (b^v_i,d_i) = \Sigma_{i  \in I}  \  {\rm cor}_{F^v_i/K_v} (a^v_i,d_i).$ This implies that {\rm (iv)} holds.
Moreover, all the modifications were made at finite places, hence {\rm (v)} holds.

\medskip
\noindent
{\bf Proposition 4.1.5.} {\it Let $a^v_i \in( F_i^v)^{\times}$ for all $v \in \Omega_K$, $i \in I$, such that :

\smallskip
\noindent
{\rm (i)} We have  $$\Sigma_{v \in \Omega_K} \Sigma_{i \in I} \ {\rm cor}_{F^v_i/K_v}  (a^v_i,d_i) = 0.$$
\noindent
{\rm (ii)} For all $x \in \sha$,  we have
$$\Sigma_{v \in \Omega_K} \Sigma_{i \in I_0(x)}  \ {\rm cor}_{F^v_i/K_v} (a^v_i,d_i) = 0.$$

Then for all $i \in I$, there exist $b_i \in F_i^{\times}$ such that

\smallskip
\noindent
{\rm (iii)} For all $v \in \Omega_K$, we have $$\Sigma_{i  \in I}  \ { \rm cor}_{F^v_i/K_v} (b_i,d_i)  = \Sigma_{i  \in I}  \ {\rm cor}_{F_i^v/K_v} (a_i^v,d_i).$$}

\noindent
{\bf Proof.} By Proposition 4.1.4. conditions {\rm (i)} and {\rm (ii)} imply that for all $v \in \Omega_K$ and all $i \in I$, there exist
$b^v_i \in (F_i^v)^{\times}$ such that 
for all $i \in I$, we have
$\Sigma_{v \in \Omega_K}  \ {\rm cor}_{F^v_i/K_v} (b^v_i,d_i) = 0,$ and that 
for all $v \in \Omega_K$, we have
$$\Sigma_{i  \in I}  \ {\rm cor}_{F^v_i/K_v} (b^v_i,d_i) = \Sigma_{i  \in I}  \  {\rm cor}_{F^v_i/K_v} (a^v_i,d_i).$$
Let $i \in I$. Since $\Sigma_{v \in \Omega_K}  \ {\rm cor}_{F^v_i/K_v} (b^v_i,d_i) = 0,$ we have $ \Sigma_{w \in \Omega_{F_i}}   (b^w_i,d_i) = 0$.
The Brauer--Hasse--Noether Theorem implies that
there exists a quaternion algebra $Q_i$ over $F_i$ such that for all $v \in \Omega_K$, we have  $Q_i \simeq (b_i^v,d_i)$. Since $Q_i^v$ splits
over $E_i^v$ for all $v \in \Omega_K$, the algebra $Q_i$ splits over $E_i$. Therefore there exists $b_i \in (F_i)^{\times}$ such that
$Q_i \simeq (b_i,d_i)$. 

\medskip
Then, for all $v \in \Omega_K$,  we have
$${\Sigma_{i  \in I}  \   \rm cor}_{F_i/K} (b_i,d_i)  =  \Sigma_{i  \in I}  \ {\rm cor}_{F^v_i/K_v} (b^v_i,d_i) = \Sigma_{i  \in I}  \  {\rm cor}_{F^v_i/K_v} (a^v_i,d_i).$$
Therefore {\rm (iii)} holds.  This completes the proof of the proposition.

\bigskip

{\bf 5. The Brauer--Manin map}

\medskip Assume that $K$ is a global field, and that
$(E^v,\sigma)$ can be embedded in $(A^v,\tau)$ for all $v \in \Omega_K$. This implies that there exists an
embedding of algebras $\epsilon : E \to A$. By Proposition 1.1.2. there exists an involution $\theta : A \to A$ of the same type as $\tau$ such that
$\epsilon$ induces an embedding
of algebras with involution $(E,\sigma) \to (A,\theta)$. Let us fix such an involution $\theta$. If $A \simeq M_n(K)$ and if $\tau$ is an
orthogonal involution, then let us chose for $\theta$ the involution induced by the quadratic form $T : E \times E$, given by $T(x,y) = {\rm Tr}_{E/K}(x \sigma(y))$ for
all $x,y \in E$. Note that this is possible by Proposition 1.4.1.

\medskip
The aim of this section is to define a map $\sha(E,\sigma)  \to {\bf Z}/2{\bf Z}$ the vanishing of which is a necessary and sufficient condition for
the existence of an embedding of algebras with involution $(E,\sigma) \to (A,\tau)$. To define this map, we need the notion of {\it embedding data} 
(cf. 5.1.-5.3.). The Brauer--Manin map is defined in 5.4. 

\bigskip
{\bf 5.1. Local embedding data -- even degree orthogonal case}

\medskip
Assume that $(A,\tau)$ is an orthogonal involution, with $A$ of degree $n$. Assume that $n$ is even, and set  $n = 2r$. 
Let us fix an isomorphism of $K$--algebras  $u : \Delta (E) \to Z(A,\theta)$, and recall  (cf. 2.5.) that for all $a^v \in (F^v)^{\times}$ this induces 
a uniquely defined  isomorphism of $K_v$--algebras $u_{a^v} : \Delta (E^v) \to Z(A,\theta_{a^v})$.

\medskip
We are assuming that for all $v \in \Omega_K$, there exists an embedding of algebras with involution $(E^v,\sigma) \to (A^v,\tau)$. This implies that the $K$--algebras $\Delta(E)$ and $Z(A,\tau)$ are
isomorphic. 
Let us fix an isomorphism of $K$--algebras $$\nu : \Delta (E) \to Z(A,\tau).$$

\medskip
Let us denote by ${\cal O}(E,A)$ the set of $(a) = (a^v)$, with $a^v \in (F^v)^{\times}$, such that 
for all $v \in \Omega_K$, there exists $\alpha^v \in (A^v)^{\times}$ with the properties :

\smallskip  \noindent
{\rm (a)}  ${\rm Int}(\alpha) : (A^v,\theta_{a^v}) \to (A^v,\tau)$ is an isomorphism
of $K_v$--algebras with involution.

\smallskip  \noindent
{\rm (b)} The induced automorphism  $c(\alpha) :  Z(A^v,\theta_{a^v})  \to Z(A^v,\tau)$ satisfies $$c(\alpha) \circ   u_{a^v} = \nu.$$

In other words, $({\rm Int}(\alpha) \circ \epsilon, a^v,\alpha^v,\nu)$ are parameters of an oriented embedding.

\medskip
\noindent
{\bf Proposition 5.1.1.} {\it Let $(a) =  (a^v) \in {\cal O}(E,A)$. Then we have :

\smallskip
{\rm (i)} \ \ \  $ {\rm res}_{\Delta (E^v)/K_v } \ {\rm cor}_{F^v/K_v}(a^v,d) = 0$ for almost all $v \in \Omega_K$, and $$ \Sigma_{v \in \Omega_K} {\rm res}_{\Delta (E^v)/K_v } \ {\rm cor}_{F^v/K_v}(a^v,d) = 0.$$

\smallskip
{\rm (ii)} \ Let $\Omega'$ be the set of places $v \in \Omega_K$ such that $\Delta (E^v) \simeq K_v \times K_v$. Then we have ${\rm cor}_{F^v/K_v}(a^v,d) = 0$ for almost all $v \in \Omega'$, and 
$$ \Sigma_{v \in \Omega'} \ {\rm cor}_{F^v/K_v}(a^v,d) = 0.$$}

\noindent
{\bf Proof.} By Lemma 2.5.4.  $C(A^v,\theta_{a^v}) = C(A^v,\theta) + {\rm res}_{\Delta (E^v)/K_v } \ {\rm cor}_{F^v/K_v}(a^v,d) $ in ${\rm Br}(\Delta (E^v)$
for all $v \in \Omega_K$. Since  $(a^v) \in {\cal O}(E,A)$ we have $C(A^v,\theta_{a^v}) = C(A^v,\tau)$ for all $v \in \Omega_K$.
Therefore we have $$C(A^v,\tau)  - C(A^v,\theta) = {\rm res}_{\Delta (E^v)/K_v } \ {\rm cor}_{F^v/K_v}(a^v,d),$$ hence 
 {\rm (i)} holds.  If $\Delta (E^v) \simeq K_v \times K_v$, then ${\rm res}_{\Delta (E^v)/K_v } $ is injective, and this implies {\rm (ii)}.

\medskip
\noindent
{\bf Proposition 5.1.2.} {\it Let $(a^v), (b^v)  \in {\cal O}(E,A)$. Then, for all $v \in \Omega_K$ 

\smallskip
{\rm (i)} \ \ $ {\rm res}_{\Delta (E^v)/K_v } \ {\rm cor}_{F^v/K_v}(a^v,d) = {\rm res}_{\Delta (E^v)/K_v } \ {\rm cor}_{F^v/K_v}(b^v,d)$.

\smallskip
{\rm (i)} \ \ If moreover $\Delta(E^v) \simeq K_v \times K_v$, then $ {\rm cor}_{F^v/K_v}(a^v,d) = {\rm cor}_{F^v/K_v}(b^v,d)$.}

\medskip
\noindent
{\bf Proof.} We have $C(A^v,\theta_{a^v}) = C(A^v,\theta) + {\rm res}_{\Delta (E^v)/K_v } \ {\rm cor}_{F^v/K_v}(a^v,d),$
and  $C(A^v,\theta_{b^v}) = C(A^v,\theta) + {\rm res}_{\Delta (E^v)/K_v } \ {\rm cor}_{F^v/K_v}(a^v,d) $ in ${\rm Br}(\Delta (E^v))$ (cf. Lemma 2.5.4.). Since 
$(a^v), (b^v) \in {\cal O}(E,A)$ we have $C(A^v,\theta_{a^v}) = C(A^v,\theta_{b^v})$, and this implies {\rm (i)}.
If $\Delta(E^v) \simeq K_v \times K_v$, then ${\rm res}_{\Delta (E^v)/K_v } $ is injective, hence  {\rm (ii)}.

\bigskip

A {\it local embedding datum} will be a set  $(a) =  (a^v) \in {\cal O}(E,A)$ such that

\smallskip
\noindent
$\bullet$ If  $v \in \Omega_K$ is such that $\Delta (E_v)$ is a quadratic extension of $K_v$, then there exist only finitely
many $v \in \Omega_K$ such that  ${\rm cor}_{F^v/K_v}(a^v,d) \not = 0$.

\smallskip
\noindent
$\bullet$ We have $$\Sigma _{v \in \Omega_K} \ {\rm cor}_{F^v/K_v}(a^v,d)  = 0.$$

\medskip
We denote by ${\cal L}(E,A)$  the set of local embedding data. 

\medskip
\noindent
{\bf Remark.} Let $(a^v) \in {\cal L}(E,A)$. Then we have ${\rm cor}_{F^v/K_v}(a^v,d) = 0$ for almost all $v \in \Omega_K$.  Indeed, by hypothesis
this is true if $v$ is such that $\Delta (E_v)$ is a quadratic extension of $K_v$, and  by Proposition 5.1.1.  {\rm (ii)} it also holds if $v$ is such that $\Delta (E_v) \simeq K_v \times K_v$.

\medskip
Recall that the notion of oriented embedding was defined in 2.6. 

\bigskip
\noindent
{\bf Proposition 5.1.3.} {\it Assume that for all $v \in \Omega_K$, there exists an oriented embedding  $(E^v,\sigma) \to (A^v,\tau)$ with respect to $\nu$. Then
there exists a local embedding datum $(a) = (a^v) \in {\cal L}(E,A)$ such that 
for all $v \in \Omega_K$ there exist  $\iota_v $ and $\alpha^v$ 
such that $(\iota_v,a^v,\alpha^v, \nu)$ are  parameters of an oriented embedding.}

\medskip
\noindent
{\bf Proof.}
{\rm Case 1.} Assume that $\Delta (E^v)/K_v$ is a quadratic extension.  Let $(b^v) \in {\cal O}(E,A)$. Then 
$C(A^v,\tau) = C(A^v,\theta)  + {\rm res}_{\Delta (E^v)/K_v } \ {\rm cor}_{F^v/K_v}(b^v,d) = C(A^v,\theta)$ in ${\rm Br} (\Delta (E^v))$, since $\Delta (E^v)/K_v$
is a quadratic extension. Moreover, we have ${\rm disc}(A^v,\tau) =  {\rm disc}(A^v,\theta_{b^v} ) = {\rm disc}(A^v,\theta)$. Hence $(A^v, \theta)$ and $(A^v,\tau)$ are isomorphic. 
By Corollary 2.7.3. {\rm (ii)}  there exist $\iota_v$ and $\alpha^v$ such that $(\iota_v,1,\alpha^v,\nu)$ are parameters of an oriented embedding.

\medskip
Case 2.  Assume now that we have $\Delta(E^v) \simeq K_v \times K_v$. Let 
$(\iota_v,a^v, \alpha^v,\nu)$ be parameters for an oriented embedding.

\medskip Let $(a) = (a^v)$, where for $v \in \Omega_K$ the element $a^v$ is chosen as above, in each of the two cases. We claim that  
$(a) = (a^v) \in {\cal L}(E,A)$. Since $a^v = 1$ when $\Delta (E^v)/K_v$ is a quadratic extension, we have ${\rm cor}_{F^v/K_v}(a^v,d) = 0$  for
all such $v$. Let $\Omega'$ be the set of $v \in \Omega_K$ such that $\Delta(E^v) \simeq K_v \times K_v$.  Then we have 
$\Sigma _{v \in \Omega_K} \ {\rm cor}_{F^v/K_v}(a^v,d) = \Sigma _{v \in \Omega'} \ {\rm cor}_{F^v/K_v}(a^v,d),$ and by Proposition 5.1.1. {\rm (ii)}
this sum is zero. Therefore 
we have $(a)  \in {\cal L}(E,A)$.

\bigskip
\noindent
{\bf Proposition 5.1.4.} {\it Let $(a) = (a^v), (b) = (b^v)   \in {\cal L}(E,A)$. Then there exists $\lambda \in K^{\times}$ such that for all $v \in \Omega_K$
we have ${\rm cor}_{F^v/K_v}  (\lambda b^v,d) = {\rm cor}_{F^v/K_v}  (a^v,d)$.}

\medskip
\noindent
{\bf Proof.} 
We have $ {\rm res}_{\Delta (E^v)/K} \ {\rm cor}_{F^v/K_v}(a^v,d) = {\rm res}_{\Delta (E^v)/K} \ {\rm cor}_{F^v/K_v}(b^v,d)$
for all $v \in \Omega_K$, and if  $\Delta (E^v) \simeq K_v \times K_v$, then ${\rm cor}_{F^v/K_v}  ( b^v,d) = {\rm cor}_{F^v/K_v}  (a^v,d)$ (cf. Proposition 5.1.2.).

\medskip
Let $\Omega' = \{v \in \Omega_K \ | \  {\rm cor}_{F^v/K_v}  ( b^v,d) \not = {\rm cor}_{F^v/K_v}  (a^v,d) \}$. The above argument shows that if $v \in \Omega'$, then
$\Delta (E^v)$ is a quadratic extension of $K_v$. 
It follows from  the definition of ${\cal L}(E,A)$ that there exist only finitely many $v \in \Omega_K$ such that $ {\rm cor}_{F^v/K_v}  (a^v,d) \not = 0$
or $ {\rm cor}_{F^v/K_v}  (b^v,d)  \not = 0$, hence $\Omega'$ is a finite set. 

\medskip
Let $v \in \Omega'$. Then $\Delta (E^v)$ splits $ {\rm cor}_{F^v/K_v}  (b^v,d) - {\rm cor}_{F^v/K_v}  (a^v,d)$. Recall that $\Delta (E^v) = K_v(\sqrt D)$, where
$D = (-1)^r {\rm N}_{E/K}(\sqrt d) =  {\rm N}_{F/K}(d)$. Then
we have $ {\rm cor}_{F^v/K_v}  (b^v,d) - {\rm cor}_{F^v/K_v}  (a^v,d) = (\lambda^v,D)$ for some $\lambda^v \in K_v^{\times}$. Since 
 $(a) , (b)  \in {\cal L}(E,A)$, by definition we have $\Sigma_{v \in \Omega_K}  \ {\rm cor}_{F^v/K_v}  (a^v,d)  = \Sigma_{v \in \Omega_K}  \ {\rm cor}_{F^v/K_v}  (b^v,d)  =  0.$ This implies that $\Sigma_{v \in \Omega_K} (\lambda^v,D) = 0.$ Hence by the Brauer--Hasse--Noether theorem, there exists $\lambda \in K^{\times}$ such
 that $(\lambda,D) = (\lambda^v,D) \in {\rm Br}(K_v)$ for all $v \in \Omega'$, and $\lambda$ has the required property.

\bigskip

{\bf 5.2. Local embedding data -- odd degree orthogonal case}

\medskip
In this section, we assume that $A \simeq M_n(K)$  and that $\tau$ is induced by an $n$--dimensional  quadratic form $q$. 
We are primarily interested in the case where $n$ is odd, but  we also need to consider the case where $n$ is even.

\medskip
Let us assume that there exists an embedding of algebras with involution $(E^v,\sigma) \to (A^v,\tau)$ for all $v \in \Omega_K$. By
Proposition 1.4.1. this implies that for all $v \in \Omega_K$ there exists $a^v \in (F^v)^{\times}$ such that $q \simeq T_{a^v}$. Let
us write $a^v = (a_1^v,\dots,a^v_m)$ with $a^v_i \in (F_i^v)^{\times}$. The set of $(a) = (a_i^v)$ with this property will
be denoted by ${\cal L}'(E,A)$.

\bigskip \noindent
{\bf Proposition 5.2.1.} {\it Let $(a) \in {\cal L}'(E,A)$, with $(a) = (a_i^v)$. Then the following properties hold :

\medskip
{\rm (i)}  \ \  ${\rm cor}_{F^v/K_v} (a^v,d) = 0$ for almost all $v \in \Omega_K$, and  $$\Sigma _{v \in \Omega_K} \ {\rm cor}_{F^v/K_v}(a^v,d)  = 0.$$

\medskip
{\rm (ii)} Let $(b) \in {\cal L}'(E,A)$, with $(b) = (b_i^v)$. Then for all $v \in \Omega_K$, we have
$${\rm cor}_{F^v/K_v} (a^v,d) =
{\rm cor}_{F^v/K_v} (b^v,d).$$}  

\noindent
{\bf Proof.} 
Let us first assume that $n$ is even. Since $(a) \in {\cal L}'(E,A)$,  we have  $q \simeq T_{a^v}$, and hence $w(T_{a^v})  = w(q)$ for all $v \in \Omega_K$.
By Lemma 1.5.1. we have
$w(T_{a^v})  = w(T) + {\rm cor}_{F^v/K_v} (a^v,d).$
Hence for all $v \in \Omega_K$, we have 
$w(q) = w(T_{a^v})  = w(T) + {\rm cor}_{F^v/K_v} (a^v,d)$.
Note that $\Sigma_{v \in \Omega_K} w(q) = \Sigma_{v \in \Omega_K} w(T) = 0$. Therefore we have
$\Sigma_{v \in \Omega_K} {\rm cor}_{F^v/K_v} (a^v,d) = 0$, and this proves {\rm (i)}.

\medskip

Let us prove {\rm (ii)}. Since $(b) \in {\cal L}'(E,A)$, for all $v \in \Omega_K$ we have $w(T_{b^v})  = w(q)$.  By Lemma 1.5.1. we have
$w(T_{b^v})  = w(T) + {\rm cor}_{F^v/K_v} (b^v,d)$
for all $v \in \Omega_K$. Therefore, for all $v \in \Omega_K$, we have $w(T) + {\rm cor}_{F^v/K_v} (a^v,d) =  w(q) = w(T) + {\rm cor}_{F^v/K_v} (b^v,d).$ Hence we have ${\rm cor}_{F^v/K_v} (a^v,d) = {\rm cor}_{F^v/K_v} (b^v,d)$, and this implies {\rm (ii)}. 

\medskip
Suppose that $n$ is odd, and set $A' = M_{n-1}(K)$. Then by [PR 10], Proposition 7.2. there exists a $\sigma$--invariant  \'etale subalgebra $E'$ of rank $n-1$ of $E$ such that $E = E' \times K$,
an $(n-1)$--dimensional quadratic form $q'$ and a 1--dimensional quadratic form $q''$ over $K$  such that $q \simeq q' \oplus q''$, and that the \'etale algebra with involution $(E',\sigma)$ can be embedded
in the central simple algebra $(A',\tau')$ over $K_v$ for all $v \in \Omega_K$, where $\tau' : A' \to A'$ is the involution induced by $q'$. Moreover, there exists an embedding of $(E,\sigma)$ into $(A,\tau)$ if and only if there exists an embedding of $(E',\sigma)$ into $(A',\tau')$. Note that we have ${\cal L}'(E,A) = {\cal L}' (E',A') \times {\cal L}'(K,K)$. We may suppose that $E_m = K$. Then we have
$d_m = 1$. Set $J = \{1,\dots,m-1 \}$, 
and note that  for all $v \in \Omega_K$, we have $\Sigma_{i \in I} {\rm cor}_{F_i^v/K_v}(a_i^v,d_i) = \Sigma_{i \in J} {\rm cor}_{F_i^v/K_v}(a_i^v,d_i)$. Since $n-1$
is  even, statements {\rm (i)} and {\rm (ii)} easily follow.

\bigskip
If $n$ is odd, then we set ${\cal L}(E,A) = {\cal L}'(E,A)$, and an element $(a) \in {\cal L}(E,A)$ will be called {\it local embedding datum}. 

\medskip
If $n$ is even, then the set of embedding data ${\cal L}(E,A)$ was defined in the previous section. The relationship between ${\cal L}(E,A)$
and ${\cal L}'(E,A)$ is as follows :

\medskip
\noindent
{\bf Proposition 5.2.2.} {\it Assume that $n$ is even. Then we have

\medskip
{\rm (i)} \ \ ${\cal L}'(E,A) \subset {\cal L}(E,A)$.

\medskip
{\rm (ii)} Let $(a) \in {\cal L}(E,A)$. Then there exists $\lambda \in K^{\times}$ such that $(\lambda a) \in {\cal L}'(E,A)$.}

\medskip
\noindent
{\bf Proof.} Let $(a) \in {\cal L}'(E,A)$. Then  ${\rm cor}_{F^v/K_v} (a^v,d) = 0$ for almost all $v \in \Omega_K$, and 
$\Sigma_{v \in \Omega_K} \ {\rm cor}_{F^v/K_v} (a^v,d) = 0$ (cf. Proposition 5.2.1. {\rm (i)}). Since $q \simeq T_{a^v}$ for all $v \in \Omega_K$, 
the algebras with involution $(A^v,\tau)$ and $(A^v,\theta_{a^v})$ are isomorphic. Since $A$ is split, Corollary 2.7.3. implies that for all $v \in \Omega_K$  there exist 
$\iota_v$ and $\alpha^v$ such that  $(\iota_v,a^v,\alpha^v,\nu)$ are parameters of an oriented embedding $(E^v,\sigma) \to (A^v,\tau)$.
This implies that $(a) \in {\cal L}(E,A)$, hence {\rm (i)} is proved.

\medskip

Let us prove {\rm (ii)}. Let $S$ be the finite set of places of $K$ at which $q$ or $T$ is not hyperbolic, or $(a^v,d) \not = 0$. Since $(a) \in {\cal L}(E,A)$,  there exists $\lambda^v \in K_v^{\times}$
such that $q$ and $\lambda^v T_{a^v}$ are isomorphic over $K_v$ for all $v \in S$. There exists $\lambda \in K^{\times}$ such that $\lambda (\lambda^v)^{-1} \in 
(K_v)^{\times 2}$ for all $v \in S$. Then $q$ and $\lambda T_{a^v}$ are isomorphic over $K_v$ for all $v \in S$. For $v \not \in S$, both $q$ and $T_{a^v}$ are
hyperbolic over $K_v$, hence we have $q \simeq \lambda T_{a^v}$. Since $\lambda T_{a^v} =  T_{\lambda a^v}$, we have $(\lambda a) \in {\cal L}'(E,A)$.

\bigskip

{\bf 5.3. Local embedding data -- the unitary case}

\medskip
Let us assume that  $(A,\tau)$ is a unitary involution.

\medskip

The set of $(a) = (a^v)$, with $a^v \in (F^{v})^{\times}$,  such that for all $v \in \Omega_K$  we have
$(A_v,\tau) \simeq (A_v,\theta_{a^v})$, is called a {\it local embedding datum}. We denote by ${\cal L}(E,A)$ the set of local embedding data.

\medskip
\noindent
{\bf Proposition 5.3.1.} {\it Let $(a) \in {\cal L}(E,A)$ be an embedding datum, with $(a) = (a_i^v)$. Then the following properties hold :

\medskip
{\rm (i)} We have $$\Sigma_{v \in \Omega_K} {\rm cor}_{F^v/K_v} (a^v,d)  = 0.$$

\medskip
{\rm (ii)} Let $(b) \in {\cal L}(E,A)$ be an embedding datum, with $(b) = (b_i^v)$. Then for all $v \in \Omega_K$, we have
$${\rm cor}_{F^v/K_v} (a^v,d) = {\rm cor}_{F^v/K_v} (b^v,d).$$}

\noindent
{\bf Proof.} 
Since $(a) \in {\cal L}(E,A)$, we have
$(A^v,\theta_{a^v}) \simeq (A^v,\tau)$ for all $v \in \Omega_K$. Hence for all $v \in \Omega_K$, we have $D(A^v,\theta_{a^v})  = D(A^v,\tau)$. By Lemma 1.5.2. we have $D(A^v,\theta_{a^v}) = D(A^v,\theta) + {\rm cor}_{F^v/K^v}(a^v,d)$
for all $v \in \Omega_K$.  
We have $\Sigma_{v \in \Omega_K} D(A^v,\tau)  = \Sigma_{v \in \Omega_K} D(A^v,\theta) = 0$, hence 
$\Sigma_{v \in \Omega_K} {\rm cor}_{F^v/K_v} (a^v,d) = 0$.
This proves {\rm (i)}.

\medskip Let us prove {\rm (ii)}.  Let $v \in \Omega_K$. Since $(b) \in {\cal L}(E,A)$, by Lemma 1.5.2. we have $D(A^v,\theta) + {\rm cor}_{F^v/K^v}(a^v,d) = D(A^v,\theta_{a^v}) = D(A^v,\tau) = 
D(A^v,\theta_{b^v}) = D(A^v,\theta) + {\rm cor}_{F^v/K^v}(b^v,d)$. 
Hence we have ${\rm cor}_{F^v/K_v} (a^v,d) = {\rm cor}_{F^v/K_v} (b^v,d)$, as claimed.

 \bigskip
 
 {\bf 5.4. The Brauer--Manin map}
 
 \bigskip
 Let $(a) \in {\cal L}(E,A)$ be an embedding datum, with $(a) = (a_i^v)$. Let us consider the map
 
 $$f_{(a)}  : \sha (E,\sigma) \to {\bf Z}/2{\bf Z}$$ defined by

$$f_{(a)} (I_0,I_1) = \Sigma_{i \in I_0} \Sigma_{v \in \Omega_K} {\rm cor}_{F^v_i/K_v}(a_i^v,d_i).$$

\medskip

Note that this is well--defined, since $\Sigma_{i \in I} \Sigma_{v \in \Omega_K} {\rm cor}_{F^v_i/K_v}(a_i^v,d_i) = 0.$
As we will see, this map is independent of the choice of $(a)$. In other words, we have

\bigskip
\noindent
{\bf Theorem 5.4.1.} {\it Let $(a), (b) \in {\cal L}(E,A)$ be two local embedding data. Then  we have $f_{(a)} = f_{(b)}$. }

\medskip
\noindent
{\bf Proof.} Suppose that $(a), (b) \in {\cal L}(E,A)$ are such that $f_{(a)} \not = f_{(b)}$. Note that for all $\lambda \in K^{\times}$, we have
$(\lambda b) \in {\cal L}(E,A)$, and $f_{(b)} = f_{(\lambda b)}$. Since there exists $\lambda \in  K^{\times}$ such that for all $v \in \Omega_K$ we
have $\Sigma_{i \in I} {\rm cor}_{F^v_i/K_v}(a_i^v,d_i) = \Sigma_{i \in I} {\rm cor}_{F^v_i/K_v}(\lambda b_i^v,d_i)$ (cf. Proposition 5.1.4. Proposition 5.2.1. {\rm (ii)} and
Proposition 5.3.1. {\rm (ii)}), we may assume that for all $v \in \Omega_K$, we have $$\Sigma_{i \in I} {\rm cor}_{F^v_i/K_v}(a_i^v,d_i) = \Sigma_{i \in I} {\rm cor}_{F^v_i/K_v}(b_i^v,d_i).$$

Let $(I_0,I_1) \in \sha(E,\sigma)$ be such that $f_{(a)}(I_0,I_1) \not  = f_{(b)}(I_0,I_1)$. Then there exists $v \in \Omega_K$ such that 
$\Sigma_{i \in I_0} {\rm cor}_{F^v_i/K_v}(a_i^v,d_i) \not = \Sigma_{i \in I_0} {\rm cor}_{F^v_i/K_v}(b_i^v,d_i).$ This implies that $v \not \in \Sigma(L/K)$, and
$v \not \in  \underset{i\in I_0} \cap \Sigma_i$. Since $\Sigma_{i \in I} {\rm cor}_{F^v_i/K_v}(a_i^v,d_i) = \Sigma_{i \in I} {\rm cor}_{F^v_i/K_v}(b_i^v,d_i)$, there
exists $j \in I_1$ such that $${\rm cor}_{F^v_j/K_v}(a_j^v,d_j) \not = {\rm cor}_{F^v_j/K_v}(b_j^v,d_j).$$  Therefore $v \not \in  \underset{i\in I_1} \cap \Sigma_i$,
and this contradicts $\Sigma (L/K) \cup \underset{i\in I_0} \cap \Sigma_i \cup  \underset{i\in I_1} \cap \Sigma_i = \Omega_K$. Hence we have $f_{(a)} = f_{(b)}$
for all $(a), (b) \in {\cal L}(E,A)$.

\bigskip
Since $f_{(a)}$ is independent of $(a)$, we obtain a map

$$f : \sha(E,\sigma) \to {\bf Z}/2{\bf Z}$$ defined by $$f(I_0,I_1)  = \Sigma_{i \in I_0} \Sigma_{v \in \Omega_K} {\rm cor}_{F^v_i/K_v}(a_i^v,d_i)$$
for any  $(a) = (a_i^v) \in {\cal L}(E,A)$. Note that $f$ is a group homomorphism. 

\bigskip
Recall that we are fixing an embedding $\epsilon : E \to A$, and an involution $\theta : A \to A$ such that $\epsilon : (E,\sigma) \to (A,\tau)$ is
an embedding of algebras with involution. If $(A,\theta)$ is orthogonal, then we also fix an orientation $u : \Delta (E) \to Z(A,\theta)$. Our next aim
is to discuss the dependence of $f$ on these choices. We first introduce some notation.

\medskip Recall that for all $a \in F^{\times}$, we set $\theta_a = \theta \circ {\rm Int}(\epsilon (a))$. Similarly, if 
$\tilde \theta : A \to A$  is an involution and if $\tilde \epsilon : (E,\sigma) \to (A,\theta)$ is an embedding of algebras with involution, then
we set $\tilde \theta_a = \tilde  \theta \circ {\rm Int}(\tilde \epsilon (a))$. Then $\tilde \theta_a : A \to A$ is an involution, and $\tilde \epsilon : (E,\sigma) \to (A,\tilde \theta)$
is an embedding of algebras with involution.

\medskip
\noindent
{\bf Definition 5.4.2.} Let $\tilde \epsilon : E \to A$ be an embedding, and let $\tilde \theta : A \to A$ be an involution such that $\tilde \epsilon : (E,\sigma) \to (A,\tilde \theta)$ is
an embedding of algebras with involution. Let $\tilde u : \Delta (E) \to Z(A,\tilde \theta)$ be an orientation. We say that $(\epsilon,\theta,u)$ and $(\tilde \epsilon,
\tilde \theta, \tilde u)$ are {\it compatible} if there exists $\alpha \in A^{\times}$ and $c \in F^{\times}$ such that the following two conditions are satisfied 

\medskip

{\rm (a)} ${\rm Int}(\alpha) : (A,\tilde \theta) \to (A,\theta_c)$   is an isomorphism of algebras with involution such that ${\rm Int}(\alpha)  \circ   \tilde \epsilon = \epsilon$.

\medskip

{\rm (b)} The induced automorphism  $c(\alpha) :  Z(A,\tilde \theta)  \to Z(A,\theta_c)$ satisfies $$c(\alpha) \circ   \tilde u = u_c.$$

Recall that  if $(A,\tau)$ is orthogonal, then we are fixing an orientation $\nu : \Delta (E) \to Z(A,\tau)$. 

\medskip
\noindent
{\bf Proposition 5.4.3.} {\it Assume that $(\epsilon,\theta,u)$ and $(\tilde \epsilon, \tilde \theta, \tilde u)$ are compatible.  Let $\tilde {\cal L}(A,E)$ be the set of local embedding
data defined with respect to $(\tilde \epsilon, \tilde \theta, \tilde u)$, and let $(a) \in \tilde {\cal L}(A,E)$. Let 
$$f'_{(a)} : \sha (E,\sigma) \to {\bf Z}/2{\bf Z}$$ be defined by
$$f'_{(a)} (I_0,I_1) = \Sigma_{i \in I_0} \Sigma_{v \in \Omega_K} {\rm cor}_{F^v_i/K_v}(a_i^v,d_i).$$
Then $f_{(a)}' = f$.}

\medskip
\noindent
{\bf Proof.} Let $\alpha \in A^{\times}$ and $c \in F^{\times}$ be such that ${\rm Int}(\alpha) : (A,\tilde \theta) \to (A,\theta_c)$   is an isomorphism of algebras with involution such that ${\rm Int}(\alpha)  \circ   \tilde \epsilon = \epsilon$, and that if $\theta$ is orthogonal, then we have $c(\alpha) \circ   \tilde u = u_c.$

\medskip
Let $(a) = (a^v)  \in  \tilde {\cal L}(A,E)$.
We claim that  $(ca) \in {\cal L}(E,A)$. A straightforward computation shows that ${\rm Int}(\alpha^{-1}) : (A,\theta_{c a^v}) \to (A,\tilde \theta_{a^v})$ is an
isomorphism of algebras with involution for all $v \in \Omega_K$.

\medskip

For all $v \in \Omega_K$, let $({\rm Int} (\beta^v) \circ \tilde \epsilon, a^v,\beta^v,\nu)$ be parameters of an oriented embedding.
Since $\tilde \epsilon = {\rm Int}(\alpha ^{-1}) \circ \epsilon$ and $c(\alpha) \circ   \tilde u_a = u_{ca}$, we see that 
$({\rm Int} (\beta^v \alpha^{-1}) \circ \epsilon c a^v ,\beta^v \alpha^{-1}, \nu)$ are  parameters of an oriented embedding with respect to $(\epsilon,\theta,u)$. 
Therefore we have  $(ca) \in {\cal L}(E,A)$. 

\medskip

Let $c = (c_1,\dots,c_m)$ with
$c_i \in F_i^{\times}$. We have
$$f'_{(a)} (I_0,I_1) = \Sigma_{i \in I_0} \Sigma_{v \in \Omega_K} {\rm cor}_{F^v_i/K_v}(a_i^v,d_i) = $$
$$=  \Sigma_{i \in I_0} \Sigma_{v \in \Omega_K} {\rm cor}_{F^v_i/K_v}(a_i^v,d_i) + \Sigma_{i \in I_0} \Sigma_{v \in \Omega_K} {\rm cor}_{F^v_i/K_v}(c_i,d_i) = $$
$$= \Sigma_{i \in I_0} \Sigma_{v \in \Omega_K} {\rm cor}_{F^v_i/K_v}(c_ia_i^v,d_i)
= f(I_0,I_1),$$  since $(ca) \in {\cal L}(E,A)$.

\medskip
\noindent
{\bf Corollary 5.4.4.} {\it Suppose that there exists an embedding of algebras with involution $(E,\sigma) \to (A,\tau)$. Then 
we have $f = 0$.}

\medskip
\noindent
{\bf Proof.} Since there exists an embedding $(E,\sigma) \to (A,\tau)$, there exists $a \in F^{\times}$ such that $\tau \simeq \theta_a$. We have
$a = (a_1,\dots,a_m)$ with $a_i \in F^{\times}_i$. For all $v \in \Omega_K$,
set $a^v_i = a_i$, and let $(a) = (a^v_i)$. By Theorem 5.4.1.
it suffices to show that $f_{(a)} = 0$. Let $(I_0,I_1) \in \sha(E,\sigma)$. Then we have $$f_{(a)} (I_0,I_1) = \Sigma_{v \in \Omega_K} \Sigma_{i \in I_0}  {\rm cor}_{F^v_i/K_v}(a_i,d_i) = \Sigma_{v \in \Omega_K} \Sigma_{i \in I_0}  {\rm cor}_{F_i/K}(a_i,d_i) = 
0.$$ Therefore $f = f_{(a)} = 0$, as claimed.

\bigskip
{\bf 5.5. Hasse principle}

\medskip
The main result of the paper is the following :

\medskip
\noindent
{\bf Theorem 5.5.1.} {\it   Let $\nu : \Delta(E) \to Z(A,\tau)$ be an orientation. Suppose that for all $v \in \Omega_K$ there exists an oriented  embedding
$(E^v,\sigma) \to (A^v,\tau)$ with respect to $\nu$. 
Then 
there exists an embedding $(E,\sigma) \to (A,\tau)$ if and only if $f  = 0$.}

\medskip
This will be proved in Sections 6--8.

\bigskip
{\bf \S 6. Orthogonal involutions}

\medskip
Suppose that $K$ is a global field,  that $(A,\tau)$ is orthogonal, and that all the factors of $E$ split $A$. The
aim of this section is to give a criterion for the Hasse principle for the existence of an embedding of $(E,\sigma)$ into $(A,\tau)$ :  in other
words, to prove Theorem 5.5.1. for orthogonal involutions.
Moreover, based on the results of \S 2, we give necessary and
sufficient conditions for such an embedding to exist everywhere locally.

\medskip

{\bf  6.1. The even degree case  -- Hasse principle}

\medskip
Suppose that  ${\rm deg}(A) = n = 2r$. 
We fix an embedding $\epsilon : (E,\sigma) \to (A,\theta)$
and an isomorphism of $K$--algebras $u : \Delta(E) \to Z(A,\theta)$.

\medskip
Let us assume that for all $v \in \Omega_K$, there exists an embedding of algebras with involution
$(E^v,\sigma) \to (A^v,\tau)$. This implies that the $K$--algebras $\Delta(E)$ and $Z(A,\tau)$ are
isomorphic. 
Let us fix an isomorphism of $K$--algebras $\nu : \Delta (E) \to Z(A,\tau)$.

\medskip
The Brauer--Manin map
$f : \sha(E,\sigma)  \to {\bf Z}/2{\bf Z}$  was defined in 5.4.

\bigskip
\noindent
{\bf Theorem 6.1.1.} {\it   Suppose that for all $v \in \Omega_K$ there exists an oriented  embedding
$(E^v,\sigma) \to (A^v,\tau)$  with respect to $\nu$. 
Then 
there exists an embedding $(E,\sigma) \to (A,\tau)$ if and only if $f  = 0$.}

\medskip
\noindent
{\bf Proof.} By Corollary 5.4.4. we already know that the existence of a global embedding $(E,\sigma) \to (A,\tau)$  implies that $f = 0$. Let us prove the
converse. Let $(a) = (a_i^v )\in {\cal L}(E,A)$, and 
let $(I_0,I_1) \in \sha$. Then by hypothesis we have  $f(I_0,I_1)  = f_{(a)} (I_0,I_1) = 0$, hence 

$$\Sigma_{v \in \Omega_K} \Sigma_{i \in I_0}  \ {\rm cor}_{F^v_i/K_v} (a^v_i,d_i) = 0.$$
By Proposition  4.1.5. there exists $b \in F^{\times}$ such that  $${\rm cor}_{F^v/K_v}(b,d) = {\rm cor}_{F^v/K_v}(a^v,d)$$ for all $v \in \Omega_K$.
Applying Lemma 2.5.4. we see that $C(A^v,\theta_{a^v}) = C(A^v,\theta_b)$ in  ${\rm Br}(\Delta(E^v))$ for all $v \in \Omega_K$. Since the
embedding is oriented with respect to $\nu$, we have $C(A^v,\tau) = C(A^v,\theta_{a^v})$ in  ${\rm Br}(\Delta(E^v))$ for all $v \in \Omega_K$. 
Therefore for all $v \in \Omega_K$, we have $C(A^v,\tau) = C(A^v,\theta_b)$ in  ${\rm Br}(\Delta(E^v))$.
Then by the Brauer--Hasse--Noether Theorem, we have $C(A,\tau) = C(A,\theta_b)$ in ${\rm Br}(\Delta(E))$, hence $C(A,\tau)$ and $C(A,\theta_b)$ are
isomorphic over $K$. 
Note that $(A^v,\tau) \simeq (A^v,\theta_b)$ over $K_v$ if $v$ is a real place. Hence by [LT 99], Theorems A and B, we conclude that
$(A,\tau) \simeq (A,\theta_b)$. By  Proposition 1.1.3. there exists an embedding of $(E,\sigma)$ into $(A,\tau)$. 

\medskip
\noindent
{\bf Corollary  6.1.2.} {\it  Assume that for all $v \in \Omega_K$ there exists an embedding
$(E^v,\sigma) \to (A^v,\tau)$. Suppose moreover that one of the following holds :

\medskip
{\rm (i)} $A$ is split.

\smallskip
{\rm (ii)} If $v \in \Omega_K$ is such that $A^v$ is non-split, then ${\rm disc}(A^v,\tau) \not = 1$ in $K_v^{\times}/K_v^{\times 2}$.

\medskip
{\rm (iii)} ${\rm deg}(A) = 2r$ with $r$ odd.

\medskip

Then 
there exists an embedding $(E,\sigma) \to (A,\tau)$ if and only if $f  = 0$.}

\medskip
\noindent
{\bf Proof.} This follows from Theorem 6.1.1. together with Corollary 2.7.3. (in cases {\rm (i)} and {\rm (ii)}), and Corollary 2.8.3.  (in case {\rm (iii)}).

\bigskip
{\bf 6.2. The odd degree case -  Hasse principle}

\medskip
Suppose that $A = M_n(K)$, and that $\tau$ is induced by an $n$--dimensional  quadratic form $q$. Recall that $f  : \sha \to {\bf Z}/2{\bf Z}$ was defined in 4.4.

\medskip
\noindent
{\bf Theorem 6.2.1.} {\it Suppose that $n$ is odd, and that for all $v \in \Omega_K$ there exists an embedding of algebras with involution  $(E^v,\sigma) \to (A^v,\tau)$.
Then there exists an embedding $(E,\sigma) \to (A,\tau)$ if and only if $f = 0$.}

\medskip
\noindent
{\bf Proof.} We already know that if there exists  an embedding $(E,\sigma) \to (A,\tau)$, then we have  $f = 0$ (cf. Corollary 5.4.3). Let us show that the
converse also holds. Assume that we have $f = 0$.

\medskip

If $n = 1$ then $E = A = K$, hence  $(E,\sigma)$ can  be embedded into $(A,\tau)$. Let us assume that $n \ge 3$. 
Set $A' = M_{n-1}(K)$. Then by [PR 10], Proposition 7.2. there exists a $\sigma$--invariant  \'etale subalgebra $E'$ of rank $n-1$ of $E$ such that $E = E' \times K$,
an $(n-1)$--dimensional quadratic form $q'$ and a 1--dimensional quadratic form $q''$ over $K$  such that $q \simeq q' \oplus q''$, and that the \'etale algebra with involution $(E',\sigma)$ can be embedded
in the central simple algebra $(A',\tau')$ over $K_v$ for all $v \in \Omega_K$, where $\tau' : A' \to A'$ is the involution induced by $q'$. Moreover, there exists an embedding of $(E,\sigma)$ into $(A,\tau)$ if and only if there exists an embedding of $(E',\sigma)$ into $(A',\tau')$. Note that we have ${\cal L}(E,A) = {\cal L}' (E',A') \times {\cal L}(K,K)$. We may suppose that $E_m = K$. Then we have
$d_m = 1$. Set $J = \{1,\dots,m-1 \}$.

\medskip Let $f'  : \sha(E',\sigma) \to {\bf Z}/2{\bf Z}$ be the Brauer--Manin map associated to $(E',\sigma)$ and $(A',\tau')$. Let $(a) = (a^v_i) \in {\cal L}(E,A)$. 
Set $b^v_i = a^v_i$ if $i = 1,\dots,m-1$. Then $(b) = (b^v_i)$ is an element of ${\cal L}' (E',A')$. By Proposition 5.2.2. {\rm (i)} we have 
${\cal L}' (E',A') \subset {\cal L} (E',A')$, hence $(b) \in {\cal L} (E',A')$.

\medskip For all $(J_0,J_1) \in \sha (E',A')$ we have $f'(J_0,J_1) = f'_{(b)}(J_0,J_1) = f_{(a)}(I_0,I_1)$, where $I_0 = J_0$ and $I_1 = J_1 \cup \{m \}$.
Since $f_{a} = f = 0$ by hypothesis, this implies that $f' = 0$. By Corollary 6.1.2. {\rm (i)} this implies that $(E',\sigma)$ can be embedded into $(A',\tau')$.
Therefore $(E,\sigma)$ can be embedded into $(A,\tau)$.

\bigskip
{\bf  6.3. Orthogonal involutions -- local conditions}

\medskip

An infinite place $w$ of $F$ is said to be {\it ramified in $E$} if $w$ is a real place that extends to a complex place of $E$.  For all $v \in \Omega_K$, let $\rho_v$ be the number
of places of $F$ above $v$ which are not ramified. 

\medskip
\noindent
{\bf Definition 6.3.1.} We say that the {\it signature conditions} hold if
 for every real prime $v$ of $K$ such that $A^v \simeq M_n(K_v)$, the signature of $q$ at $v$ is of the shape $(r_v + \rho_v, s_v+\rho_v)$ for
some non--negative integers $r_v$ and $s_v$ such that  $r_v$ and $s_v$ are even if $n$ is even.

\medskip
\noindent
{\bf Definition 6.3.2.} We say that the {\it hyperbolicity condition} is satisfied if for all $v \in \Omega_K$ such that the \'etale algebra with involution $(E^v,\sigma)$ is split, the algebra with
involution $(A^v,\tau)$ is hyperbolic.

\medskip
Note that as $(A^v,\tau)$ is hyperbolic for all but a finite number of places $v \in \Omega$, we only need to check the hyperbolicity condition at  finitely many places.

\medskip
The following is a consequence of the results of \S 3, in particular  Propositions 3.1.1. and 3.1.2. (see also [B 12], Proposition 2.4.1. and [B 13], Theorem 12.1.).

\bigskip
\noindent
{\bf Proposition 6.3.3.} {\it Suppose that $n$ is even. The \'etale algebra with involution $(E^v,\sigma)$ can be embedded in the algebra with involution $(A^v,\tau)$ for
all $v \in \Omega_K$ if and only if the following conditions hold :

\medskip
{\rm (i)} We have ${\rm disc}(A,\tau) = {\rm disc}(E) \in K^{\times}/K^{\times 2}$.

\smallskip 
{\rm (ii)} The hyperbolicity condition is
satisfied.

\smallskip
{\rm (iii)} The signature conditions are satisfied.}

\medskip Assume now that $n$ is odd. 
We have $E \simeq E' \times K$, where $E'$ is an \'etale algebra of rank $n-1$ stable by $\sigma$. Note that if $n = 1$, then an embedding of
$(E,\sigma)$ into $(A,\tau)$ always exists, hence we may assume that $n \ge 3$. The following is a consequence of Propositions 3.2.1. and 3.2.2. 

\medskip
\noindent
{\bf Proposition 6.3.4.} {\it Suppose that $n$ is odd and $\ge 3$. The  \'etale algebra with involution $(E^v,\sigma)$ can be embedded in the algebra with involution $(A^v,\tau)$ for
all $v \in \Omega_K$ if and only if the following conditions hold :

\medskip
{\rm (i)} For all $v \in \Omega_K$ such that $(E')^v,\sigma)$ is split, we have $q \simeq q' \oplus q''$, where $q'$ is hyperbolic, and $q''$ is a $1$--dimensional
quadratic form over $K_v$. 

\smallskip
{\rm (ii)} The signature conditions are satisfied.}

\bigskip
{\bf \S 7. Symplectic involutions}

\medskip
\medskip
Suppose that $K$ is a global field,  that all the factors of $E$ split $A$, 
and that $(A,\tau)$ is a symplectic involution. 
Prasad and Rapinchuk proved that the Hasse principle holds in this case (cf.  [PR 10], Theorem 5.1).

\medskip
An infinite place $w$ of $F$ is said to be {\it ramified in $E$} if $w$ is a real place that extends to a complex place of $E$.  For all $v \in \Omega_K$, let $\rho_v$ be the number
of places of $F$ above $v$ which are not ramified. 

\medskip
\noindent
{\bf Definition 7.1.1.} We say that the {\it signature condition} holds if
 for every real prime $v$ of $K$ such that $A^v$ is non--split, the signature of  $(A^v,\tau)$  is of the shape $(r_v + \rho_v, s_v+\rho_v)$ for
some non--negative integers $r_v$ and $s_v$.

\medskip The following result is a consequence of the Hasse principle, and of Proposition 3.3.1.

\medskip
\noindent
{\bf Theorem 7.1.2.} {\it 
The \'etale algebra with involution $(E,\sigma)$ can be embedded in the central simple algebra with involution $(A,\tau)$ if and only if the signature condition holds.}

\bigskip

{\bf \S 8. Unitary involutions}

\medskip
Suppose that $K$ is a global field,  that all the factors of $E$ split $A$, 
and that $(A,\tau)$ is a unitary involution. 

\bigskip
{\bf 8.1. Unitary involutions -- Hasse principle}

\bigskip
Suppose that $(E^v,\sigma)$ can be embedded in $(A^v,\tau)$ for all $v \in \Omega_K$, and recall that 
$f : \sha(E,\sigma)  \to {\bf Z}/2{\bf Z}$ is the Brauer--Manin map defined in 5.4. 

\medskip
\noindent
{\bf Theorem 8.1.1.} {\it Suppose that for all $v \in \Omega_K$ there exists an embedding of algebras with involution $(E^v,\sigma) \to (A^v,\tau)$. Then there exists an embedding  of algebras with involution $(E,\sigma) \to (A,\tau)$ if and only if $f = 0$.}

\medskip
\noindent
{\bf Proof.} If there exists a global embedding, then we have $f = 0$ (cf. Corollary 5.4.4.). Let us prove the converse.
For all $v \in \Omega_K$ there exist $a^v \in (F^v)^{\times}$ such that $(A^v,\tau) = (A^v,\theta_{a^v})$, hence we have  ${\rm D}(A^v,\tau) =
 {\rm D}(A^v,\theta_{a^v})$ in ${\rm Br}(K_v)$. On the other hand, by 
Lemma 1.5.2. we
have ${\rm D}(A^v,\theta_{a^v}) = {\rm D}(A^v,\theta) + {\rm cor}_{F^v/K_v}(a^v,d).$ Hence we have
$\Sigma_{v \in \Omega_K} {\rm cor}_{F^v/K_v} (a^v,d) =0.$
Let $(I_0,I_1) \in \sha$. Then by hypothesis we have  $f_{(a)} (I_0,I_1) = 0$, therefore 
$$\Sigma_{v \in \Omega_K} \Sigma_{i \in I_0}  \ {\rm cor}_{F^v_i/K_v} (a^v_i,d_i) = 0.$$
By Proposition 4.1.5.  there exists $b \in F^{\times}$ such that for all $v \in \Omega_K$, we have ${ \rm cor}_{F^v/K_v} (b,d)  = {\rm cor}_{F^v/K_v} (a^v,d),$
and that  $b^v = a^v$  if $v$ is a real place. On the other hand, we have
${\rm D}(A,\tau) =  {\rm D}(A,\theta_b)$. Hence we have $(A,\tau) \simeq (A,\theta_b)$. By  Proposition 1.1.3. there exists an embedding of $(E,\sigma)$ into $(A,\tau)$. 

\bigskip
{\bf 8.2. Unitary involutions -- local conditions}

\medskip

An infinite place $w$ of $F$ is said to be {\it ramified in $E$} if $w$ is a real place that extends to a complex place of $E$.  For all $v \in \Omega_K$, let $\rho_v$ be the number
of places of $F$ above $v$ which are not ramified. 

\medskip
\noindent
{\bf Definition 8.2.1.} We say that the {\it signature condition} holds if
 for every real prime $v$ of $K$, the signature of  $(A^v,\tau)$  is of the shape $(r_v + \rho_v, s_v+\rho_v)$ for
some non--negative integers $r_v$ and $s_v$.

\medskip
\noindent
{\bf Definition 8.2.2.} We say that the {\it hyperbolicity condition} is satisfied if for all $v \in \Omega_K$ such that the \'etale algebra with involution $(E^v,\sigma)$ is split, the algebra with
involution $(A^v,\tau)$ is hyperbolic.

\medskip
The following is a consequence of Propositions 3.4.1. and 3.4.2. 

\medskip
\noindent
{\bf Proposition 8.2.3.} {\it The \'etale algebra with involution $(E^v,\sigma)$ can be embedded in the algebra with involution $(A^v,\tau)$ for
all $v \in \Omega_K$ if and only if the following conditions hold :

\medskip
{\rm (i)} For all $v \in \Omega_K$, we have $${\rm det}(A^v,\tau) {\rm disc}(E^v,\sigma)^{-1} \in {\rm N}_{F^v/K_v}((F^v)^{\times}) 
 {\rm N}_{L^v/K_v}((L^v)^{\times}).$$
 
 \smallskip {\rm (iii)} The signature condition is satisfied. 

\smallskip 
{\rm (iii)} The hyperbolicity condition is
satisfied.}

\bigskip
{\bf 9. Applications and examples}

\medskip
The aim of this section is to describe some special cases in which the Hasse principle for the embedding problem holds, and to give some examples.
We keep the notation of the previous sections. In particular, $K$ is a global field, $(E,\sigma)$ is an \'etale algebra with involution, and
$(A,\tau)$ is a central simple algebra with involution. 

\medskip
{\bf 9.1. The group $\sha(E',\sigma)$}

\medskip
Let us write $E = E_1 \times \dots \times E_{m_1} \times E_{m_1 + 1} \times \dots \times E_m$, where $E_i/F_i$ is a quadratic extension
for all $i = 1, \dots, m_1$ and $E_i = F_i \times F_i$ or $E_i = K$ if $i = m_1 + 1, \dots, m$. Recall that $I = \{1,\dots,m \}$, and set $I({\rm split}) = \{m_1+1, \dots, m \}$,
$I' = I({\rm nonsplit}) =  \{1,\dots,m_1 \}$. If $I'$ is empty, then we set $\sha(E',\sigma) = 0$.

\medskip
Let $\pi : \sha (E,\sigma) \to \sha (E',\sigma)$ be the map that sends the class of $(I_0,I_1)$ to the class of $(I_0 \cap I',I_1 \cap I')$. Then $\pi$ is surjective,
and ${\rm Ker}(\pi)$ is the subgroup of $\sha (E,\sigma)$ consisting of the classes of partitions $(I_0,I_1)$ such that $I_0 \subset I({\rm split})$ or
$I_1 \subset I({\rm split})$. 

\medskip
Let $f : \sha(E,\sigma) \to {\bf Z}/2 {\bf Z}$ be the Brauer--Manin map (cf. \S 5.4.). Note that ${\rm Ker}(\pi) \subset {\rm Ker}(f)$, since if $i \in I({\rm split})$, then $d_i = 1$.
Hence $f$  induces a map $\overline f : \sha (E',\sigma) \to {\bf Z}/2{\bf Z}$
such that $f = \overline f \circ \pi$. 

\medskip
\noindent
{\bf Proposition 9.1.1.} {\it We have $f = 0$ if and only if $\overline f = 0$.}

\medskip
\noindent
{\bf Proof.} This follows immediately from the definitions.

\medskip
{\bf 9.2. Sufficient conditions}

\medskip

Assume that for all $v \in \Omega_K$ there exists an embedding
$(E^v,\sigma) \to (A^v,\tau)$, and let $\nu : \Delta(E) \to Z(A,\tau)$ be an orientation. The results of Sections 6--8 imply the following :

\medskip
\noindent
{\bf Theorem 9.2.1.} {\it Suppose that the following conditions hold : 

\medskip

{\rm (i)} For all $v \in \Omega_K$, there exists an oriented embedding  $(E^v,\sigma) \to (A^v,\tau)$ with respect to $\nu$.

\smallskip

{\rm (ii)}   $\sha(E',\sigma)$ is trivial.

\medskip
Then there exists an embedding  $(E,\sigma) \to (A,\tau)$.}

\medskip
\noindent
{\bf Proof.} This follows from Theorems 6.1.1. 6.2.1.  7.1. 8.1.1. and Proposition 9.1.1. 

\medskip Note that the existence of an {\it oriented} embedding is only necessary if $(A,\tau)$ is orthogonal, $A$ is non-split and ${\rm deg}(A) = 2r$ with $r$ even (cf. Corollary 6.1.2.). Note also that this
implies Theorem A of Prasad and Rapinchuk (cf. [PR 10], page 584) -- indeed, if $E$ is a field extension of $L$, then $\sha (E,\sigma)$ (= $\sha(E',\sigma)$ in
this case) is obviously trivial.
Theorem 9.2.1. also has the following application :

\medskip
\noindent
{\bf Corollary 9.2.2.} {\it Suppose that the following conditions hold :

\medskip

{\rm (i)} For all $v \in \Omega_K$, there exists an oriented embedding  $(E^v,\sigma) \to (A^v,\tau)$ with respect to $\nu$.

\smallskip

{\rm (ii)}   There exists $i_0 \in I$ such that  for all $i \in I$, we have $$\Sigma(L/K) \cup \Sigma_{i_0} \cup \Sigma_i \not = \Omega_K.$$

Then there exists an embedding  $(E,\sigma) \to (A,\tau)$.}

\bigskip
This generalizes the Hasse principle results of  [PR 10], [Lee 12] and [B 12]. The Corollary is a consequence of Theorem 9.2.1. and the following Lemma :

\medskip
\noindent
{\bf Lemma 9.2.3.} {\it Assume that there exists $i_0 \in I$ such that  for all $i \in I$, we have $\Sigma(L/L) \cup \Sigma_{i_0} \cup \Sigma_i \not = \Omega_K$. 
Then the group $\sha(E,\sigma)$ is trivial. Therefore $\sha(E',\sigma)$ is trivial.}

\medskip
\noindent
{\bf Proof.} 
Suppose that the group $\sha(E,\sigma)$ is not trivial, and let $(I_0,I_1)$ be a partition of $I$ representing a non--trivial element
of $\sha(E,\sigma)$. Then we have 
$$\Sigma(L/K) \cup 
(\underset{i\in I_0}{\cap}\Sigma_i)\cup(\underset{j\in I_1}{\cap}\Sigma_j)=\Omega_K.$$ 
Assume that $i_0 \in I_0$. Then we have $\Sigma(L/K) \cup 
\Sigma_{i_0}\cup(\underset{j\in I_1}{\cap}\Sigma_j)=\Omega_K,$ hence for all $j \in I_1$, we have $\Sigma(L/K) \cup 
\Sigma_{i_0}\cup \Sigma_j=\Omega_K$, contradicting the hypothesis.

\medskip
\noindent
{\bf Corollary 9.2.4.}  {\it Suppose that the following conditions hold :

\medskip

{\rm (i)} For all $v \in \Omega_K$, there exists an oriented embedding  $(E^v,\sigma) \to (A^v,\tau)$ with respect to $\nu$.

\smallskip

{\rm (ii)}   There exists a real place $u \in \Omega_K$ such that $u \not \in  \Sigma_i$ for all $i \in I$.

\medskip
Then there exists an embedding  $(E,\sigma)  \to (A,\tau)$.}

\medskip
\noindent
{\bf Proof.} By  {\rm (ii)}, 
condition {\rm (ii)} of Corollary 9.2.2. holds, hence there exists an embedding $(E,\sigma) \to (A,\tau)$.

\medskip
Assume now that $K = {\bf Q}$. Recall that $(E,\sigma)$ is a CM \'etale algebra if $E$ is a product of CM fields, and if $\sigma$ is the complex conjugation. 
Then we have

\medskip
\noindent
{\bf Corollary 9.2.5.}  {\it Suppose $K = {\bf Q}$, and that $(E,\sigma)$ is a CM \'etale algebra. Assume that for all $v \in \Omega_K$, there exists an oriented embedding  $(E^v,\sigma) \to (A^v,\tau)$. Then there exists an embedding  $(E,\sigma)  \to (A,\tau)$.}

\medskip
\noindent
{\bf Proof.} This follows from Corollary 9.2.4. since condition {\rm (ii)} holds for CM \'etale algebras.

\bigskip
{\bf 9.3. An example}

\medskip
As we have seen in Corollary 9.2.5. above, the local--global principle holds for oriented embeddings when $(E,\sigma)$ is a CM  \'etale algebra with involution. The
aim of this section is to show that this is not the case for not necessarily orientated local embeddings. More precisely, there exist CM \'etale
algebras with involution $(E,\sigma)$ and (non--split) central simple algebras with orthogonal involution $(A,\tau)$ such that $(E,\sigma)$ embeds into $(A,\tau)$ everywhere
locally, but not globally. 

\medskip
Let $v_1,v_2,v_3$ and $v_4$ be four distinct places of $K$. Let $a \in K^{\times}$ be such that
$a \not \in K_{v_i}^{\times 2}$ for $i = 1,\dots,4$, and let $b \in K^{\times}$ such that $b \not \in K^{\times 2}$ and that $b \in K_{v_i}^{\times 2}$ for $i = 1,\dots,4$.
Let $E_1 = K(\sqrt a)$, and let $\sigma_1 : E_1 \to E_1$ be the $K$--linear involution such that $\sigma_1(\sqrt a) = - \sqrt a$. Set 
$E_2 = K(\sqrt b)$,  let $\sigma_2 : E_2 \to E_2$ be the $K$--linear involution such that $\sigma_2(\sqrt b) = - \sqrt b$. Set
$E = E_1 \otimes E_2$ and $\sigma = \sigma_1 \otimes \sigma_2$. Then $(E,\sigma)$ is a rank 4  \'etale $K$--algebra with involution, and
$F = E^{\sigma} = K(\sqrt {ab})$. 

\medskip
Let $H_1$ be the quaternion skew field over $K$ ramified exactly at $v_1$ and $v_2$, and $H_2$ the quaternion skew field over $K$ ramified
exactly at $v_3$ and $v_4$. Let $\tau_i : H_i \to H_i$ be the canonical involution for $i = 1,2$, and set $(A,\tau) = (H_1,\tau_1) \otimes (H_2,\tau_2)$.
Since $\tau_1$ and $\tau_2$ are both symplectic involutions, their tensor product $\tau$ is an orthogonal involution. We have $H_1 \otimes H_2 \simeq
M_2(H)$, where $H$ is a quaternion skew field over $K$.

\medskip
\noindent
{\bf Proposition 9.3.1.} {\it For all $v \in \Omega_K$, there exists an embedding of algebras with involution $(E^v,\sigma) \to (A^v,\tau)$}.

\medskip
\noindent
{\bf Proof.} Since $E_1$ splits $H_1$ and $H_2$ locally everywhere, it splits $H$ locally everywhere too, and hence $E$ embeds in $H$ as a maximal
subfield globally. Let $\tau_0$ be the canonical involution of $H$. Since $\tau_0$ restricts to the non--trivial automorphism on any maximal subfield, it
follows that there exists an embedding of algebras with involution of $(E_1,\sigma_1)$ into $(H,\tau_0)$. 

\medskip

Let $w \in \Omega_K$. By hypothesis, either $H_1$ or $H_2$ is split over $K_w$. Hence either $(H_1^w,\tau_1) \simeq (M_2(K_w),\sigma_0)$ or
$(H_2^w,\tau_2) \simeq (M_2(K_w),\sigma_0)$, where $\sigma_0$ denotes the symplectic involution of $M_2(K_w)$. Therefore we have
$$(M_2(H^w),\tau) \simeq (H_1^w \otimes H_2^w,\tau_1 \otimes \tau_2) \simeq (M_2(K_w),\sigma_0) \otimes (H^w,\tau_0).$$
The algebra with involution $(E^w_1,\sigma_1)$ can be embedded into $(H^w,\tau_0)$, and the
algebra with involution $(E^w_2,\sigma_2)$ can be embedded into $(M_2(K_w),\sigma_0)$. 
Hence $(E^w,\sigma)$ embeds into $(M_2(H^w),\tau)$.

\medskip
\noindent
{\bf Proposition 9.3.2.} {\it There is no global embedding  $(E,\sigma) \to (A,\tau)$}. 

\medskip
\noindent
{\bf Proof.} Let us denote by $H_i^0$ the skew elements of $H_1$, for $i = 1,2$. Then every skew element of $H_1 \otimes H_2$ belongs to the direct
sum $H_1^0 \oplus H_2^0$. Moreover, if a skew element is square central, then it has to be in $H_1^0$ or in $H_2^0$.
Assume by contradiction that there exists an embedding of algebras with involution $f : (E,\sigma) \to (A,\tau)$. 
Note that $f(\sqrt b)$ is a square central skew element. Therefore it has to belong to $H_1^0$ or to $H_2^0$. But this
contradicts the fact that $E_2 = K(\sqrt b)$ does not split $H_1$ nor $H_2$.

\bigskip
In the above example, we can take $K = {\bf Q}$ and can choose $a$ and $b$ such that $E$ is a $CM$ \'etale algebra. This provides the
desired counter--example to the Hasse principle.

\bigskip

\centerline  {\bf Appendix A}

\medskip

{\bf Embedding functor, Tate--Shafarevich group and orientation}

\medskip
The purpose of this appendix is to recall some of the results of [Lee 14], and to outline the relationship
of these results with those of the present paper.

\medskip

{\bf \S  A1. The embedding functor}

\medskip

Let $K$ be a field of characteristic $\not = 2$, let $K_s$ be a separable closure of $K$, and let
$\Gamma_K = {\rm Gal}(K_s/K)$. Let
$G$ be a reductive group over $K$. Let $T$ be a torus and let $\Psi$
be a root datum attached to $T$ (see [SGA 3], Exp. XXI, 1.1.1.). For a maximal torus $T'$ in
$G$, we let $\Phi(G,T')$ be the root datum of $G$ with
respect to $T'$. If $\Phi(G,T')_{K_s}$ and $\Psi_{K_s}$
are isomorphic, then we say that 
$G$ and $\Psi$ have the same type.

\medskip

Assume that $G$ and $\Psi$ have the same type. Let
$\underline {\rm Isom} (\Psi,\Phi(G,T'))$ be the scheme of isomorphisms
between the root data $\Psi$ and $\Phi(G,T')$. Define
$$\underline {\rm Isomext} (\Psi,\Phi(G,T')) = \underline {\rm Isom} (\Psi,\Phi(G,T'))/ {\rm W}(\Psi),$$
where ${\rm W}(\Psi)$ is the Weyl group of $\Psi$. The scheme $\underline {\rm Isomext} (\Psi,\Phi(G,T'))$
is independent of the choice of the maximal torus $T'$, and we denote it by $\underline {\rm Isomext} (\Psi, {G})$.
An {\it orientation} is by definition an element of $\underline {\rm Isomext} (\Psi, {G})(K)$.

\medskip

The \emph{embedding functor} $E(G,\Psi)$ is defined as follows:
for any $K$-algebra $C$, let $E(G,\Psi)(C)$ be the set of embeddings $f : T_{C} \to G_{C}$ such that
$f$ is both a closed immersion and a group homomorphism which induces an isomorphism
$f^{\Psi}:\Psi_{C}\xrightarrow{\sim}\Phi(G_{C},f(T_{C}))$ such that $f^{\Psi}(\alpha)=\alpha\circ f^{-1}|_{f(T_{C'})}$ for
all the $C'$-roots $\alpha$ in $\Psi_{C'}$ for each $C$-algebra $C'$ (see [Lee 14], 1.1.)
Given an orientation
$\nu \in\ul{\Isomext}(\Psi,G)(K)$, we define the \emph{oriented
embedding functor} as follows (cf. [Lee 14], 1.2.) : for any $K$-algebra $C$, set
\[E(G,\Psi,\nu)(C)=\left\{\begin{array}{l}\mbox{$f:T_{C}\hookrightarrow G_{C}$}\left|
\begin{array}{l}\mbox{$f\in E(G,\Psi)(C)$,
and the image of $f^{\Psi}$ } \\
\mbox{ in $\underline{\Isomext}(\Psi,G)(C)$ is
$\nu$.}\end{array}\right.\end{array}\right\}.\]

The oriented embedding functor is a homogeneous space under the
adjoint action of $G$. For each root datum $\Psi$, we can
associate a simply connected root datum $\sico(\Psi)$ to it (cf. 
[SGA3], Exp. XXI, 6.5.5 {\rm (iii)]). Let $\sico(T)$ be the torus
associated to $\sico(\Psi)$. 

\bigskip

{\bf \S A2. Algebras with involution and the embedding functor}

\medskip Let $L$ be a field of characteristic $\not = 2$, and let  $A$
be a central simple algebra over $L$ with involution $\tau$. Let
$E$ be an \'etale algebra over $L$ with involution $\sigma$, and suppose that $L^{\sigma} = K$. Given
$(A,\tau)$ and $(E,\sigma)$, we always assume that
\begin{center}
$\dim_L(E)=\mathrm{deg}_K(A)$ and $\tau|_L=\sigma|_L$.
\end{center}
The unitary groups $\rU(A,\tau)$ and $\rU(E,\sigma)$ are defined as follows. For
any  commutative $K$-algebra $C$, set 
$$\rU(A,\tau)(C)=\{x\in A\otimes_KC|\ x\tau(x)=1\},$$ and
$$\rU(E,\sigma)(C)=\{x\in E\otimes_KC|\ x\sigma(x)=1\}.$$ 
Let $G=\rU(A,\tau)^{\circ}$ be the connected component of
$\rU(A,\tau)$ containing  the neutral element, and let $T = \rU(E,\sigma)^{\circ}$
be the connected component of 
$\rU(E,\sigma)$ containing the neutral element.

\medskip Set $F = E^{\sigma}$. 
Let us  suppose furthermore that
\begin{center}
$\dim_K (F) =\left\{
         \begin{array}{ll}
             \lceil\frac{\dim_L (E)}{2}\rceil, & \hbox{if $\tau$ is of the first kind;} \\
             \dim_L (E), & \hbox{if $\tau$ is of the second kind.}
           \end{array}
         \right.$
\end{center}
Then one can associate a root datum $\Psi$ to the torus $T$ such that
$G$ is of type $\Psi$  (see [Lee 14], 1.3.).
Moreover, except for $A$ of degree 2 with $\tau$ orthogonal, there
exists a $K$-embedding from $(E,\sigma)$ to $(A,\tau)$ if and only
if there exists an orientation $\nu$ such that $E(G,\Psi,\nu)(K)$ is
nonempty (see [Lee 14],  Theorem 1.15. and Proposition 1.17.).

\bigskip
{\bf \S A3. {Orientations in terms of algebras}}

\medskip

Let $(E,\sigma)$ and $(A,\tau)$ be as above. Assume moreover that $(A,\tau)$ is orthogonal, and that the degree of $A$ is even. Recall that $\Delta(E)$ is the discriminant of
the \'etale algebra $E$, and that $Z(A,\tau)$ is the center of the Clifford algebra of $(A,\tau)$.  In 1.8. an orientation is defined as the choice of
an isomorphism $\Delta(E) \to Z(A,\tau)$. This is equivalent to the definition of A 1. More precisely, we have

\medskip
\noindent
{\bf Proposition A.3.1.} {\it We have an isomorphism

$$\ul{\Isom}(\Delta(E),Z(A,\tau)) \simeq\ul{\Isomext}(\Psi,G).$$}

\noindent
{\bf Proof.}
Let $E_\tau$ be a
maximal $\tau$-invariant \'etale subalgebra of $A$. Let
$T_\tau=\rU(E_\tau,\tau)^{\circ}$; then $T_\tau$ is a
maximal torus of $G$. Let $\Phi(G,T_\tau)$ be the root datum of $G$
with respect to $T_\tau$. Then we have a natural map
$\alpha:\ul{\Isom}((E,\sigma),(E_\tau,\tau))\ra\ul{\Isom}(\Psi,\Phi(G,T_\tau))$.
Using the identification of $\ul{\Aut}(E,\sigma)$ and
$\ul\Aut ({\Psi})$, we see that $\alpha$ is equivariant under the action of
$\ul{\Aut}(E,\sigma)$. Let $\Gamma_0$ be the subgroup of
$\ul{\Aut}(E,\sigma)$ corresponding to 
the Weyl group of $\Psi$ under this identification.
Let us
consider the following commutative diagram:
\begin{center}
\xymatrix{
  \ul{\Isom}((E,\sigma),(E_\tau,\tau))\ar[d]  \ar[r]
                & \ul{\Isom}(\Psi,\Phi(G,T_\tau)) \ar[d] \\
  \ul{\Isom}((E,\sigma),(E_\tau,\tau))/\Gamma_0  \ar[r]
                & \ul{\Isom}(\Psi,\Phi(G,T_\tau))/\rW(\Psi). }
\end{center}
Recall that $\ul{\Isom}(\Psi,\Phi(G,T_\tau))/\rW(\Psi) = 
\ul{\Isomext}(\Psi,\Phi(G,T_\tau))$,  and note that we have 
$\ul{\Isom}((E,\sigma),(E_\tau,\tau))/\Gamma_0\simeq
\ul{\Isom}(\Delta(E),\Delta(E_\tau))$.

\medskip

If we pick another maximal \'etale subalgebra  $E_\tau'$ of $A$ invariant by $\tau$, then the method used for $\ul{\Isomext}(\Psi,\Psi_\tau)$ in  [Lee 14] 1.2.1.
shows that 
we have a canonical isomorphism
between $\ul{\Isom}(\Delta(E),\Delta(E'_\tau))$ and
$\ul{\Isom}(\Delta(E),\Delta(E_\tau))$.

\medskip

Let us fix an isomorphism $\Delta(E_\tau) \to
Z(A,\tau)$ as in 1.8. This gives an isomorphism
$\ul{\Isom}(\Delta(E),\Delta(E_\tau)) \to
\ul{\Isom}(\Delta(E),Z(A,\tau))$. Hence, we have
$$\ul{\Isom}(\Delta(E),Z(A,\tau))\simeq\ul{\Isomext}(\Psi,\Phi(G,T_\tau)) = \ul{\Isomext}(\Psi,G),$$
as claimed.

\bigskip

{\bf  \S A4. Tate--Shafarevich group}

\medskip
Assume now that $K$ is a global field. Then, using Borovoi's results
(cf. [Bo 99]), it is shown in [Lee 14] that 
the Brauer-Manin obstruction is the only obstruction to
the local-global principle for $E(G,\Psi,u)$ and the obstruction
lies in the Tate--Shafarevich group $\sha^2(K,\sico(T))$ (cf.[ Lee 14], Proposition 2.8). 
Note that $\sha^2(K,\sico(T))$ is isomorphic to $\sha^1(K,\sico(\hat T))^{\ast}$ by
Poitou-Tate duality (cf. [NSW 08], Chap. VIII, Thm. 8.6.9).

\medskip
In the following, we determine the group  $\sha^1(K,\sico(\hat T))$ explicitly, and
show that it is isomorphic to the group $\sha (E',\sigma)$ defined in \S 9 :

\medskip
\noindent
{\bf Proposition A.4.1.} {\it The groups  $\sha^1(K,\sico(\hat T))$  and $\sha (E',\sigma)$ are isomorphic.}

\medskip
The proof of this proposition is different according as $L = K$ or $L \not = K$. Let us start by introducing some notation
that will be used in both proofs. For any finite separable field extension $N/N'$ and any discrete $\Gamma_N$--module $M$, set
$\rI_{N/N'} (M) = {\rm Ind}_{\Gamma_N}^{\Gamma_{N'}}(M)$. Note that $\rI_{N/N'}({\bf Z})$ is the character group of $\rR_{N/N'}({\bf G}_m)$. Let 
$\hat\rS_{N/N'}$ be the
character group of the norm-one torus $\rR^{(1)}_{N/N'}({\bf G}_m)$.

\medskip
\noindent
{\bf Proof of Proposition A.4.1. when $L = K$}

\medskip

Let us consider the following diagram:

\begin{equation}
\xymatrix{
  & 1\ar[d] & 1 \ar[d]\\
  1 \ar[r] & \rR^{(1)}_{F/K}({\bf G}_m) \ar[r] \ar[d] & \rR^{(1)}_{E/K}({\bf G}_m) \ar[r]\ar[d] & \sico(T) \ar[r] &
  1\\
  1 \ar[r] & \rR_{F/K}({\bf G}_m) \ar[r] \ar[d] & \rR_{E/K}({\bf G}_m)\ar[d]
  \\
  1 \ar[r] & {\bf G}_m \ar[r]^{\times 2} \ar[d]& {\bf G}_m \ar[d]\\
  & 1  & 1
  }
\end{equation}
where the first row (cf. [Lee 14], Lemma 3.16.) and the columns are exact. Then consider the
corresponding  diagram of  character groups:
\begin{equation}\label{e2}
\xymatrix{
  && 0 \ar[d] & 0 \ar[d]\\
  & &  {\bf Z} \ar[r]^{\times 2}\ar[d] & {\bf Z}  \ar[d]\\
  && \rI_{E/K}({\bf Z})\ar[r]^{\pi} \ar[d]&
  \rI_{F/K}({\bf Z})\ar[r]\ar[d] & 0\\
  0 \ar[r] & \sico(\hat{T}) \ar[r] & \hat\rS_{E/K} \ar[r]^{\pi} \ar[d]&
  \hat\rS_{F/K} \ar[r] \ar[d]& 0\\
  && 0 & 0
  }
\end{equation}

Note that we have
$\rI_{E/K}({\bf Z})=\underset{i=1}{\overset{m}{\oplus}}\rI_{E_i/K}({\bf Z})$
and
$\rI_{F/K}({\bf Z})=\underset{i=1}{\overset{m}{\oplus}}\rI_{F_i/K}({\bf Z})$.
The module $\rI_{E_i/K}({\bf Z})$ can also be written as
$\rI_{F_i/K}(\rI_{E_i/F_i}({\bf Z}))$. Let $d$ be the degree map from
$\rI_{E_i/F_i}({\bf Z})\simeq {\bf Z} \oplus {\bf Z}$ to ${\bf Z}$, which sends
$(x,y)$ to $x+y$. Then on each $\rI_{F_i/K}(\rI_{E_i/F_i}({\bf Z}))$,
the map $\pi$ is the map induced by the degree map from
$\rI_{E_i/F_i}({\bf Z})$ to ${\bf Z}$.

\medskip

Set $\Gamma = \Gamma_K$. We derive the following long exact sequence from diagram (2) :

$$ 0 \to {\rm sc}(\hat T)^{\Gamma}  \to (\hat S_{E/K})^{\Gamma} \buildrel {\pi} \over \longrightarrow \ (\hat S_{F/K})^{\Gamma} \to H^1(K,{\rm sc}(\hat T))
\to H^1(K,\hat S_{E/K}).$$
Thus we have the exact sequence
$$0 \to (\hat S_{F/K})^{\Gamma} / \pi ((\hat S_{E/K})^{\Gamma}) \buildrel \delta \over  \longrightarrow  H^1(K,{\rm sc}(\hat T))  \to H^1(K,\hat S_{E/K} )$$

Note that 
$\rH^2(K,\rR^{(1)}_{E/K}({\bf G}_m))$ injects into
$\rH^2(K,\rR_{E/K}({\bf G}_m))$ by Hilbert's Theorem 90. 
By the Brauer-Hasse-Noether
Theorem, $\sha^2(K,\rR_{E/K}({\bf G}_m))$ vanishes, hence so does
$\sha^2(K,\rR^{(1)}_{E/K}({\bf G}_m))$. 
By Poitou-Tate duality, we
have
$$\sha^1(K,\hat\rS_{E/K})\simeq\sha^2(K,\rR^{(1)}_{E/K}({\bf G}_m)))^\ast=0.$$

Therefore, $\sha^1(K,\sico(\hat {T}))$ is in the image of $(\hat S_{F/K})^{\Gamma} / \pi ((\hat S_{E/K})^{\Gamma})$.

\medskip

Since the  $F_i's$ are field extensions of $K$,  we have $\rI_{F_i/K}({\bf Z})^{\Gamma} \simeq {\bf Z}.$ Thus, we have
$\rI_{F/K}({\bf Z})^{\Gamma} \simeq \overset{m} {\underset{i} \oplus} \rI_{F_i/K}({\bf Z})^{\Gamma} \simeq {{\bf Z}}^m$, and 
$(\hat\rS_{F/K})^{\Gamma}\simeq {\bf Z}^m/(1,...,1)$.

\medskip

If $E_i = F_i \times F_i$, then $\pi$ sends $\rI_{E_i/K}({\bf Z})^{\Gamma}\simeq 
\rI_{F_i/K}({\bf Z})^{\Gamma} \times \rI_{F_i/K}({\bf Z})^{\Gamma}$ surjectively onto $\rI_{F_i/K}({\bf Z})^{\Gamma}\simeq {\bf Z}$. If $E_i = K$, then 
$\rI_{E_i/K}({\bf Z}) \simeq {\bf Z} \simeq \rI_{F_i/K}({\bf Z}).$ 
If $E_i$
is a quadratic field extension of $F_i$, the map $\pi$ sends
$\rI_{E_i/K}({\bf Z})^{\Gamma}\simeq {\bf Z}$ to
$\rI_{F_i/K}({\bf Z})^{\Gamma}\simeq {\bf Z}$ by multiplication by 2. Recall that  $m = m_1 + m_2$, where $m_1$ is the number of indices $i$
such that $E_i$
is a quadratic field extension of $F_i$, and $m_2$ the number of indices $i$ such that either $E_i = F_i \times F_i$ or $E_i = K$. Then
we have
$$(\hat S_{F/K})^{\Gamma} / \pi ((\hat S_{E/K})^{\Gamma})  \simeq
({\bf Z} /2 {\bf Z})^{m_1}/(1,...,1).$$

\medskip

We claim that  $\delta : (\hat S_{F/K})^{\Gamma} / \pi ((\hat S_{E/K})^{\Gamma} )  \to H^1(K,{\rm sc}(\hat T))$ sends bijectively $\sha(E',\sigma)$ to $\sha^1(K,\sico(\hat T))$.

\medskip
Let $(I_0,I_1) \in \sha (E',\sigma)$, let $a$ be the corresponding element in $$(\hat S_{F/K})^{\Gamma} / \pi ((\hat S_{E/K})^{\Gamma} )$$
 and let $x$ be the image of $a$ in $\rH^{1}(K,\sico(\hat{T}))$. We
 claim that $x$ is in $\sha^{1}(K,\sico(\hat {T}))$. It suffices to prove that for any $v\in\Omega_K$, we have
 $a^v=0$.

\medskip 

For a place $v\in\underset{i\in\rI_1}\cap\Sigma_{i}$, we have that
$E_i^v$ splits over $F_i^v$ for all $i\in\rI_1$. Hence, $\pi$
maps
$\rI_{E_i^v/K_v}({\bf Z})^{\Gamma_v}\simeq\rI_{F_i^v/K_v}({\bf Z})^{\Gamma_v}\oplus\rI_{F_{i}^v/K_v}({\bf Z})^{\Gamma_v}$
onto $\rI_{F_{i}^v/K_v}({\bf Z})^{\Gamma_v}$ for each $i\in\rI_1$, so $(\hat S_{F/K})^{\Gamma_v} / \pi ((\hat S_{E/K})^{\Gamma_v}) = 0$
for each $i\in\rI_1$ and $a^v_{i}=0$. On the other hand, for each
$i\in\rI_0$, $a_i=0$ by definition. Therefore, we have $a^v=0$.

\medskip

For a place $v\in\underset{i\in\rI_0}{\cap}\Sigma_{i}$, we replace
$(a_1,...,a_{m_1})$ by $(a_1,...,a_{m_1})+(1,...,1).$ Note that
$(a_1,...,a_{m_1})+(1,...,1)$ and $(a_1,...,a_{m_1})$  represent the same
class $a$ in $(\hat S_{F/K})^{\Gamma} / \pi ((\hat S_{E/K})^{\Gamma})$.
By the same argument  as above, we have $a_v=0$.  Since 
$(\underset{i\in\rI_0}{\cap}\Sigma_i)\cup(\underset{j\in\rI_1}{\cap}\Sigma_j)=\Omega_K$,
we have $a^v=0$ for all $v\in\Omega_K$, which proves that $x$ is in
$\sha^{1}(K,\sico(\hat{ T}))$.

\medskip
This proves that $\delta$ induces a map $\sha(E',\sigma) \to \sha^1(K,\sico(\hat T))$. We already know that this map is injective. Let
us prove that it is also surjective.

\medskip

Let $0\neq x\in\sha^1(K,\sico(\hat T))$. Let
$a\in (\hat S_{F/K})^{\Gamma} / \pi ((\hat S_{E/K})^{\Gamma})$
be the preimage of $x$, let $a^v$ be the localization of $a$ at the
place $v$, and let $(a_1,...,a_{m_1})$ be a lift of $a$ in $({\bf Z} /2 {\bf Z})^{m_1}$.
Let $(I_0,I_1)$ be the corresponding partition. 
Now we claim that 
$(\underset{i\in\rI_0}{\cap}\Sigma_i)\cup(\underset{j\in\rI_1}{\cap}\Sigma_j)=\Omega_K$.
Suppose that
$(\underset{i\in\rI_0}{\cap}\Sigma_i)\cup(\underset{j\in\rI_1}{\cap}\Sigma_j)\neq\Omega_K$, and
let
$v\in\Omega_K\backslash(\underset{i\in\rI_0}{\cap}\Sigma_i)\cup(\underset{j\in\rI_1}{\cap}\Sigma_j)$.
Therefore there exist $i_0\in\rI_0$ and $i_1\in\rI_1$ such
that $E_{i_0}^v$ is not split over $F_{i_0}^v$ and $E_{i_1}^v$ is
not split over $F_{i_1}^v$. Let
$F_{i}^v=\overset{n_{i}}{\underset{j=1}\prod} L_{i,j}$, where the
$L_{i,j}$'s are field extensions of  $K_v$. Let
$E_{i}^v =\overset{n_{i}}{\underset{j=1}\prod} M_{i,j}$, where
$M_{i,j}$ is a quadratic \'etale algebra over $L_{i,j}$. Set $\Gamma_v = \Gamma_{K_v}$. Then we
have
$$\rI_{F_i^v/K_v}({\bf Z})^{\Gamma_v}/\pi(\rI_{{E_i^v}/K_v}({\bf Z})^{\Gamma_v})=\overset{n_i}{\underset{j=1}{\oplus}}\rI_{L_{i,j}/K_v}({\bf Z})^{\Gamma_v}/\pi(\rI_{M_{i,j}/K_v}({\bf Z})^{\Gamma_v}).$$
If $M_{i,j}$ is split over $L_{i,j}$, then
$$\rI_{M_{i,j}/K_v}({\bf Z})^{\Gamma_v}=\rI_{L_{i,j}\times
L_{i,j}/K_v}({\bf Z})^{\Gamma_v}=\rI_{L_{i,j}/K_v}({\bf Z})^{\Gamma_v}\oplus\rI_{L_{i,j}/K_v}({\bf Z})^{\Gamma_v},$$
so the map $\pi$ sends $\rI_{M_{i,j}/k_v}({\bf Z})^{\Gamma_v}$ surjectively
to  $\rI_{L_{i,j}/k_v}({\bf Z})^{\Gamma_v}$. On the other hand, if
$M_{i,j}$ is a field extension over $L_{i,j}$, then
 $\pi$ maps
$\rI_{M_{i,j}/K_v}({\bf Z})^{\Gamma_v}\simeq {\bf Z}$ to
$2 {\bf Z} \subseteq {\bf Z} \simeq\rI_{L_{i,j}/K_v}({\bf Z})^{\Gamma_v}$ and we
have
$$\rI_{L_{i,j}/K_v}({\bf Z})^{\Gamma_v}/\pi(\rI_{M_{i,j}/K_v}({\bf Z})^{\Gamma_v}) \simeq {\bf Z}/ 2 {\bf Z}.$$

For
$a_i\in\rI_{F_i/K}({\bf Z})^{\Gamma}/\pi(\rI_{E_i/K}({\bf Z})^{\Gamma})\simeq {\bf Z}/2{\bf Z},$
the localization map sends $a_i$ diagonally into to
$\rI_{F_i^v/K_v}({\bf Z})^{\Gamma_v}/\pi(\rI_{E_i^v/K_v}({\bf Z})^{\Gamma_v})\simeq\underset{\overset{j,\
 where M_{i,j}}{\ is\ non-split}}{\bigoplus} {\bf Z}/2{\bf Z}$. Let
 $a_{i}^{v}$ be the image of $a_i$ in
 $\rI_{F_i^v/K_v}({\bf Z})^{\Gamma_v}/\pi(\rI_{E_i^v/K_v}({\bf Z})^{\Gamma_v})$.
 By our choice of $v$, we have
 $\rI_{F_{i_0}^v/K_v}({\bf Z})^{\Gamma_v}/\pi(\rI_{E_{i_0}^v/K_v}({\bf Z})^{\Gamma_v})$
 (resp.
 $\rI_{F_{i_1}^v/K_v}({\bf Z})^{\Gamma_v}/\rI_{E_{i_1}^v/K_v}({\bf Z})^{\Gamma_v}$)
 non-trivial. In particular, $a_{i_1}^{v}$ is non-zero as $a_{i_1}$ is non-zero.
 Note that $$\underset{i}{\oplus} (\hat\rS_{F_i^v/K_v})^{\Gamma_v}/\pi((\hat\rS_{E_i^v/K_v})^{\Gamma_v} ) =
 \frac{\underset{i}{\oplus}\rI_{F_i^v/K_v}({\bf Z})^{\Gamma_v}/\pi(\rI_{E_i^v/K_v}({\bf Z})^{\Gamma_v})}{(\ol{1},...,\ol{1})},$$
 where $\ol{1}$ denotes the image of the diagonal element
  of
 $\rI_{F_i^v/K_v}({\bf Z})^{\Gamma_v}$ in
$$\rI_{F_i^v/K_v}({\bf Z})^{\Gamma_v}/\pi(\rI_{E_i^v/K_v}({\bf Z})^{\Gamma_v}).$$
 Since
 $a^v=0$, either $a_{i}^{v}=0\in\rI_{F_i^v/K_v}({\bf Z})^{\Gamma_v}/\pi(\rI_{E_i^v/K_v}({\bf Z})^{\Gamma_v})$ for all $i$,
 or $a_{i}^{v}=\ol{1}\in\rI_{F_i^v/K_v}({\bf Z})^{\Gamma_v}/\pi(\rI_{E_i^v/K_v}({\bf Z})^{\Gamma_v})$ for all $i$. In particular, this implies that
 $a_{i_0}^{v}$ and $a_{i_1}^{v}$ are both $0$ or both $1$, which is a contradiction.
 Therefore we have 
 $(\underset{i\in\rI_0}{\cap}\Sigma_i)\cup(\underset{j\in\rI_1}{\cap}\Sigma_j)=\Omega_K$
 and $(I_0,I_1) \in \sha(E',\sigma)$. This completes the proof of the Proposition.

\bigskip
{\bf Proof of Proposition A 4.1. when $L \not = K$.}

\medskip
In
this case, the torus $\sico(T)$ fits in the following exact
sequence:
\begin{equation}\label{e5}
\xymatrix@C=0.5cm{
  1 \ar[r] & \sico( T) \ar[r] & \rR_{F/K}(\rR^{(1)}_{E/F}({\bf G}_m)) \ar[r] & \rR^{(1)}_{L/K}({\bf G}_m) \ar[r] & 1
  }
\end{equation}

\medskip
We take the dual sequence of
exact sequence (3) :
\begin{equation}\label{e6}
\xymatrix{
  0 \ar[r] & \hat{\rS}_{L/K} \ar[r]^{\iota} & \rI_{F/K}(\hat{\rS}_{E/F}) \ar[r]^{p} & \sico(\hat T) \ar[r] & 0
  },
\end{equation}

from which we derive the long exact sequence

\begin{equation}\label{e7}
\xymatrix{
  ... \ar[r] &  \rH^1(K,\hat{\rS}_{E/K}) \ar[r]^{\iota^1} & \rH^1(K,\rI_{F/K}(\hat{\rS}_{E/F}))
  \ar[r]^{p^1} & \rH^1(K,\sico(\hat T)) \ar[r] & \rH^2(K,\hat{\rS}_{E/K})
  }.
\end{equation}
By Poitou-Tate duality, we have
$\sha^2(K,\hat\rS_{E/K})\simeq\sha^1(K,\rR^{(1)}_{E/K}({\bf G}_m))^{\ast}$.
We claim that
$\sha^2(K,\hat\rS_{E/K})\simeq\sha^1(K,\rR^{(1)}_{E/K}({\bf G}_m))^{\ast}=0$.
To see this, we consider the following exact sequence:
\begin{equation*}
\xymatrix{
  1 \ar[r] & \rR^{(1)}_{L/K}({\bf G}_m) \ar[r] & \rR_{L/K}({\bf G}_m) \ar[r] & {\bf G}_m \ar[r] & 1
  },
\end{equation*}
By Hilbert Theorem 90, we have
$\rH^1(K,\rR^{(1)}_{L/K}({\bf G}_m))=K^{\times}/\rN_{L/K}(L^{\times})$,
where $\rN_{L/K}$ is the norm map from $L$ to $K$. Since the norms of
the quadratic extension $L$ over $K$ satisfy the local-global
principle, we have $\sha^1(K,\rR^{(1)}_{L/K}({\bf G}_m))=0$. Hence
$\sha^2(K,\hat\rS_{L/K})=0$. Therefore the Tate--Shafarevich group
$\sha^1(K,\sico(\hat T))$ lies in the image of
$\rH^1(K,\rI_{F/K}(\hat{\rS}_{E/F}))$.

\medskip

Let us consider the following exact sequence:
\begin{equation}\label{e9}
\xymatrix{
  1 \ar[r] & {\bf G}_m \ar[r] & \rR_{L/K}({\bf G}_m) \ar[r]^{\pi} & \rR^{(1)}_{L/K}({\bf G}_m) \ar[r] & 1
  },
\end{equation}
where $\pi(x)=x/\sigma(x)$. Considering the dual sequence, we get
\begin{equation}\label{e10}
\xymatrix{
  0 \ar[r] & \hat\rS_{L/K} \ar[r] & \rI_{L/K}({\bf Z}) \ar[r]^{d} & {\bf Z} \ar[r] & 0
  },
\end{equation}
where $d$ is the degree map which maps
$(a,b)\in {\bf Z} \oplus {\bf Z} \simeq\rI_{L/K}({\bf Z})$ to $a+b$. Taking the
long exact sequence of (\ref{e10}), we have
\begin{equation}\label{e11}
\xymatrix{
  \rI_{L/K}({\bf Z})^{\Gamma}  \ar[r]^{d} & {\bf Z} \ar[r] & \rH^1(K,\hat\rS_{L/K}) \ar[r] & \rH^1(K,\rI_{L/K}({\bf Z}))=0
  }.
\end{equation}
Since $L$ is a quadratic field extension of $K$, we have
$$\rH^1(K,\hat\rS_{L/K})\simeq {\bf Z} /d(\rI_{L/K}({\bf Z})^{\Gamma})= {\bf Z}/2{\bf Z}.$$

Similarly, we have 
\begin{center}
$\rH^1(K,\rI_{F/K}(\hat\rS_{E/F}))=\rH^1(F,\hat \rS_{E/F})=
\underset{i=1}{\overset{m}{\prod}}\rH^1(F_i,\hat\rS_{E_i/F_i})$
\end{center}

\smallskip
If $E_i = F_i \times F_i$, then $\rH^1(F_i,\hat\rS_{E_i/F_i}) = 0$ since $d$ is surjective. If $E_i$ is a quadratic extension of $F_i$, then $ \rH^1(F_i,\hat\rS_{E_i/F_i})
= {\bf Z}/2 {\bf Z}.$ Recall that  $m = m_1 + m_2$ where $m_1$ is the number of indices $i$ such that $E_i$ is a quadratic extension of $F_i$, and $m_2$ the number
of indices $i$ such that $E_i = F_i \times F_i$. Then we have $\rH^1(K,\rI_{F/K}(\hat\rS_{E/F})) \simeq ({\bf Z}/2 {\bf Z})^{m_1}$. 

\smallskip
The map
$\iota^1:\rH^1(K,\hat\rS_{L/K})\ra\rH^1(K,\rI_{F/K}(\hat\rS_{E/F}))$
maps ${\bf Z} / 2 {\bf Z}$ diagonally into $({\bf Z} / 2 {\bf Z})^{m_1}$.
Therefore, we have
$$\sha^1(k,\sico(\hat T))\subseteq\img(p^1)
\simeq({\bf Z} / 2 {\bf Z})^{m_1}/(1,...,1).$$

\medskip
Let us show that $p^1$ maps $\sha (E',\sigma)$ bijectively to $\sha (K, {\rm sc}(\hat T))$.

\medskip
Let $(I_0,I_1) \in \sha (E',\sigma)$, and let $a$ in
$\rH^1(K,\rI_{F/K}(\hat{\rS}_{E/F}))$ be the corresponding element. We
want to show that $p^1(a)\in\sha^1(K,\sico(\hat T)).$ Let
$v\in\Omega_K$. If $v\in\Sigma(L/K)$ or
$v\in\underset{j\in\rI_1}{\cap}\Sigma_j$, then $a^v=0$. Hence, it
suffices to prove that for
$v\in\Omega_K \setminus(\Sigma(L/K) \cup(\underset{i\in\rI_1}{\cap}{\Sigma}_i))$,
we have $a^v =\iota_v^1(1)=\iota^1(1)_v$. Since
$\Sigma (L/K) \cup(\underset{i\in\rI_0}{\cap}\Sigma_i)\cup(\underset{j\in\rI_1}{\cap}\Sigma_j)=\Omega_K,$
we have $v\in(\underset{i\in\rI_0}{\cap}{\Sigma}_i)$. Consequently,
for all $i\in\rI_0$,  we have $\rH^1(F_{i}^{},\hat{\rS}_{E_i^v/F_i^v})=0$
and the projection of $\iota_v^1(1)$ to these components are
trivial. For $i\in\rI_1$, we have that $a_i$ and the $i-th$ coordinate of
$\iota^{1}(1)$ are both 1, so their images in
$\rH^1(F_{i}^{v},\hat{\rS}_{E_i^v/F_i^v})$ are equal. This proves that
$a^v=\iota_v^1(1)$, hence $p^1(a^v) = 0$. 

\medskip
We next show that the restriction of $p^1$ to $\sha (E',\sigma)$ is surjective onto $\sha^1(K,\sico(\hat T)).$

\medskip

Let
$a=(a_1,...,a_{m_1})\in({\bf Z} / 2 {\bf Z})^{m_1} \simeq\rH^1(K,\rI_{E^{\sigma}/k}(\hat\rS_{E/F}))$
and let  $(I_0,I_1)$ be the associated partition. If $a=0$ or
$a=(1,...,1)$, then $a$ is in the image of $\iota^1$ and
$p^1(a)=0\in\sha^1(K,\sico(\hat T))$. Hence, we may assume that $I_0$
and $I_1$ are non-empty. We claim that
\begin{center}
$0\neq p^1(a)\in\sha^1(K,\sico(\hat T))$ if and only if $\rI_0$ and
$\rI_1$ are non-empty and

\medskip

$\Sigma (L/K) \cup(\underset{i\in\rI_0}{\cap}\Sigma_i)\cup(\underset{j\in\rI_1}{\cap}\Sigma_j)=\Omega_K$.
\end{center}

Suppose that $0\neq p^1(a)\in\sha^1(K,\sico(\hat T)).$  Let
$v\in\Omega_K\setminus(\Sigma (L/K) \cup(\underset{i\in\rI_0}{\cap}{\Sigma}_i))$.
Since $v\notin\Sigma (L/K) $, we have
$\rH^1(L^{v},\hat{\rS}_{L^{v}/K_{v}})= {\bf Z} / 2 {\bf Z}$. Let $a^v$ denote
the localization of $a$ at $v$. Since
$p^1(a)\in\sha^1(K,\sico(\hat T))$, we have $a^v$ in the image of
$\iota_v^1$, so either $a^v=0$ or $a^v=\iota_v^1(1)$. It suffices to
show that  $v\in\underset{i\in\rI_1}{\cap}\Sigma_i$. Consider the $i$-th
component of $({\bf Z} / 2 {\bf Z})^{m_1}$, which corresponds to
$\rH^1(K,\rI_{F_i/K}(\hat{\rS}_{E_i/F_i}))=\rH^1(F_i,\hat{\rS}_{E_i/F_i})$.
If $E_i$ splits over $F_i$ at  a place $v\in\Omega_K$, then by the exact
sequence (\ref{e11}), the map $d$ is surjective and
$\rH^1(F_{i}^{v},\hat{\rS}_{E_i^v/F_i^v})=0$, which means that the
$i$-th component vanishes at place $v$. Since
$v\notin\underset{i\in\rI_0}{\cap}{\Sigma}_i$, there exists $i\in\rI_0$
such that $E_{i}^{v}$ is not split over $F_{i}^{v}$. Let
$F_{i}^{v}=\overset{n_{i}}{\underset{j=1}\prod} L_{i,j}$, where
$L_{i,j}$'s are field extensions of $K_v$. Let
$E_{i}^{v}=\overset{n_{i}}{\underset{j=1}\prod} M_{i,j}$, where
$M_{i,j}$ is a quadratic \'etale algebra over $L_{i,j}$. Then
$$\rH^1(F_{i}^{v},\hat{\rS}_{E_i^v/F_i^v})=\underset{j}{\prod}\rH^1(L_{i,j},\hat{\rS}_{M_{i,j}/L_{i,j}}).$$
By the choice of $i$, there is $j$ such that $M_{i,j}$ is not split
over $L_{i,j}$, and hence
$\rH^1(L_{i,j},\hat{\rS}_{M_{i,j}/L_{i,j}})\neq 0$. Then the
projection of $\iota_v^1(1)$  to
$\rH^1(L_{i,j},\hat{\rS}_{M_{i,j}/L_{i,j}})$ is 1. On the other
hand, the projection of $a^v$ to the same component is 0 since
$i\in\rI_0$. Therefore, $a^v=0$ which means that 
$\rH^1(F_{i}^{v},\hat{\rS}_{E_i^v/F_i^v})=0$ for all $i\in\rI_1$,
hence  $v\in\underset{i\in\rI_1}{\cap}\Sigma_i$. This proves that $a \in \sha (E',\sigma)$. 

\bigskip

\centerline {\bf Appendix B}

\medskip

\centerline {\bf Orthogonal involutions and maximal subfields}

\medskip

Let $K$ be a global field of characteristic $\not = 2$, and let $\Omega_K$ be the set of places of $K$. 
The purpose of this appendix is to give a new proof of Theorem B of Prasad and Rapinchuk (see [PR 10], Introduction, page 586)  using some of the results of the present paper, in particular the local embedding conditions of \S 3. Theorem B has two parts. The proof of the first part will be presented in the first section of this appendix,
and the proof of the second one in the second section. The application to maximal subfields determining orthogonal involutions of degree divisible by 4 is
given in \S B3.

\medskip
We thank Gopal Prasad for encouraging us to include our proofs of these results in our paper.

\medskip

{\bf \S B1. First part of Theorem B}

\medskip

We start by recalling the dimension condition :

\medskip
\noindent
{\bf Definition B.1.1.} Let $E$ be an $n$--dimensional  commutative \'etale $K$--algebra, and let $\sigma : E \to E$ be a $K$--linear involution. Let $F = E^{\sigma}$. 
We say that $(E,\sigma)$ satisfies the {\it dimension condition} if  ${\rm dim}_K(F) = [{{n+1} \over 2}]$.

\medskip
The following result was proved by Prasad and Rapinchuk (cf. [PR 10], Theorem B {\rm (i)}, Introduction, page 586; see also [PR 10], Theorem 9.4.) :

\medskip
\noindent
{\bf Theorem B.1.} {\it Let $(A_1,\tau_1)$ and $(A_2,\tau_2)$ be two central simple algebras with orthogonal involutions with ${\rm deg}(A_1) = {\rm deg}(A_2) = n$,
with $n \ge 3$.
Suppose that the $A_i$'s have the same isomorphism classes of $n$--dimensional commutative \'etale algebras
invariant under the involutions satisfying the dimension condition. Then for all $v \in \Omega_K$, we have $(A^v,\tau_1) \simeq (A^v,\tau_2).$}

\medskip
\noindent
{\bf Definition B.1.2.} Let $(A,\tau)$ be an orthogonal involution, with ${\rm deg}(A) = n$. 
 Let $v \in \Omega_K$ and let $E$ be a $\tau$--invariant
rank $n$ \'etale subalgebra of $A^v$ satisfying the dimension condition. An $n$--dimensional  $\tau$--invariant
subalgebra $\tilde E$  of $A$ is called a 
{\it lifting} of $E$ if $(\tilde E,\tau) \simeq (E,\tau)$ over $K_v$. 

\medskip
\noindent
{\bf Lemma B.1.3.} {\it  Let  $(A,\tau)$ be an orthogonal involution, with ${\rm deg}(A) = n$.
Let $v \in \Omega_K$ and let $E$ be an $n$--dimensional  $\tau$--invariant \'etale
subalgebra of $A^v$ satisfying the dimension condition. 
Then $E$ has a lifting in $A$.}

\medskip
\noindent
{\bf Proof.} See for instance [PR 10] , Proposition 2.4. 

\medskip
\noindent
{\bf Lemma B.1.4.} {\it Let $(A_1,\tau_1)$ and $(A_2,\tau_2)$ be two central simple algebras with orthogonal involutions with ${\rm deg}(A_1) = {\rm deg}(A_2) = n$.
Suppose that the $A_i$'s have a common  $n$--dimensional commutative \'etale subalgebra
invariant under the involutions satisfying the dimension condition. Then we have
${\rm disc}(A_1,\tau_1) = {\rm disc}(A_2,\tau_2)$.}

\medskip
\noindent
{\bf Proof.} Let $E$ be an $n$--dimensional  $\tau_1$--invariant subalgebra of $A_1$ satisfying the dimension condition such that $(E,\tau_1)$ embeds
into $(A_2,\tau_2)$. By Proposition 1.6.1. we have ${\rm disc}(A_1,\tau_1) =
{\rm disc}(E) =  {\rm disc}(A_2,\tau_2)$.

\medskip
\noindent
{\bf Proposition B.1.5.} {\it Let $(A_1,\tau_1)$ and $(A_2,\tau_2)$ be two central simple algebras with orthogonal involutions with ${\rm deg}(A_1) = {\rm deg}(A_2) = n$,
with $n \ge 3$. 
Suppose that the $A_i$'s have the same isomorphism classes of $n$--dimensional commutative \'etale algebras
invariant under the involutions satisfying the dimension condition. Then we have the following

\medskip
{\rm (i)}  ${\rm disc}(A_1,\tau_1) = {\rm disc}(A_2,\tau_2)$.

\medskip
{\rm (ii)} The algebras $A_1$ and $A_2$ are isomorphic.

\medskip
{\rm (iii)} Let $v \in \Omega_K$ be such that $A_1^v$ and $A_2^v$ are split, and let $q_1$ and $q_2$ be quadratic forms inducing the
involutions $\tau_1$ and $\tau_2$. Then the Witt indices of $q_1$ and $q_2$ are equal.}

\medskip
\noindent
{\bf Proof.} {\rm (i)} This follows from Lemma B.1.4.

\medskip
{\rm (ii)} If $n$ is odd, then $A_1$ and $A_2$ are both split, hence they are isomorphic. Let us assume that $n$ is even, and set $n = 2r$. Let $v \in \Omega_K$. Then
for $i = 1,2$, either $A^v_i$ is split or it is isomorphic to $M_r(D)$, where $D$ is the unique quaternion division algebra over $K_v$. Hence it suffices to
prove that for all $v \in \Omega$, the algebra $A^v_1$ splits if and only if $A_2^v$ splits.

\medskip
Let $v \in \Omega_K$ be a place such that $A_1^v \simeq M_r(D)$. Let $\tau_1$ be induced by a hermitian form $h_1$ over $D$. Suppose that
$A_2^v$ splits, and that $\tau_2$ is induced by a quadratic form $q_2$ over $K_v$. We claim that $q_2$ is isotropic.

\medskip
Assume first that $v$ is a finite place. For $n \ge 5$ all quadratic forms are isotropic (cf. [Sch 85], 6.4.2.). For $n = 4$, if $q_2$ is
anisotropic, then the determinant of $q_2$ is trivial. Since ${\rm disc}(A_1,\tau_1) = {\rm disc}(A_2,\tau_2) = 1$, the $2$--dimensional hermitian form $h_1$
over $D$ is hyperbolic (cf. [T 61], Theorem 3.).  Then there exists an $n$--dimensional  $\tau_1$--invariant  \'etale subalgebra $E_1$ of $A_1^v$ satisfying
the dimension condition such that $(E_1,\tau_1)$ is split. Let $E$ be a lifting of $E_1$ in $A_1$. Then $(E,\tau_1)$ can be embedded into $(A_2,\tau_2)$. 
Since $(E^v,\tau_1) \simeq (E_1^v,\tau_1)$ is split, by
Proposition 3.1.1. this implies that $(A^v_2,\tau_2)$ is hyperbolic, which contradicts the assumption that $q_2$ is anisotropic. Therefore $q_2$ is
isotropic.

\medskip
Suppose now that $v$ is a real place. Let  $k = [{r  \over 2}]$, and let $E_1 = {\bf C}^k \times {\bf C}^k \times {\bf C}^{r - 2k}$, and let $\sigma : E_1 \to E_1$ be
the involution which exchanges the two copies of ${\bf C}^k$, and acts on ${\bf C}^{r - 2k}$ as the complex conjugation. Then $(E_1,\sigma)$
can be embedded into $(A_1^v,\tau_1)$ (cf. Proposition 3.1.2.), hence we can assume that the restriction of $\tau_1$ to $E_1$ is $\sigma$. Let $E$ be a
lifting of $E_1$ in $A_1$. Since $(E,\tau_1)$ can be embedded into $(A_2,\tau_2)$, by Proposition 3.1.2. the signature of $q_2$ is of the form
$(2k + 2s, 2k + 2s')$ for some non--negative integers $s, s'$. Note that $k \ge 1$, since $n \ge 3$. Hence $q_2$ is isotropic. 

\medskip

Let $v \in \Omega_K$, and let us write $q_2 \simeq q_0 \oplus q$, with $q_0$ hyperbolic and $q$ anisotropic. Let ${\rm dim}(q_0) = 2m$; since $q_2$ is isotropic, we have $m \ge 1$. 
Then there exists a $\tau_2$--invariant commutative \'etale subalgebra $E_2 = K_v^m \times K_v^m \times E'$ of $A^v_2$ satisfying the dimension condition such
that $\tau_2$ permutes the two copies of $K_v^m$. Let $E$ be a lift of $E_2$. Since $(E,\tau_2)$ can be embedded into $(A_1,\tau_1)$, we see that
$K_v$ splits $A_1^v$, which is a contradiction.

\medskip Therefore for all $v \in \Omega_K$, the algebra $A_1^v$ splits if and only if $A^v_2$ splits. This implies that $A_1$ and $A_2$ are isomorphic.
If $v$ is a place such that $A_1^v$ and $A_2^v$ are split,  the above argument also shows that the Witt indices of $q_1$ and $q_2$ are equal, and this
proves {\rm (iii)}. 

\medskip
\noindent
{\bf Proposition B.1.6.} {\it Let $A$ be a central simple $K$--algebra of degree $n$, and let $\tau_1 : A \to A$ and $\tau_2 : A \to A$ be two orthogonal 
involutions. Assume that we have

\medskip
{\rm (i)}  ${\rm disc}(A,\tau_1) = {\rm disc}(A,\tau_2)$.

\medskip
{\rm (ii)} Let $v \in \Omega_K$ be such that $A^v$ is split, and let $q_1$ and $q_2$ be quadratic forms inducing the
involutions $\tau_1$ and $\tau_2$. Then the Witt indices of $q_1$ and $q_2$ are equal.

\medskip
Then for all $v \in \Omega_K$, the algebras with involution $(A^v,\tau_1)$ and $(A^v,\tau_2)$ are isomorphic.}

\medskip
\noindent
{\bf Proof.}  If $v \in \Omega_K$ is a real place such that $A^v$ splits, 
 then having the same Witt index implies that $q_1$ and $q_2$ have the same signature, hence they are isomorphic.
Therefore we have $(A^v,\tau_1) \simeq (A^v,\tau_2)$. For a real place $v$ such that $A^v$ is not split, we have $(A^v,\tau_1) \simeq (A^v,\tau_2)$ (cf.
[Sch 85], 10.3.7). Therefore for all real places $v \in \Omega_K$, we have $(A^v,\tau_1) \simeq (A^v,\tau_2)$.

\medskip
Let $v \in \Omega_K$ be a finite place. Assume first that $A^v$ is non-split. By {\rm (i)} we have  ${\rm disc}(A^v,\tau_1)  = {\rm disc}(A,\tau_2)$,
and this implies that $(A^v,\tau_1) \simeq (A^v,\tau_2)$ (cf. [T 61], Theorem 3.).

\medskip Assume now that $v$ is a finite place such that $A^v$ is split. Since $q_1$ and $q_2$ have the same Witt index, it remains to prove that their
anisotropic parts are similar. Let $n_0$ be the dimension of the anisotropic parts of $q_1$ and $q_2$. We are reduced to the case where $n_0 \le 4$,
and the quadratic forms $q_1$ and $q_2$ are anisotropic of dimension $n_0$. If $n_0 = 4$, then there is only one isomorphism class of anisotropic forms,
hence $q_1 \simeq q_2$. Recall that we have ${\rm det}(q_1) = {\rm det}(q_2)$. Two anisotropic forms of dimension $\le 3$ having the same determinant are
similar. Therefore $q_1$ and $q_2$ are similar, hence $(A^v,\tau_1) \simeq (A^v,\tau_2)$. This completes the proof of the Proposition.

\medskip
\noindent
{\bf Proof of Theorem B.1.} By Proposition B.1.5. we can assume that $A_1 = A_2 = A$, and that conditions {\rm (i)}  and {\rm (ii)} of Proposition B.1.6. are
fulfilled. Hence Proposition B.1.6. implies the Theorem. 

\medskip

If $n$ is even,  having the same invariant {\it subfields} is enough. In order to prove this, we need a few lemmas :

\medskip
\noindent
{\bf Lemma B.1.7.} {\it Let $k$ be a local field, let $r \ge 1$ be an integer, and let $\delta \in k^{\times}/k^{\times 2}$. Assume that one of the following holds

\medskip
{\rm (i)} $\delta \not = 1$ in $k^{\times}/k^{\times 2}$.

\medskip
{\rm (ii)} $r$ is even, and $\delta  = 1$ in $k^{\times}/k^{\times 2}$.

\medskip
Then we have the following :

\medskip

{\rm (iii)} There exists a degree $r$ field extension $M$ of $k$ and $x \in M^{\times}$ such that $x \not \in M^{\times 2}$ and that ${\rm N}_{M/k}(x) = \delta \in k^{\times}/k^{\times 2}$.

\medskip

{\rm (iv)} There exists a tower of field extensions $N / M / k$ with $[M:k] = r$ and $[N:M] = 2$ such that the discriminant of $N$ is ${\delta}$.}

\medskip
\noindent
{\bf Proof.} Let us prove {\rm (iii)}. Assume first that ${\rm (i)}$ holds. Suppose that $\delta$ is a unit, and let $M$ be the unique unramified extension of $k$ of degree $r$.
Let $x$ be a unit of $M$ such that ${\rm N}_{M/k}(x) = \delta$. Then $x \not \in M^{\times 2}$.

\medskip
Suppose that $\delta = \pi$, where $\pi$ is a uniformizer of $K$. Let $f(X) = X^r  +(-1)^{r} \pi$, and let $M = k[X]/(f)$. Let $x$ be the image of $X$ in $M$. Then we have ${\rm N}_{M/k}(x) = \pi$, and $x \not \in M^{\times 2}$.

\medskip
Assume that {\rm (ii)} holds, and let $r = 2m$. Let $M/k$ be unramified of degree $2m$. Let $\pi$ be a uniformizer of $k$. Then $\pi$ is also a uniformizer of $M$,
hence $\pi \not \in M^{\times 2}$. We have ${\rm N}_{M/K}(\pi) = \pi^{2m}$, hence ${\rm N}_{M/K}(\pi)  \in k^{\times 2}$. Set $x = \pi$. 

\medskip
Therefore {\rm (iii)} is proved. Let us prove {\rm (iv)}. With $M$ and $x$ as in {\rm (iii)}, let us set $N = M(\sqrt x)$. Then $N$ and $M$ have the required properties,
hence {\rm (iv)} is proved. 

\medskip
\noindent
{\bf Lemma B.1.8.} {\it  Let  $(A,\tau)$ be an orthogonal involution, with ${\rm deg}(A) = 2r$ with $r \ge 2$.
Let $S$ be a finite subset of $\Omega_K$ and for all $v \in S$, let $E(v)$ be an $n$--dimensional  $\tau$--invariant \'etale
subalgebra of $A^v$ satisfying the dimension condition. 
Then the algebras $E(v)$ for $v \in S$ have a common lifting in $A$ which is a
degree $2r$ field extension of $K$.}

\medskip
\noindent
{\bf Proof.} Assume first that $r$ is even, or ${\rm disc}(A,\tau) \not = 1  \in K^{\times}/K^{\times 2}$. Let $\delta = {\rm disc}(A,\tau) \in K^{\times}/K^{\times 2}$. Let $w \in \Omega_K$ be a finite place such that  $w \not \in S$, and that if $\delta \not = 1 \in K^{\times}/K^{\times 2}$, then
$\delta \not \in K_w^{\times 2}$. Let $N$ and $M$ be as in Lemma B.1.7. {\rm (iv)}
and let $\sigma : N \to N$ be the $K_w$--linear involution with fixed field $M$. Note that since ${\rm disc}(N) = \delta$ and $(N,\sigma)$ is not split, Proposition 
3.1.1. {\rm (ii)} implies that $(N,\sigma)$ can be embedded into $(A^w,\tau)$. By [PR 10], Proposition 2.4. there exists
an \'etale subalgebra  $E$ of $A$ which is a common lifting of $N$ and $E(v)$, for all $v \in S$.  Since $N$ is a field, $E$ is a field as well. 

\medskip
Suppose now that $r$ is odd and that $ {\rm disc}(A,\tau) = 1  \in K^{\times}/K^{\times 2}$. Then by hypothesis we have $r \ge 3$; let us write $r = m + 3$.
Let $v_1$ and $v_2$ be two distinct finite places of $K$ such that $v_1, v_2 \not \in S$, that $A^{v_1}$ and $A^{v_2}$ are split, and that
$(A^{v_2},\tau)$ is hyperbolic. 

\medskip
Let $\pi$ be a uniformizer at $v_1$, and let $\mu$ be an unit such that $\mu \not \in K_{v_1}^{\times 2}$. Then $K_{v_1}$ has exactly four square classes, and
they are represented by $1$, $\mu$, $\pi$ and $\mu \pi$. Set $E_1 = K_{v_1}(\sqrt \mu)$, $E_2 = K_{v_1}(\sqrt \pi)$, and $E_3 = K_{v_1}(\sqrt \mu \pi)$. Let
$\sigma_i : E_i \to E_i$ be the non-trivial automorphism of $E_i/K_{v_1}$ for all $i = 1,2,3$. If $m > 0$,
let $E_4 = K_{v_1}^m \times K_{v_1}^m$, and let $\sigma_4 : E_4 \to E_4$ be the involution which exchanges the two copies of $K_{v_1}$. Set $E (v_1)  = E_1 \times E_2 \times E_3 \times E_4$, and let $\sigma(v_1) : E(v_1)  \to E(v_1)$ be the involution which is equal to
$\sigma_i$ on $E_i$. Then $(E(v_1),\sigma(v_1))$ is a non--split rank $2r$ \'etale algebra of discriminant 1 satisfying the dimension condition. Hence 
$(E(v_1),\sigma(v_1))$ can be embedded into $(A,\tau)$ by Proposition 3.1.1. 
{\rm (ii)}. Let us denote by $(E(v_1),\tau)$ the image of $(E(v_1),\sigma(v_1))$ in  $(A^{v_1},\tau)$.

\medskip
Let $L$ be the unramified extension of degree $r$ of $K_{v_2}$.  Since $r$ is odd, we have
${\rm disc}(L) = 1 \in K^{\times}_{v_2} / K^{\times 2}_{v_2}$. Let $E(v_2) = L \times L$, and let $\sigma (v_2) : E(v_2) \to E(v_2)$ be the
involution exchanging the two copies of $L$. Then $(E(v_2),\sigma (v_2))$ is a split rank $2r$ \'etale algebra with involution satisfying the dimension condition.
Since $A^{v_2}$ is split and $(A^{v_2},\tau)$ is hyperbolic, $(E(v_2),\sigma (v_2))$ can be embedded into $(A^{v_2},\tau)$ by Proposition 3.1.1. {\rm (i)}. Let
us denote by $(E(v_2),\tau)$ the image of $(E(v_2),\sigma(v_2))$ in  $(A^{v_2},\tau)$.

\medskip
Let $(E,\tau)$ be a common lifting of $(E(v_1),\tau)$, $(E(v_2),\tau)$ and of $(E(v),\tau)$ for $v \in S$; such a lifting exists by [PR 10], Proposition 2.4. 
Let $F$ be the subalgebra of $\tau$--symmetric elements of $E$. Then $F$ is a field, since $F^{v_2}$ is a field. Moreover, $(E,\tau)$ is not split, since
$(E(v_1),\tau)$ is not split. Therefore $E$ is a degree $2r$ field extension of $K$. This completes the proof of the proposition.

\medskip
We have the following strengthening of Proposition B.1.5. :

\medskip
\noindent
{\bf Proposition B.1.9.} {\it Let $(A_1,\tau_1)$ and $(A_2,\tau_2)$ be two central simple algebras with orthogonal involutions with ${\rm deg}(A_1) = {\rm deg}(A_2) = 2r$
with $r \ge 2$.
Suppose that the $A_i$'s have the same isomorphism classes of $2r$--dimensional  subfields
invariant under the involutions satisfying the dimension condition. Then

\medskip
{\rm (i)}  ${\rm disc}(A_1,\tau_1) = {\rm disc}(A_2,\tau_2)$.

\medskip
{\rm (ii)} The algebras $A_1$ and $A_2$ are isomorphic.

\medskip
{\rm (iii)} Let $v \in \Omega_K$ be such that $A_1^v$ and $A_2^v$ are split, and let $q_1$ and $q_2$ be quadratic forms inducing the
involutions $\tau_1$ and $\tau_2$. Then the Witt indices of $q_1$ and $q_2$ are equal.}

\medskip
\noindent
{\bf Proof.} {\rm (i)} follows from Lemma B.1.4. The proofs of {\rm (ii)} and {\rm (iii)} follow the pattern of the proof of Proposition B.1.5. with the following
modifications. By Lemma B.1.8. the algebras $E_1$ and $E_2$ appearing in the proof of Proposition B.1.5. have liftings that are subfields of $A_1$ respectively $A_2$. Since we
are assuming that $A_i$'s have the same isomorphism classes of $2r$--dimensional  subfields
invariant under the involutions satisfying the dimension condition, the arguments of the proof of Proposition B.1.5. apply and we get the
desired conclusion. 

\medskip
Therefore we obtain the following :

\medskip
\noindent
{\bf Theorem B.1'.} {\it Let $(A_1,\tau_1)$ and $(A_2,\tau_2)$ be two central simple algebras with orthogonal involutions with ${\rm deg}(A_1) = {\rm deg}(A_2) = 2r$,
with $r \ge 2$. 
Suppose that the $A_i$'s have the same isomorphism classes of maximal subfields 
invariant under the involutions satisfying the dimension condition.  Then for all $v \in \Omega_K$, we have $(A^v,\tau_1) \simeq (A^v,\tau_2).$}

\medskip
\noindent
{\bf Proof.} By Proposition B.1.9.  we may assume that $A_1 = A_2 = A$, and that conditions {\rm (i)} and {\rm (ii)} of Proposition B.1.6. are fulfilled.
Hence Proposition B.1.6. implies the Theorem.

\bigskip
{\bf \S B2. Second part of Theorem B}

\medskip
We now prove the second part of Theorem B of Prasad and Rapinchuk (cf. [PR 10], Introduction, page 586, and Theorem 8.1.).
Let $A$ be a central simple $K$--algebra,
let $\tau : A \to A$ be an orthogonal involution of degree $4m$.

\medskip
Let ${\cal J} = {\cal J}(A,\tau)$ be the set of orthogonal involutions $\eta : A \to A$ such that 
$(A^v,\eta) \simeq (A^v,\tau)$ for all $v \in \Omega_K$.  Let $\Omega'$ be the set of  places $v \in \Omega_K$ such that $A^v$ is not split and that
$Z(A^v,\tau) = K_v \times K_v$. 

\medskip
\noindent
{\bf Theorem B.2.} {\it We have the following :

\medskip

{\rm (i)} Let $\eta \in {\cal J}$. Then there exists an $\eta$--invariant  maximal subfield $E_{\eta}$  of $A$  satisfying the
dimension condition such that $(E^v_{\eta},\eta)$ is split  for all $v \in \Omega'$. 

\medskip

{\rm (ii)} Let $\eta \in {\cal J}$, and let $E_{\eta}$ be an $\eta$--invariant  subalgebra of $A$ of rank $4m$ such that $(E^v_{\eta},\eta)$ is split for all $v \in \Omega'$.  Let $\gamma \in {\cal J}$. If $(E_{\eta},\eta)$ embeds into $(A,\gamma)$, then the algebras with involution $(A,\eta)$ and $(A,\gamma)$ are isomorphic.}

\medskip
\noindent
{\bf Proof.} {\rm (i)}  Let $v \in \Omega'$. Then we have $A \simeq M_{2m}(D(v))$, where $D(v)$ is the unique quaternion division algebra over $K_v$. By hypothesis,
we have ${\rm disc}(A,\eta) = 1 \in K^{\times}/K^{\times 2}$. By  [T 61], Theorem 3. this implies that $(A,\eta)$ is hyperbolic. Let $L(v)$ be a quadratic extension
of $K_v$ splitting
$D(v)$ and set $E(v) = L(v)^m \times L(v)^m$. Let us endow $E(v)$ with the involution $\sigma(v)  : E(v) \to E(v)$ which exchanges the two copies
of $L(v)$. By Propositions 3.1.1. and 3.1.2. there exists an embedding of algebras with involution $(E(v),\sigma (v) ) \to (A^v,\eta)$. 
By Lemma B.1.8.  there exists a $\eta$--invariant subfield $E_{\eta}$ of $A$ such that $(E_{\eta}^v,\eta) \simeq (E(v),\sigma(v))$ for all $v \in \Omega'$. 

\medskip

{\rm (ii)} Let $E_{\eta}$ be an $\eta$--invariant  \'etale subalgebra of $A$ satisfying the dimension condition such that $(E_{\eta},\eta)$ is split over $K_v$ for all $v \in \Omega'$.
Let $F_{\eta}$ be the subalgebra of $\eta$--symmetric elements of $E_{\eta}$, and let $d \in F_{\eta}^{\times}$
such that $E_{\eta} = F_{\eta}(\sqrt d)$. Let $\gamma \in {\cal J}$, and assume that there exists an embedding of algebras with involution $(E_{\eta},\eta) \to (A,\gamma)$.
By Proposition 1.1.3. there exists $a \in F_{\eta}^{\times}$ such that $(A,\eta_a) \simeq (A,\gamma)$. We claim that $(A,\eta_a) \simeq (A,\eta)$.

\medskip

Let $\Omega_1$ be the set of places of $K$ such that $A^v$ is non-split and that $Z(A^v,\tau)$ is a field. Let $\Omega_2$ be the set of places of $K$ such
that $A^v$ is split. Note that we have $\Omega_K = \Omega' \cup \Omega_1 \cup \Omega_2$. We have compatible orientations $u : \Delta (E_{\eta}) \to Z(A,\eta)$
and $u_a : \Delta (E_{\eta}) \to Z(A,\eta_a)$. We regard $C(A,\eta_a)$ and $ C(A,\eta)$ as $\Delta (E_{\eta})$--algebras via $u_a$ and $u$.
By Proposition 2.5.4. we have $C(A,\eta_a)  = C(A,\eta) + {\rm res}_{\Delta (E_{\eta})/K}{\rm cor}_{F_{\eta}/K}(a,d)$ in ${\rm Br}(\Delta (E_{\eta}))$.

\medskip
Let $v \in \Omega'$. Then $(E^v_{\eta},\eta)$ is split  by hypothesis, hence we have $d \in (F_{\eta}^v)^{\times 2}$. Therefore ${\rm cor}_{F^v_{\eta}/K}(a,d) = 0$,
and hence we have $C(A^v,\eta_a)  = C(A^v,\eta)$ in ${\rm Br}(\Delta (E^v_{\eta}))$.

\medskip

Let $v \in \Omega_1$. Then $Z(A^v,\eta)$ is a field, hence $\Delta(E^v_{\eta})$ is also a field. Thus  ${\rm res}_{\Delta (E^v_{\eta})/K}{\rm cor}_{F^v_{\eta}/K}(a,d) = 0$,
hence we have $C(A^v,\eta_a)  = C(A^v,\eta)$ in ${\rm Br}(\Delta (E^v_{\eta}))$.

\medskip

Let $v \in \Omega_2$. Then $A^v$ is split, and hence $(A^v,\eta)$ admits improper similitudes. Hence there exists an isomorphism of algebras with
involution ${\rm Int}(\alpha) : (A^v,\eta_a) \to (A^v,\eta)$ such that $u = c(\alpha) \circ u_a$. Therefore $C(A^v,\eta_a)$ is isomorphic to $C(A^v,\eta)$ as 
$\Delta(E_{\eta}^v)$--algebras. 

\medskip
Hence we have $C(A^v,\eta_a)  = C(A^v,\eta)$ in ${\rm Br}(\Delta (E^v_{\eta}))$ for all $v \in \Omega_K$.  By the Hasse--Brauer--Noether theorem, this
implies that $C(A,\eta_a)  = C(A,\eta)$ in ${\rm Br}(\Delta (E_{\eta}))$. Note that ${\rm disc}(A,\eta_a) = {\rm disc}(A,\eta)$, and that $(A^v,\eta_a)$ and
$(A^v,\eta)$ are isomorphic for all real places $v$ of $K$.  By [LT 99], Theorems A and B this implies that $(A,\eta_a) \simeq (A,\eta)$. Since
$(A,\eta_a) \simeq (A,\gamma)$, we obtain $(A,\eta) \simeq (A,\gamma)$.

\bigskip
{\bf \S B3. Application}

\medskip
The results of \S B1 and \S B2 have the following application (see [PR 10], Introduction, page 586, last line before the statement of Theorem B) :

\medskip
{\bf Theorem B.3.} {\it Let $(A_1,\tau_1)$ and $(A_2,\tau_2)$ be two central simple algebras with orthogonal involutions with ${\rm deg}(A_1) = {\rm deg}(A_2) = 4m$.
Assume that the $A_i$'s have the same isomorphism classes of
maximal subfields invariant under the involutions satisfying the dimension condition. Then we have $(A,\tau_1) \simeq (A,\tau_2).$}

\medskip
\noindent
{\bf Proof.} By Theorem B.1.'  we have $(A^v_1,\tau_1) \simeq (A^v_2,\tau_2)$ for all $v \in \Omega_K$. 
Let ${\cal J}  = {\cal J}(A,\tau_1)$ be the set of orthogonal involutions $\eta : A \to A$ such that 
$(A^v,\eta) \simeq (A^v_2,\tau_1)$ for all $v \in \Omega_K$. Then we have $\tau_2 \in {\cal J}$. By Theorem B.2. {\rm (i)} there exists 
a $\tau_1$--invariant  maximal subfield $E_{\tau_1}$ of $A$ satisfying the dimension condition and such that $(E_{\tau_1}^v,\tau_1)$ is split over $K_v$ for all $v \in \Omega'$
(where $\Omega'$ is the set of places $v \in \Omega_K$ such that $A^v$ is not split and $Z(A^v,\tau_1) = K_v \times K_v$.
Since $(A,\tau_1)$ and $(A,\tau_2)$ have the
same isomorphism classes of maximal subfields invariant by the involutions and satisfying the dimension condition, $(E_{\tau_1},\tau_1)$ embeds
into $(A,\tau_2)$. Therefore by Theorem B.2. {\rm (ii)} we have $(A,\tau_1) \simeq (A,\tau_2).$

\bigskip

{\bf Bibliography}

\bigskip \noindent
[B 12] E. Bayer--Fluckiger, Embeddings of maximal tori in orthogonal groups, {\it Ann. Inst. Fourier} (to appear).

\medskip  \noindent
[B 13] E. Bayer--Fluckiger, Isometries of quadratic spaces, {\it J. Eur. Math. Soc.} (to appear).

\medskip
\noindent
[Bo 91] A. Borel, {\it Linear algebraic groups}, Second enlarged edition, Graduate texts in mathematics  {\bf 126}, Springer-Verlag, New York, 1991.

\medskip 
\noindent
[Bo 99]  M. Borovoi, A cohomological obstruction to the Hasse principle for homogeneous spaces, {\it Math. Ann.}
\textbf{314} (1999), 491-504.

\medskip
\noindent
[BCM 03] R. Brusamarello, P. Chuard--Koulmann and J. Morales, Orthogonal groups containing
a given maximal torus, {\it J. Algebra} {\bf 266} (2003), 87--101. 

\medskip \noindent
[SGA3] M. Demazure and A. Grothendieck, \textit{Sch\'{e}mas en Groupes (SGA 3)} tome III, 
Documents Math\'ematiques \textbf{{8}}, SMF, 2011 (second edition, revised and completed by Ph. Gille and P. Polo).

\medskip \noindent
[F 12] A. Fiori,  Special points on orthogonal symmetric spaces, {\it J. Algebra}, {\bf 372} (2012),
397-419.

\medskip \noindent
[G 12] S. Garibaldi, Outer automorphisms of algebraic groups and determining groups by their maximal tori, 
{\it Michigan Math. J.} {\bf 61} (2012), 227--237.

\medskip \noindent
[GR 12]  S. Garibaldi and A. Rapinchuk, Weakly commensurable S--arithmetic subgroups
in almost simple algebraic groups of types B and C, {\it Algebra Number Theory}
{\bf 7}  (2013), 1147-1178.

\medskip \noindent
[K 69] M. Kneser, {\it Lectures on Galois cohomology of classical groups}, Tata Institute of Fundamental Research
Lectures on Mathematics {\bf 47}, TIFR Bombay (1969). 

\medskip \noindent
[KMRT 98] M. Knus, A. Merkurjev, M. Rost and  J--P. Tignol, {\it The Book of Involutions}, AMS Colloquium Publications {\bf 44}, 1998.

\medskip \noindent
[Lee 14] T-Y. Lee, Embedding functors and their arithmetic properties, {\it Comment. Math. Helv.}, {\bf 89} (2014), 671--717.

\medskip \noindent
[LT 99] D. W. Lewis, J-P. Tignol, Classification theorems for central simple algebras with involution. With an appendix by R. Parimala. 
{\it Manuscripta Math.} {\bf 100} (1999), 259-276. 

\medskip \noindent
[MT 95] A. Merkurjev, J-P. Tignol, Multiples of similitudes and the Brauer group of homogeneous varieties, {\it J. reine angew. Math.} {\bf 461} (1995), 13-47.

\medskip \noindent
[NSW 08] J. Neukirch, A. Schmidt and K. Wingberg,
\textit{Cohomology of number fields}, Grundlehren der mathematischen
Wissenschaften \textbf{323}, Springer, 2008

\medskip \noindent
[PR 10] G. Prasad and A.S. Rapinchuk, Local--global principles for embedding of fields with
involution into simple algebras with involution, {\it Comment. Math. Helv.} {\bf 85} (2010),
583--645. 

\medskip \noindent 
[PT 04], R. Preeti and J-P. Tignol,  Multipliers of improper similitudes, {\it Doc. Math} {\bf 9} (2004), 183-204 

\medskip \noindent
[Sch 85] W. Scharlau, {\it Quadratic and hermitian forms}, Grundlehren der Mathematischen Wissenschaften {\bf 270}, Springer-Verlag, Berlin, 1985.

\medskip \noindent
[SGA 3] M. Demazure, A. Grothendieck, {\it Sch\'emas en groupes}, Documents Math\'ematiques {\bf 7, 8}, SMF, 2011.

\medskip \noindent
[T 61] T. Tsukamoto, On the local theory of quaternionic anti-hermitian forms, {\it J. Math. Soc. Japan} {\bf 13} (1961), 387Ð400.

\bigskip
\bigskip
Eva Bayer--Fluckiger and Ting-Yu Lee

EPFL-FSB-MATHGEOM-CSAG

Station 8

1015 Lausanne, Switzerland

\medskip

eva.bayer@epfl.ch

ting-yu.lee@epfl.ch

\bigskip

Raman Parimala

Department of Mathematics $ \&$ Computer Science

Emory University

Atlanta, GA 30322, USA.

\medskip
parimala@mathcs.emory.edu

 \end{document}